%% file: utor.tex
\documentclass[10pt]{amsart}

\input{preamble.tex}

\begin{document}

\author{David Mehrle}
	\address{Department of Mathematics, 
			University of Kentucky, 
			Lexington, Kentucky, U.S.A.}
	\email{davidm@uky.edu}
\author{J.D. Quigley}
	\address{Department of Mathematics, 
			University of Virginia, 
			Charlottesville, VA, U.S.A.}
	\email{mpb6pj@virginia.edu}
\author{Michael Stahlhauer}
	\address{Department of Mathematics, 
			Universit\"at Bonn, Bonn, Germany}
	\email{mstahlhauer@uni-bonn.de}
\title
	[Koszul Resolutions over Free Incomplete Tambara Functors]
	{Koszul resolutions over free incomplete Tambara functors 
		for cyclic $p$-groups}
\maketitle

\begin{abstract}
In equivariant algebra, Mackey functors replace abelian groups and incomplete Tambara functors replace commutative rings. 
In this context, we prove that equivariant Hochschild homology can sometimes be computed using Mackey functor-valued Tor. 
To compute these Tor Mackey functors for odd primes $p$, we define cyclic-$p$-group-equivariant analogues of the Koszul resolution which resolve the Burnside Mackey functor (the analogue of the integers) as a module over free incomplete Tambara functors (the analogue of polynomial rings). 
We apply these Koszul resolutions to compute Mackey functor-valued Hochschild homology of free incomplete Tambara functors for cyclic groups of odd prime order and for the cyclic group of order 9. 
\end{abstract}

\setcounter{tocdepth}{1}
\tableofcontents

\section{Introduction}

Hochschild homology is a homology theory for associative algebras that has found many applications, for example in algebraic topology, algebraic geometry, and number theory. 
One fundamental computation is the Hochschild homology of $k[x]$ as a $k$-algebra, which follows from the  isomorphism 
\[ 
	\HH_*(k[x]) \cong k[x] \otimes_k \Tor^{k[x]}_*(k,k). 
\]
There is a natural resolution of $k$ as a $k[x]$-module known as the Koszul complex. 

One might ask whether or not the same calculation can be made in equivariant algebra, where commutative rings are replaced by incomplete Tambara functors \cite{BH2018a} and abelian groups are replaced by Mackey functors. 
Equivariant versions of Hochschild homology have been studied extensively in recent years, see for instance \cite{AGH2021, AGHKK2021, BGHL2019, DMPR2017}, since they arise as algebraic approximations to equivariant and Real topological Hochschild homology. 
The equivariant case is complicated by the fact that the isomorphism above relies on the fact that $k[x]$ is flat as a $k$-module, but flatness is rare in the equivariant context, even for free algebras \cite{HMQ21}. 

In this paper, we calculate Mackey-functor-valued Hochschild homology of free incomplete Tambara functors for cyclic $p$-groups using Koszul resolutions. 
We find that these resolutions are much more complicated than their classical counterparts, even in the case where our free incomplete Tambara functors are flat, and doubly so when they are not. 
In contrast to \cite{BGHL2019}, where the authors study a twisted Hochschild homology of Green functors, we compute untwisted Hochschild homology defined via a cyclic bar construction. 

\subsection{Results} 
Throughout this paper, we work with Hochschild homology of incomplete Tambara functors. 
In \cref{section:HHisTOR}, we set up the definitions and prove that in certain cases, Hochschild homology of incomplete Tambara functors can be computed using Mackey-functor-valued $\Tor$, largely following the proof of the same fact in ordinary algebra \cite{CE56}. 

\begin{maintheorem}[\cref{corollary:HHTorIdenfitication}]
	\label{maintheoremA}
	If $\uR$ is a free incomplete Tambara functor whose underlying Mackey functor is projective, then there is an isomorphism of graded Mackey functors 
	\[ 
		\uHH_*(\uR) \cong \uR \boxtimes \uTor_*^{\uR}(\uA,\uA),
	\]
	where $\uA$ is the Burnside Mackey functor, which becomes an $\uR$-module by letting the incomplete Tambara functor generators act trivially. 
\end{maintheorem}

Our calculations make use of an explicit description of the free incomplete Tambara functors for cyclic $p$-groups via generators and relations. 
We record such a description for incomplete Tambara functors over cyclic $p$-groups in \cref{section:description_Tambara_functor}.

\begin{maintheorem}[]
	Let $\uR$ be a free incomplete Tambara functor for $C_{p^n}$. 
	If $\uR$ is flat as a Mackey functor, generators and relations for $\uR$ as an incomplete Tambara functor are given in \cref{thm:FreeIncompleteTambaraDescription}. 
\end{maintheorem}

We use the previous two theorems to compute the Hochschild homology of free $C_p$-Green functors and free $C_p$-Tambara functors which are flat as Mackey functors. 
For $G = C_p$, there are two free incomplete Tambara functors which are flat as Mackey functors: the free Green functor on a fixed generator, and the free Tambara functor on an underyling generator. 
In the case of a free Green functor on a fixed generator (for any group $G$), the Hochschild homology follows the same pattern as the Hochschild homology of $\bbZ[x]$ as a $\bbZ$-algebra. This follows from \cref{Prop:GreenFixed}. 
The case of a free $C_p$-Tambara functor on an underlying generator is the interesting one:

\begin{maintheorem}
	Let $\uR = \uA^{\cO^\top}[x_{C_p/e}]$ be a free $C_p$-Tambara functor on an underlying generator. Then 
	\[ 
		\uHH_*(\uR) \cong \uR \boxtimes \uTor^{\uR}_*(\uA, \uA),
	\]
	where $\uTor^{\uR}_*(\uA, \uA)$ is given in \cref{Prop:TambaraUnderlying}. 
\end{maintheorem}

To make this computation, we lift the Koszul resolution of $\z$ as a $\z[x_{1},\ldots,x_{p}]$-module to the world of Mackey functors to obtain an augmented complex of free $\uA^{\cO^{\top}}[x_{C_p/e}]$-modules. 
The underlying level will be the aforementioned Koszul resolution, but at the fixed level, there will be some nontrivial homology arising from the norm $\nm_e^{C_p}(x)$. 
We therefore add a second complex to kill this norm and obtain our Koszul resolution by taking a mapping cone. 

In addition to the above calculations in the $C_p$ case, we also make progress towards calculations over arbitrary cyclic $p$-groups for odd primes $p$. 
In general, we construct a free resolution of $\uA$ as a trivial $\uA^{\cO^\top}[x_{C_{p^n}/e}]$-module, where $\uA^{\cO^\top}[x_{C_{p^n}/e}]$ is the free Tambara functor on an underlying generator. 

\begin{maintheorem}[\cref{Thm:Koszulpn}]
	\label{maintheoremE}
	Let $p$ be an odd prime, $G = C_{p^n}$, and let $\uR = \uA^{\cO^{\top}}[x_{G/e}]$ be the free Tambara functor on an underlying generator. 
	The complex of $\uR$-modules $\bar{\uK}_\bullet$ defined in \cref{Constr:General} is a free $\uR$-module resolution of the Burnside Mackey functor $\bar{\uK}_\bullet \to \uA$, where the map $\bar{\uK}_0 \to \uA$ is the quotient map sending $x \mapsto 0$. The length of this resolution is $\sum_{i=0}^n p^i$. 
\end{maintheorem}

We apply this free resolution to compute $\uTor^{\uR}(\uA, \uA)$ for $G = C_9$ when $\uR$ is a free Tambara functor on an underlying generator. 

\begin{maincorollary}[\cref{C9CaseFinalCalculation}]
	Let $\uR = \uA^{\cO^\top}[x_{C_9/e}]$. 
	Then $\uTor_*^{\uR}(\uA, \uA)$ is given in \cref{C9CaseFinalCalculation}. 
\end{maincorollary}

Constructing our analogue of the Koszul complex over $C_{p^n}$ is involved, but essentially builds on the construction sketched for $C_p$ above. 
Roughly speaking, our resolution has to kill homology from polynomial generators, similar to the classical Koszul complex, but also has to kill homology from norms of generators at levels above the underlying level.
These norms must be killed by Koszul-type complexes, and so we obtain not just a single map between complexes, but an $(n+1)$-dimensional multicomplex of free $\uA^{\cO^{\top}}[x_{C_{p^n}/e}]$-modules. 
Our Koszul resolution is then obtained by taking the total complex. 

\subsection{Outline}

In \cref{section:background}, we review the definition of incomplete Tambara functors and their homological algebra. 
In \cref{section:description_Tambara_functor}, we give a generators-and-relations description of free incomplete Tambara functors for cyclic groups of prime power order. 
In \cref{section:HHisTOR}, we define Hochschild homology of Tambara functors and show that it can be computed using Mackey-functor-valued Tor, following the classical proof \cite{CE56}. 
We also discuss the relation of this Hochschild homology to Hill's genuine equivariant K\"ahler differentials \cite{Hil2017} in this section. 
In \cref{SS:KoszulCp}, we construct Koszul resolutions for free $C_p$-Tambara functors and use this to compute Hochschild homology in light of \cref{corollary:HHTorIdenfitication}. 
In \cref{section:c9}, we extend this to a Koszul resolution for a free Tambara functor on a fixed generator for $C_9$ and use this to compute Mackey-functor-valued Hochschild homology. 
Finally, in \cref{section:KoszulCpn}, we extend the techniques of the previous two sections to describe a Koszul resolution for free $C_{p^n}$-Tambara functors generated at the fixed level. 


\subsection{Acknowledgments}

The first and second authors were supported by NSF RTG grant DMS-2135884. 
The second author was supported by an AMS-Simons Travel Grant and NSF grants DMS-2039316 and  DMS-2414922 (formerly DMS-2203785 and DMS-2314082). 
The second and third authors were supported by the Max Planck Institute for Mathematics in Bonn. 

\section{Background}
\label{section:background}

\subsection{Mackey and Green Functors}

We assume the reader is familiar with Mackey functors, but we include a brief review to set notation. 
Fix a finite group $G$. 

Let $\cA^G$ be the Burnside category of $G$. 
Objects are finite $G$-sets, and morphisms are isomorphism classes of diagrams (spans) of the form 
\[ 
	A \from X \to B
\]
in finite $G$-sets. 
Composition of diagrams is given by pullback. 
The disjoint union of finite $G$-sets is the categorical product in $\cA^G$. 
The disjoint union also turns the hom-sets $\cA^G(X,Y)$ into commutative monoids. 

\begin{definition}
	A $G$-\emph{Mackey functor} is a product-preserving functor from the group completion of $\cA^G$ to abelian groups. 
	A morphism of Mackey functors is a natural transformation. 
	Denote by $\Mack(G)$ the category of $G$-Mackey functors. 
\end{definition}

Practically speaking, we can define a Mackey functor $\uM$ with a finite amount of data: 
\begin{itemize}
	\item for each transitive finite $G$-set $G/L$, an abelian group $\uM(G/L)$ with an action of the Weyl group\footnote{Recall that the Weyl group $W_G(L)$ is the group of automorphisms of $G/L$ as a finite $G$-set.} $W_G(L)$ by group homomorphisms;
	\item restrictions 
		\[
			\res_K^H := \uM([G/H \from G/K \xto{\id} G/K]),
		\] 
		which are equivariant for the actions of the Weyl groups in the sense that  
		\[
			g \cdot \res_K^H(y) = \res_K^H(y)
		\]
		for all $g \in W_H(K) \subset W_G(K)$ and all $y \in \uM(G/H)$;
	\item and transfers 
		\[
			\tr_K^H := \uM([G/K \xfrom{\id} G/K \to G/H]),
		\]
		which are equivariant for the actions of the Weyl groups in the sense that
		\[
			\tr_K^H(g \cdot x) = \tr_K^H(x)
		\]
		for all $g \in W_H(K) \subset W_G(K)$ and all $x \in \uM(G/K)$.
\end{itemize}
These data must satisfy some conditions, see \cite{Gre1971} for the original article or \cite{Web2000, Maz2013} for more recent surveys. 
Pay attention to the Weyl actions, because they will play an important role later.

For any Mackey functor $\uM$, we call the values of $\uM$ on transitive finite $G$-sets the \emph{levels of $\uM$}. 
Each level $\uM(G/H)$ has an action of the Weyl group $W_G(H)$, the automorphism group of $G/H$ as a finite $G$-set.
We display the levels as the lattice of subgroups of $G$, in a diagram with transfers and restrictions going between the levels:
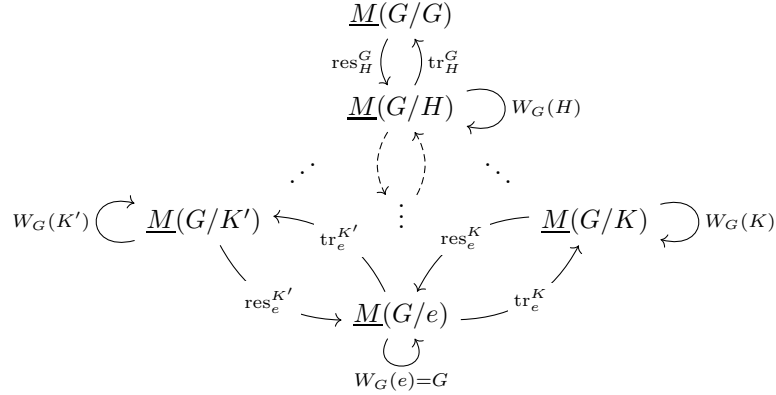
\[
	\begin{tikzcd} 
		& \uM(G/G)
			\ar[d, bend right=30, "\res^G_H"']
			\\
		& \uM(G/H) 
			\arrow[loop right, distance=2em, "W_G(H)"]
			\ar[u, bend right=30, "\tr_H^G"']
			\ar[d, dashed, bend right=30]
			\ar[dr, phantom, "\ddots" description, no head]
			\ar[dl, phantom, "{{\reflectbox{$\ddots$}}}" description, no head]
			\\
		\uM(G/K')
			\ar[dr, bend right=30, "\res_e^{K'}" description]
			\arrow[loop left, distance=2em, "W_G(K')"]
		& 
		\vdots
			\ar[u, dashed, bend right=30]
		& 
		\uM(G/K)
			\ar[dl, bend right=30, "\res_e^K" description]
			\arrow[loop right, distance=2em, "W_G(K)"]
		\\
		& 
		\uM(G/e)
			\ar[ur, bend right=30, "\tr_e^K" description]
			\ar[ul, bend right=30, "\tr_e^{K'}" description]
			\arrow[out=240,in=300,loop,swap,looseness=4,"W_G(e)=G"]	
	\end{tikzcd}
\]
Such a diagram is sometimes called a \emph{Lewis diagram}, because they were introduced in \cite{Lew1988}. 
Note that this is not a commutative diagram; neither restriction of a transfer nor transfer of a restriction is an identity. 
However, ignoring the arrows displaying the Weyl action, taking either the restrictions \emph{or} the transfers does yield a commutative diagram. 
We say that $\uM(G/G)$ is the \emph{top level} or \emph{fixed level} of the Mackey functor $\uM$, and $\uM(G/e)$ is the \emph{underlying level}. 

\begin{example}
	Let $C_n$ denote a cyclic group of order $n$. For a prime $p$, a $C_{p^2}$-Mackey functor $\uM$ has a diagram 
	\[ 
		\begin{tikzcd}[row sep=30]
			\uM(C_{p^2}/C_{p^2})
				\ar[d, bend right=30, "\res_{C_p}^{C_{p^2}}"']
			\\
			\uM(C_{p^2}/C_p)
				\ar[u, bend right=30, "\tr_{C_p}^{C_{p^2}}"']
				\ar[d, bend right=30, "\res_{e}^{C_p}"']
				\arrow[loop right, distance=2em, "C_{p^2}/C_p"]
			\\
			\uM(C_{p^2}/e)
				\ar[u, bend right=30, "\tr_{e}^{C_p}"']
				\arrow[loop right, distance=2em, "C_{p^2}"]			
		\end{tikzcd}
	\]
\end{example}

Because the category of Mackey functors is a diagram category in abelian groups, the category of Mackey functors is an abelian category. 
It is also a closed symmetric monoidal category, with monoidal product given by Day convolution of the tensor product in abelian groups with the Cartesian product of finite $G$-sets \cite{Lew1981}. 
We call this symmetric monoidal product the box product, and denote it by $\boxtimes$. 
The unit in this closed symmetric monoidal structure is the Burnside functor $\uA$.

\begin{example}
	The Burnside functor $\uA$ is the $G$-Mackey functor $\uA(G/H) = A(H)$, where $A(H)$ is the Burnside group of $H$. 
	Restriction from $H$ to $K$ is given by forgetting the $H$ action on a finite $H$-set $Y$ and only retaining the $K$-action.  
	Transfer is given by induction: the transfer from $K$ to $H$ of the class of a finite $K$-set $X$ is the class of $H \times_K X$. 
\end{example}

Because $\uA$ is the unit for the box product, it is a commutative monoid -- a commutative ring-like object in Mackey functors.  
We call such a Mackey functor a Green functor. 

\begin{definition}
A \emph{Green functor} is a commutative monoid for the box product in the category of Mackey functors. 
\end{definition}

Like Mackey functors, a Green functor $\uR$ is determined by a finite amount of data: 
\begin{itemize}
	\item for each transitive finite $G$-set $G/L$, a commutative ring $\uR(G/L)$ with an action of $W_G(L)$ via ring homomorphisms, 
	\item restrictions 
	\[
		\res^H_K \colon \uR(G/H) \to \uR(G/K),
	\]
	which are homomorphisms of commutative rings,
	\item transfers  
	\[
		\tr_K^H \colon \uR(G/K) \to \uR(G/H),
	\]
	which are homomorphisms of $\uR(G/H)$-modules, where $\uR(G/K)$ becomes a $\uR(G/H)$-module via restriction: 
	\[ 
		y \tr_K^H(x) = \tr_K^H(\res_K^H(y) x)
	\]
	for all $y \in \uR(G/H)$ and $x \in \uR(G/K)$. This is often called \emph{Frobenius reciprocity}. 
\end{itemize}
These data must satisfy the same conditions as Mackey functors. 
See \cite{Dress1971,Lew1981,Maz2013} for more on Green functors.

\subsection{Tambara Functors}

Green functors are only one kind of ring-like object in the category of Mackey functors. 
Tambara functors are Green functors with the extra data of ``multiplicative transfers" or norms.

\begin{definition}
	Denote by $\cP^G$ the \emph{category of polynomials in finite $G$-sets}. 
	Objects in this category are finite $G$-sets. 
	Morphisms in this category are \emph{polynomials}, that is, equivalence classes of diagrams 
	\[ 
		A \from X \to Y \to B,
	\]
	where two such diagrams are considered equivalent if there is a commutative diagram in finite $G$-sets of the form 
	\[ 
		\begin{tikzcd}[row sep=0]
			& 
			X \ar[r] \ar[dd, "\cong"] \ar[dl] &
			Y \ar[dr] \ar[dd, "\cong"] \\
			A & 
			& 
			& 
			B \\
			& 
			X' \ar[r] \ar[ul] &
			Y' \ar[ur] & 
		\end{tikzcd}
	\]
	Composition in this category is given by \cite[Proposition 7.1]{Tam1993}. 
\end{definition}

The composition in this category is tricky to define, so in practice one works with this category using a generating set of morphisms and relations. 

\begin{definition}
	Let $f \colon A \to B$ be a morphism of finite $G$-sets. Define morphisms in $\cP^G$
	\begin{align*}
		R_f &:= [B \xfrom f A \xto{\id} A \xto{\id} A ]\\
		N_f &:= [A \xfrom\id A \xto{f} B \xto{\id} B ]\\
		T_f &:= [A \xfrom\id A \xto\id A \xto f B]
	\end{align*}
\end{definition}

\begin{theorem}[{cf. \cite[Section 7]{Tam1993}}]\
	\label{theorem:PolynomialCategoryGens}
	\begin{enumerate}[(a)]
		\item Any morphism in $\cP^G$ can be written as a composite 
		\[ 
			T_f \circ N_g \circ R_h = [A \xfrom h X \xto g Y \xto f B]
		\]
		\item $R$, $N$, and $T$ define functors from finite $G$-sets to $\cP^G$; $R$ is contravariant, $N$ and $T$ are covariant. 
		\item Given a pullback of finite $G$-sets, 
		\[
			\begin{tikzcd}
				A' \ar[r, "g'"] \ar[d, "f'"] & 
				A \ar[d, "f"] \\
				B' \ar[r, "g"] & 
				B,
			\end{tikzcd} 
		\]
		we have 
		\begin{align*}
			N_{g'} \circ R_{f'} &= R_{f} \circ N_g \\
			T_{g'} \circ R_{f'} &= R_{f} \circ T_g 
		\end{align*}
		\item Given any diagram isomorphic to one of the form (an \emph{exponential diagram})
		\[ 
			\begin{tikzcd}
				X \ar[d, "f'"] & 
				A \ar[l, "g'"] &
				X \times_Y \prod_f A \ar[l, "h"] \ar[d, "g"] \\ 
				Y
				& 
				& 
				\prod_f A \ar[ll, "f"]
			\end{tikzcd}
		\]
		we have 
		\[ 
			T_f \circ N_g \circ R_h = N_{f'} \circ T_{g'}
		\]
	\end{enumerate}
\end{theorem}

In practice, one works with morphisms in $\cP^G$ by writing them in the form $T_f \circ N_g \circ R_h$, and using parts (c) and (d) of the theorem to commute the three types of generating morphisms so that they are in this form. 

The category $\cP^G$ has a product given by the disjoint union of finite $G$-sets. 
The disjoint union of finite $G$-sets also turns the hom-sets $\cP^G(X,Y)$ into commutative monoids. 

\begin{definition}[{cf. \cite[Theorem 6.2]{Tam1993}, \cite[Proposition 4.2]{BH2018a}}]
	A \emph{Tambara functor} $\uT$ is a product-preserving functor from the group completion of $\cP^G$ to sets such that each $\uT(X)$ is an abelian group. 
	A morphism of Tambara functors is a natural transformation.
\end{definition}

A Tambara functor is a Green functor with extra structure, cf.~\cite[Theorem 4.13]{BH2018a}. 
In addition to the data of a Green functor, a Tambara functor $\uT$ has \emph{norms} 
\[ 
	\nm_K^H \colon \uT(G/K) \to \uT(G/H) 
\]
for all $K \leq H \leq G$, which are homomorphisms of the multiplicative monoids of the commutative rings $\uT(G/K)$ and $\uT(G/H)$. 
The norms must satisfy certain conditions for the norm of a sum and norm of a transfer, which can be deduced from \cref{theorem:PolynomialCategoryGens} (this is not easy; the general formula is due to Mazur, \cite{HM2019}). 

\begin{example}
	The Burnside functor $\uA$ is a Tambara functor. 
	The norm from $K$ to $H$ is given by coinduction of finite $K$-sets. 
	Explicitly, if $X$ is a finite $K$-set, then 
	\[ 
		\nm_K^H(X) = [\Set^K(H,X)],
	\]
	where $H$ acts on the domains of functions. 
\end{example}

\subsection{Incomplete Tambara Functors}

By taking subcategories of the category of polynomials $\cP^G$, we can produce variations on Tambara functors.

\begin{definition}
	Let $\cD$ be a wide (i.e. contains all objects), pullback stable \cite[Definition 2.8]{BH2018a} subcategory of the category of finite $G$-sets. 
	Let $\cP^G_{\cD}$ denote the wide subcategory of $\cP^G$ whose morphisms are those polynomials 
	\[ 
		A \xfrom f X \xto g Y \xto h B
	\]
	such that $g$ is a morphism in $\cD$. 
	This is the \emph{category of polynomials with exponents in $\cD$.}
\end{definition}

This is a subcategory of $\cP^G$ by {\cite[Theorem 2.10]{BH2018a}}.

\begin{example}
	The wide subcategory of finite $G$-sets whose only morphisms are isomorphisms is a pullback stable subcategory. 
	The category of polynomials with exponents in this subcategory, $\cP^G_{\Iso}$, is equivalent to the Burnside category $\cA^G$.
	Therefore, product-preserving functors from $\cP^G_{\Iso}$ to abelian groups are Mackey functors; the forgetful functor given by pullback along the inclusion $\cP^G_{\Iso} \subset \cP^G$ shows that any Tambara functor has an underlying Mackey functor. 
\end{example}

\begin{example}
	Other choices for $\cD$ give Green functors (\cite[Section 4]{BH2018a}) or even non-unital Tambara functors (\cite[Definition 2.17]{Hil2017}). 
	The forgetful functor given by pullback along the inclusion $\cP^G_{\cD} \subseteq \cP^G$ shows that every Tambara functor has an underlying Green functor. 
\end{example}

For certain subcategories, we get objects that still form commutative-ring-like structure in Mackey functors. 
These are called the incomplete Tambara functors. 

\begin{definition}
	Let $\cO$ be a wide, pullback stable, finite coproduct complete subcategory of finite $G$-sets. 
	We call $\cO$ an \emph{indexing category for $G$}. 
\end{definition}

For a fixed finite group $G$, indexing categories form a lattice under inclusion \cite[Proposition 2.16]{Rub2021}.

\begin{notation}
	We denote the least element of the lattice of indexing categories by $\cO^\bot$ and the greatest element by $\cO^\top$. 
\end{notation}

\begin{definition}
	An \emph{$\cO$-Tambara functor} $\uR$ is a product-preserving functor from the group completion of $\cP^G_\cO$ to sets such that each $\uR(X)$ is an abelian group. 
	A morphism of $\cO$-Tambara functors is a natural transformation. An \emph{incomplete Tambara functor} is an $\cO$-Tambara functor for some $\cO$. 
\end{definition}

\begin{example}
	$\cO^\top$ is the whole category of finite $G$-sets, so $\cO^\top$-Tambara functors are (complete) Tambara functors. 

	$\cO^\bot$ is the subcategory of finite $G$-sets with morphisms $f \colon X \to Y$ that preserve isotropy, i.e. the stabilizer subgroup of $f(x)$ is the same as the stabilizer subgroup of $x$. 
	$\cO^\bot$-Tambara functors are Green functors \cite[Section 4]{BH2018a}.
\end{example}

Any indexing category $\cO$ contains $\cO^\bot$, so any incomplete Tambara functor has an underlying Green functor, with the forgetful functor given by pullback along the inclusion $\cP_{\cO^\bot}^G \subset \cP^G_{\cO}$. 
In a (complete) Tambara functor, this forgetful functor amounts to forgetting the norms. 

An incomplete Tambara functor is called incomplete because it has some, but not necessarily all, of the norms in a Tambara functor \cite[Theorem 4.13]{BH2018a}. 
In general, an $\cO$-Tambara functor $\uT$ is a Green functor together with norms 
\[ 
	\nm_K^H \colon \uT(G/K) \to \uT(G/H) 
\]	
whenever $G/K \to G/H$ is a morphism in $\cO$. 
These norms must satisfy formulas for the norm of a transfer and norm of a sum, as with Tambara functors. 

\begin{notation}
	If $\uR$ is an $\cO$-Tambara functor, we write $\uR^{\cO}$ to record the indexing category for $\uR$. 
	If the indexing category is clear from context, we may omit it. 
\end{notation}

An example of an incomplete Tambara functor which is not a Tambara functor is the \emph{free $\cO$-Tambara functor} 
\[
	\uA^{\cO}[x_{G/H}] := \cP^G_{\cO}(G/H,-)
\]
where $\cO \neq \cO^\top$. 
We consider these objects in detail in \cref{section:description_Tambara_functor}. 

\subsection{Modules}

A Green functor (aka $\cO^\bot$-Tambara functor) is a commutative monoid for the box product on Mackey functors, and therefore has a well-defined notion of modules. 
For any other kind of incomplete Tambara functor, we adopt this notion of module. 

\begin{definition}
	If $\uR$ is an incomplete Tambara functor, an $\uR$-module is a module over the underlying Green functor of $\uR$. 
\end{definition}

Modules also have a nice levelwise description. 
An $\uR$-module $\uM$ is a Mackey functor such that:
\begin{itemize}
	\item each $\uM(G/H)$ is a $\uR(G/H)$-module;
	\item restrictions $\res^H_K \colon \uM(G/H) \to \uM(G/K)$ are $\uR(G/H)$-module homomorphisms, where $\uM(G/K)$ becomes a $\uR(G/H)$-module via restriction: 
	\[ 
		\res^H_K(rm) = \res^H_K(r) \res^H_K(m)
	\]
	for $r \in \uR(G/H)$ and $m \in \uM(G/H)$;
	\item transfers $\tr_K^H \colon \uM(G/K) \to \uM(G/H)$ are $\uR(G/H)$-module homomorphisms, where $\uM(G/K)$ becomes a $\uR(G/H)$-module via restriction: 
	\[
		r \cdot \tr_K^H(n) = \tr_K^H(\res^H_K(r)n)
	\]
	for $r \in \uR(G/H)$ and $n \in \uM(G/K)$;
	\item for all $s \in \uR(G/K)$ and $m \in \uM(G/H)$, 
	\[ 
		\tr_K^H(s)m = \tr_K^H(s\res_K^H(m)).
	\]
\end{itemize}

By work of Lewis \cite{Lew1981}, we know that the category of $\uR$-modules is an abelian, closed symmetric monoidal category with symmetric monoidal product 
\[ 
	\uM \boxtimes_{\uR} \uN := \coeq(\uM \boxtimes \uR \boxtimes \uN \rightrightarrows \uM \boxtimes \uN).
\]
Therefore, we can do homological algebra in this category. 

\begin{definition}
	Let $U$ be a finite $G$-set. 
	The \emph{free Mackey functor on a generator $x_U$ at level $U$} is 
	\[ 
		\uA\{x_U\} := \cA^G(U, -)
	\]
	The \emph{free $\uR$-module on a generator at level $U$} is 
	\[ 
		\uR\{x_U\} := \uR \boxtimes \uA\{x_U\}. 
	\]
	A \emph{free $\uR$-module} is any $\uR$-module $\uM$ such that 
	\[ 
		\uM \cong\bigoplus_i \uR\{x_{U_i}\} 
	\]
	for some finite $G$-sets $U_i$ (possibly infinitely many). 
\end{definition}

Free modules deserve to be called free because they admit the universal property 
\[ 
	\uR\mhyphen\Mod(\uR\{x_U\}, \uM) \cong \uM(U). 
\]

\begin{example}
	There is a canonical isomorphism of $\uR$-modules $\uR\{x_{G/G}\} \cong \uR$ given by $x_{G/G} \mapsto 1 \in \uR(G/G)$.
\end{example}

In practice, understanding a free $\uR$-module generated at level $U$ is done by breaking up $U$ into orbits $U = \bigsqcup_i G/H_i$, and then 
\[ 
	\uR\{x_U\} \cong \bigoplus_i \uR\{x_{G/H_i}\}. 
\]
To understand $\uR\{x_{G/H_i}\}$, it helps to remember that $\uR\{x_U\}(V) = \uR(U \times V)$. 
Informally, we are adding a module generator $x$ at level $G/H_i$ on which the Weyl group $W_G(H_i)$ acts freely, and then freely adding in all transfers and restrictions of this element at other levels, subject to the double coset formula for the restriction of a transfer. 
(We don't freely add norms because an $\uR$-module is a module over the underlying Green functor of $\uR$, which has no norms!)

\begin{example}
	If $\uR$ is an incomplete $C_3$-Tambara functor, then $\uR\{x_{C_3/e}\}$ is described by 
	\[ 
		\uR\{z_{C_3/e}\} = 
		\begin{tikzcd}
			\uR(C_3/e)\{\tr(z)\}
				\ar[d, bend right=50, "\res"{left}]
			\\
			\uR(C_3/e)\{z^{(0)},z^{(1)},z^{(2)}\}
				\ar[u, bend right=50, "\tr"{right}]	
				\arrow[out=240,in=300,loop,swap,looseness=4, "C_3"]
		\end{tikzcd}
	\]
	where the Weyl action of $W_{C_3}(e) = C_3$ on the underlying level is given by permuting the three variables, $\gamma \cdot z^{(i)} = z^{(i+1)}$, with indices taken mod 3. 
	The transfer is the $\uR(C_3/e)$-linear homomorphism determined by 
	\[
		\tr(z^{(i)}) = \tr(z); 
	\]	
	we think of $\tr(z)$ on the top level as a formal element that generates that level as an $\uR(C_3/e)$-module. In fact, any element $f \tr(z)$ at the top level is a transfer of $f z^{(0)}$ (or $fz^{(1)}$ or $fz^{(2)}$). 
	The restriction is determined by the rule that restriction of a transfer is the sum over the Weyl conjugates of an element. For example,  
	\[ 
		\res(\tr(z)) = z^{(0)} + z^{(1)} + z^{(2)}. 
	\]
	The top level of $\uR\{z_{C_3/e}\}$ becomes a module over $\uR(C_3/C_3)$ via restriction, while the underlying level is a free $\uR(C_3/e)$-module of rank 3 with basis $\{z^{(0)},z^{(1)},z^{(2)}\}$.
\end{example}

The following lemmas illustrate how to work with homomorphisms between free modules. 
They are also crucial ingredients in the computations in \cref{SS:KoszulCp,section:c9,section:KoszulCpn}. 

\begin{lemma}
	\label{lemma:freeModules}
Let $\uR$ be an incomplete Tambara functor for $G = C_{p^n}$. The free $\uR$-module $\uR\{x_{G/H}\}$ has levels as follows. 
	\begin{itemize}
		\item At level $G/H$, 
			\[
				\uR\{x_{G/H}\}(G/H) \cong 
				\bigoplus_{gH \in G/H} \uR(G/H) \cdot gx
			\]
			with $G/H$-action permuting the summands. The $\uR(G/H)$-module structure is induced by the action of $\uR(G/H)$ on itself. 
		\item for $L \leq H$, 
			\[
				\uR\{x_{G/H}\}(G/L) \cong 
				\bigoplus_{gH \in G/H} 
					\uR(G/L) \cdot \res^H_L(gx),
			\]
			with $G/L$-action by $aL \cdot \res^H_L(gx) = \res^H_L((ag)x)$
			Here, $\res^H_L(gx)$ is a formal restriction of a Weyl conjugate of the generator. 
			The $\uR(G/L)$-module structure is induced by the action of $\uR(G/L)$ on itself. 
		\item For $L \geq H$, 
			\[
				\uR\{x_{G/H}\}(G/L) \cong 
				\bigoplus_{gL \in G/L} 
					\uR(G/H) \cdot \tr_H^L(gx),
			\]
			with $G/L$-action permuting the summands. 
			Here, $\tr_H^L(gx)$ is a formal transfer of an element $gx$ at level $G/H$. 
			This becomes an $\uR(G/L)$-module via restriction $\res^L_H$. 
	\end{itemize}
	The transfers and restrictions in this $\uR$-module are given as follows. 
	\begin{itemize}
		\item For subgroups $K \leq L \leq H$, transfer and restriction in $\uR\{x_{G/H}\}$ are given by applying the transfer and restriction of $\uR$ on each summand. 
		\item For subgroups $L \geq K \geq H$, transfer in this free module is a sum over $L/K$-cosets 
		\[
		\begin{tikzcd}[row sep=small]
			\uR\{x_{G/H}\}(G/K) 
				\ar[r, "\tr_K^L"] 
				\ar[d, "\rotatebox{90}{$\cong$}", phantom] 
				& 
			\uR\{x_{G/H}\}(G/L) 
				\ar[d, "\rotatebox{90}{$\cong$}", phantom]
				\\
			\displaystyle
			\bigoplus_{aK \in G/K} \uR(G/H) \cdot \tr_H^K(ax)
				\ar[r, "\tr_K^L"] 
				& 
			\displaystyle
			\bigoplus_{gL \in G/L} \uR(G/H) \cdot \tr_H^L(gx)
				\\
			\displaystyle
			\sum_{aK \in G/K} t_a \tr_H^K(ax) 
				\ar[r, mapsto] 
				\ar[u, "\rotatebox{90}{$\in$}", phantom]
				& 
			\displaystyle
			\sum_{gL \in G/L} 
			\bigg(
				\sum_{\substack{aK \in G/K \\ aK \subseteq gL}} 
					t_a
			\bigg)
			\tr_H^L(gx)
				\ar[u, "\rotatebox{90}{$\in$}", phantom]
		\end{tikzcd}
		\] 
		and restriction is a diagonal, expanding the coefficient of $\tr_H^L(gx)$ to cover all $K$-cosets inside $gL$:
		\[
		\begin{tikzcd}[row sep=small]
			\uR\{x_{G/H}\}(G/K) 
				\ar[d, "\rotatebox{90}{$\cong$}", phantom] 
				& 
			\uR\{x_{G/H}\}(G/L) 
				\ar[d, "\rotatebox{90}{$\cong$}", phantom]
				\ar[l, "\res^L_K"]
				\\
			\displaystyle
			\bigoplus_{aK \in G/K} \uR(G/H) \cdot \tr_H^K(ax)
				& 
			\displaystyle
			\bigoplus_{gL \in G/L} \uR(G/H) \cdot \tr_H^L(gx)
				\ar[l, "\res^L_K"]
				\\
			\displaystyle
			\sum_{gL \in G/L} 
			\bigg(
				\sum_{\substack{aK \in G/K\\aK \subseteq gL}} 
					t_g \tr_H^K(ax) 
			\bigg)
				\ar[u, "\rotatebox{90}{$\in$}", phantom]
				& 
			\displaystyle
			\sum_{gL \in G/L} t_g \tr_H^L(gx).
				\ar[u, "\rotatebox{90}{$\in$}", phantom]
				\ar[l, mapsto]
		\end{tikzcd}
		\] 
		\item If $K < H < L$, then $\tr_K^L$ is the composite of $\tr_K^H$ as in the first case with $\tr_H^L$ as in the second case.
	\end{itemize}
\end{lemma}

\begin{example}
	Let $\gamma$ be the generator of $C_9$. 
	If $\uR$ is a $C_9$-Tambara functor, then $\uR\{x_{C_9/C_3}\}$ is depicted in the diagram below. 
	\[
	\begin{tikzcd}[row sep=huge, ampersand replacement=\&]
		\uR(C_9/C_3)\{ \tr_{C_3}^{C_9}(x) \}
			\ar[d, bend right=30,
				"\scalebox{0.75}{$\begin{bmatrix}
					\scalebox{1.33}{1} \\ 
					\scalebox{1.33}{1} \\ 
					\scalebox{1.33}{1}
				\end{bmatrix}$}"']
			\\
		\uR(C_9/C_3)\{x, \gamma x, \gamma^2 x\}
			\ar[u, bend right=30, 
				"{\big[1\ 1\ 1\big]}"']
			\ar[d, bend right=30, 
				"\scalebox{0.5}{$\begin{bmatrix} 
					\scalebox{2}{$\res_e^{C_3}$} & & \\ 
					& \scalebox{2}{$\res_e^{C_3}$} & \\
					& & \scalebox{2}{$\res_e^{C_3}$} 
				\end{bmatrix}$}"']
			\\
		\uR(C_9/e)\{	
			\res^{C_3}_e(x), 
			\res^{C_3}_e(\gamma x), 
			\res^{C_3}_e(\gamma^2 x)
			\}
			\ar[u, bend right=30, 
				"\scalebox{0.5}{$\begin{bmatrix} 
					\scalebox{2}{$\tr_e^{C_3}$} & & \\ 
					& \scalebox{2}{$\tr_e^{C_3}$} & \\
					& & \scalebox{2}{$\tr_e^{C_3}$} 
				\end{bmatrix}$}"']
	\end{tikzcd}
	\]
\end{example}

\begin{lemma}
	\label{lemma:kernels}
	Let $\uR$ be an incomplete Tambara functor for $G = C_{p^n}$. Let $K < H$ be subgroups of $G$. 
	Let $f$ be a $\uR$-module homomorphism 
	\[
		f \colon \uR\{x_{G/H}\} \to \uR\{y_{G/K}\}
	\]
	such that $f(x_{G/H})$ and its Weyl conjugates under the $G/H$-action on $\uR\{y_{G/K}\}(G/H)$ are $\uR$-linearly independent, in the sense that if for any $t_g \in \uR(G/K)$ and $H_g \leq H$ that satisfy an equation of the form
	\[
		\sum_{gH \in G/H} \tr_{H_g}^K(t_g \res^H_{H_g}(f(gx))) = 0, 
	\]
	we must have that $t_g = 0$ for all $g$.
	
	Then, the kernel of $f$ is the sub-$\uR$-module of $\uR\{x_{G/H}\}$ generated by the kernels of restrictions $\res^{(-)}_K$ to $K$ in $\uT$.
	In particular, the kernel of $f$ vanishes at level $G/K$ and below.
\end{lemma}

\begin{proof} 
	For $L \leq H$, an element in $\uR\{x_{G/H}\}(G/L)$ has the form
	\[
		\sum_{gH \in G/H} t_g\,\res^H_L(gx),
	\]
	with $t_g \in \uR(G/L)$. Assume this element lives in the kernel of $f$. 
	
	\begin{itemize}
		\item If $L = H$, then 
		\[
			f\left(\sum_{gH \in G/H} t_g\,gx\right) 
				= \sum_{gH \in G/H} \res^H_K(t_g)\,gf(x) 
				= 0,
		\]	
		where $\res^H_K(t_g)$ is restriction in $\uR$. 
		By linear independence, we have $\res^H_K(t_g) = 0$ for all $g$, so this element lies in the kernel of restriction. 

		\item If $K < L \leq H$, then 
		\[
			f\left(\sum_{gH \in G/H} t_g\,\res^H_L(g x) \right)
				= \sum_{gH \in G/H} \res^L_K(t_g) \res^H_L(g f(x))	
				= 0.
		\]
		By linear independence, $\res^L_K(t_g) = 0$ for all $g$. Hence, $t_g \in \ker(\res^L_K)$ for all $g$.		
		\item If $L \leq K$, then 
		\[
			f\left(\sum_{gH \in G/H} t_g\,\res^H_L(gx)\right)
				= \sum_{gH \in G/H} t_g \res^H_L(gf(x)) 
				= 0
		\]
		By linear independence, we have $t_g = 0$ for all $g$. Hence, the element must have been zero to begin with.
	\end{itemize}
	
	If $L > H$, an element of $\uR\{x_{G/H}\}(G/L)$ has the form
	\[
		\sum_{gL \in G/L} t_g \tr_H^L(gx)
	\]
	with $t_g \in \uR(G/H)$. If this element is in the kernel of $f$, then 
	\[
		 	f\left(\sum_{gL \in G/L} t_g \tr_H^L(gx)\right) 
		 		= \sum_{gL \in G/L} \tr_H^L(\res^H_K(t_g) f(gx))
		 		= 0. 
		\]
		Thus, again by linear independence, we have $\res^H_K(t_g) = 0$ for all $g$.
	
	Finally, it is clear from the formulas above that the submodule generated by the kernels of restrictions $\res^{(-)}_K$ is contained in the kernel of $f$. 
\end{proof}

The important conclusion of this lemma is that $i_K^*\ker(f) = 0$, where $i_K^*$ denotes the restriction of a $G$-Mackey functor to a $K$-Mackey functor. 
For an example of this lemma in action, see \cref{spiderExampleCohomology}.

\begin{remark}
Note that the notion of linear independence we used in this lemma is different from the one in \cite[Definition 2.4]{Lee2019}, since the choice of coefficients $t_g$ is from $\uR(G/K)$ instead of $\uR(G/H_g)$. This choice makes sense in our case since the target $\uR$-module is generated at level $G/K$, so all of its levels are $\uR(G/K)$-modules in a natural way. With the notion as presented in \cite{Lee2019}, the map defined by $f(x_{G/H}) = \tr_K^H(y_{G/K})$ would not be covered by the assumptions of the lemma, since the element $\tr_K^H(y_{G/K})$ is annihilated by any element in $\ker(\res_K^H)$ (which is the main point of the lemma). Since the examples of maps $f$ we consider here are mostly of this form, we use this different version of linear independence.
\end{remark}

\subsection{Homological Algebra}

\begin{definition}
	An $\uR$-module is \emph{projective} if it is a summand of a free $\uR$-module. 
	An $\uR$-module $\uM$ is \emph{flat} if $\uM \boxtimes_{\uR} -$ is an exact functor. 
\end{definition}

As is the case in ordinary homological algebra, free implies projective and projective implies flat for $\uR$-modules. 
The concept of flatness in the case of modules over an incomplete Tambara functor is much more subtle than it is for modules over a commutative ring. 
For example, for solvable groups, the constant functor $\ubbZ$ is not flat \cite[Lemma 3.2.16]{HMQ21}, nor is its dual $\ubbZ^*$. 
In related categories, projective may not even imply flat \cite{Lew1999}. 

To determine whether or not an $\uR$-module is flat, we make use of the $\uTor$ functors. 
These are $\uR$-module-valued $\Tor$ functors internal to the abelian category of $\uR$-modules. 

\begin{definition}
	Let $\uM$ be an $\uR$-module for an incomplete Tambara functor $\uR$. 
	Let 
	\[ 
		\uTor^{\uR}_n(\uM,-)
	\]
	be the $n$-th left-derived functor of $\uM \boxtimes -$. 
	For any other $\uR$-module $\uN$, $\uTor^{\uR}_n(\uM, \uN)$ is an $\uR$-module. 
\end{definition}

The properties of the functors $\uTor$, and the methods of their calculation, are the same as those for $\Tor$ in the case of modules over a ring. 
This follows formally.

\section{Free incomplete Tambara functors over cyclic groups of prime power order}
	\label{section:description_Tambara_functor}

Now that the basic properties of Mackey and incomplete Tambara functors have been recalled, we will study our main players: free incomplete Tambara functors over cyclic groups of prime power order. 
Our main result in this section (\cref{thm:FreeIncompleteTambaraDescription}) is an explicit description of a special subset of these, namely, those free incomplete Tambara functors which are also free as Mackey functors. 

\begin{definition}
	The \emph{free $\cO$-Tambara functor on a generator at level $G/H$} is the $\cO$-Tambara functor
	\[
		\uA^{\cO}[x_{G/H}] = \cP^G_{\cO}(G/H,-). 
	\]
\end{definition}

These functors satisfy the universal property that $\cO$-Tambara morphisms out of $\uA^{\cO}[x_{G/H}]$ represents evaluation at $G/H$:
\[ 
	\cO\mhyphen\Tamb_G(\uA^{\cO}[x_{G/H}], \uR) \cong \uR(G/H).
\]

\begin{remark}
	The right-hand-side of this isomorphism is a commutative ring for any $\uR$, demonstrating that $\uA^{\cO}[x_{G/H}]$ has a canonical co-ring structure. 
	In fact, \cite[Section 2]{BH2019} shows that as the finite $G$-set varies, $\uA^{\cO}[x_{(-)}]$ is a co-$\cO$-Tambara functor object in $\cO$-Tambara functors. 
\end{remark}

	Explicit descriptions of the four free incomplete $C_2$-Tambara functors $\uA^{\cO^\bot}[x_{C_2/C_2}]$, $\uA^{\cO^\bot}[x_{C_2/e}]$, $\uA^{\cO^\top}[x_{C_2/C_2}]$, and $\uA^{\cO^\top}[x_{C_2/e}]$ are given in \cite[Section 3]{BH2019}. 
	In this section, we obtain an explicit levelwise description of $\uA^{\cO}[x_{G/H}]$ -- assuming that this free Tambara functor is free as a Mackey functor, cf. \cite{HMQ21}. 
	
Let $G = C_{p^n}$, $H = C_{p^m}$ with $0 \leq m \leq n$, and $\cO$ a $C_{p^n}$-indexing category such that $\uA^{\cO}[x_{G/H}]$ is free as a Mackey functor. 
	Let $\uR := \uA^{\cO}[x_{G/H}]$. Our goal in this section is to explicitly describe the $\cO$-Tambara functor structure of $\uR$.

\begin{remark}\label{remark:indexing_category}
	In \cite[Corollary C]{HMQ21}, the indexing categories $\cO$ such that the free incomplete Tambara functor $\uA^{\cO}[x_{G/H}]$ is free as a Mackey functor are classified. 
	In the case $G= C_{p^n}$, the indexing category needs to be trivial below the subgroup $H$ and it has to contain the norm $\nm_H^G$. 
	These requirements need not determine the indexing category fully. 
	However, the free incomplete Tambara functor $\uA^{\cO}[x_{G/H}]$ generated by an element at level $H$ is completely determined, and $\cO$ may be chosen to contain all possible norms above level $H$. 
	This follows from the observation that the Mackey decomposition formula for the restriction of the total norm $\nm_H^G$ to a subgroup $K\geq H$ forces the intermediate norm $\nm_H^K$ to exist
		\footnote{The Mackey decomposition formula is usually stated for restrictions of transfers, but the evident analog for restrictions of norms holds, for instance, by observing that the formula for commuting restrictions and transfers is the same as for commuting restrictions and norms, cf. \cite[Thm. 2.13(c)]{HMQ21}.}. 
	From this, any norm $\nm_K^{K^\prime}$ for $H\leq K\leq K^{\prime}\leq G$ is determined by using the relations between transfer, norm and restriction by using that $x_{G/H}$ lives in level $H$.
\end{remark}

\begin{proposition}
	The $\cO$-Tambara functor $\uR = \uA^{\cO}[x_{C_{p^n}/C_{p^m}}]$ may be decomposed as
	\[ 
		\uR \cong 
			\bigoplus_{k = m}^{n} 
			\bigoplus_{[\vec{v}]} 			
			\uA\{ \nm_{C_{p^m}}^{C_{p^k}}(x^{\vec{v}} ) \}. 
	\]
	In this sum, we index over equivalence classes of tuples $\vec{v}$ of non-negative integers indexed by cosets in $C_{p^n}/C_{p^k} \cong C_{p^{n-k}}$, with the equivalence relation given by action of the Weyl group $W_{C_{p^n}}(C_{p^k})$ by shifting the entries, minus the set $\bar{\Delta}$ of elements $\vec{v}$ such that there is $k^\prime > k$ and $\vec{w}\in \bbZ_{\geq 0}^{|C_{p^n}/C_{p^{k^\prime}}|}$ such that $v_{\gamma} = w_{\pi(\gamma)}$ for all $\gamma \in C_{p^n}/C_{p^k}$, where $\pi\colon C_{p^n}/C_{p^k}\to C_{p^n}/C_{p^{k^\prime}}$ is the projection.
\end{proposition}


\begin{proof}
	This follows from the explicit identification of $\uR$ as a free Mackey functor given in the proof of \cite[Theorem A]{HMQ21}. 
	The argument given there identifies $\uR = \uA^{\cO}[x_{C_{p^n}/C_{p^m}}]$ as the norm $N_{C_{p^m}}^{C_{p^n}}(\bbZ[x] \otimes \uA)$, where $M \otimes \uA$ is the Green functor with 
\[
	(M \otimes \uA)(T) = M \otimes \uA(T).
\]
	This is the norm of a free Mackey functor, which itself is free by \cite[Proposition 4.5]{HMQ21}. 
	The decomposition into single generators is then obtained using
	\begin{align*}
	N_{C_{p^m}}^{C_{p^n}}(\bbZ[x] \otimes \uA) & \cong N_{C_{p^m}}^{C_{p^n}}\left(\bigoplus_{k \in \bbZ_{\geq 0}} \uA\{ x^k_{C_{p^m}/C_{p^m}} \}  \right)\\
	& \cong N_{C_{p^m}}^{C_{p^n}}\left(\uA_{\amalg_{\bbZ_{\geq 0}} C_{p^m}/C_{p^m}}  \right)\\
	& \cong \uA_{ \Set^{C_{p^m}}(C_{p^n}, \amalg_{\bbZ_{\geq 0}} C_{p^m}/C_{p^m}) }.
	\end{align*}
	The statement now follows from an identification of the $C_{p^n}$-set $\Set^{C_{p^m}}(C_{p^n}, \amalg_{\bbZ_{\geq 0}} C_{p^m}/C_{p^m})$. 
	A $C_{p^m}$-equivariant map $C_{p^n}\to \amalg_{\bbZ_{\geq 0}} C_{p^m}/C_{p^m}$ is the same as a tuple of natural numbers $k_\gamma \in \bbZ_{\geq 0}$ indexed by $C_{p^n}/C_{p^m}$. 
	The corresponding element at level $C_{p^m}$ is identified as $\prod_{\gamma \in C_{p^n}/ C_{p^m}} (\gamma \cdot x)^{k_\gamma}$. 
	The $C_{p^n}$-action is by cyclic permutation of the indices. 
	Thus we observe that a tuple $\vec{v}$ has stabilizer at least $C_{p^k}$ for $m\leq k\leq n$ if and only if there is another tuple $\vec{w} \in \bbZ_{\geq 0}^{C_{p^n}/C_{p^k}}$ indexed by $C_{p^n}/C_{p^k}$ such that $v_{\gamma} = w_{\pi(\gamma)}$ for all $\gamma \in C_{p^n}/C_{p^m}$, where $\pi\colon C_{p^n}/C_{p^m}\to C_{p^n}/C_{p^{k}}$ is the projection. 
	Thus, we may index these tuples instead by the reduced form $\vec{w} \in \bbZ_{\geq 0}^{C_{p^n}/C_{p^k}}$, and denote the corresponding generator $\nm_{C_{p^m}}^{C_{p^k}}(x^{\vec{w}} )$. 
	All other tuples in the same $C_{p^n}$-orbit are obtained by cyclic permutation of the indices.
	
	This finishes the proof of the decomposition of $\uR$ as a free Mackey functor.  
\end{proof}

From this description, we can also concretely describe all the levels of this free incomplete Tambara functor individually:

\begin{corollary}\label{thm:FreeIncompleteTambaraDescription}
	The $\cO$-Tambara functor $\uR = \uA^{\cO}[x_{C_{p^n}/C_{p^m}}]$ may be described as follows:
	\begin{enumerate}[(a)]
		\item For $0 \leq \ell \leq m$, $\uR(C_{p^n}/C_{p^\ell})$ is the free $C_{p^{n-m}}$-algebra over $\uA(C_{p^n}/C_{p^\ell})$ on one generator:
		\[
			\uR(C_{p^n}/C_{p^\ell}) \cong \uA(C_{p^n}/C_{p^\ell})[x^{(0)}, x^{(1)}, \ldots, x^{(p^{n-m}-1)}].
		\]
		The Weyl group $W_{C_{p^n}}C_{p^\ell} \cong C_{p^{n-\ell}}$ acts via the projection $C_{p^{n-\ell}} \to C_{p^{n-m}}$. \\
		If $\vec{v} = (i_0,i_1,\ldots,i_{p^{n-m}-1})\in \z_{\geq 0}^{\times p^{n-m}}$, we write
		\[
			x^{\vec{v}} := (x^{(0)})^{i_0} (x^{(1)})^{i_1} \cdots (x^{(p^{n-m}-1)})^{i_{p^{n-m}-1}}.
		\]

		\item For $m+1 \leq  \ell \leq n$, we have
		\[	
			\uR(C_{p^n}/C_{p^\ell}) \cong \uA(C_{p^n}/C_{p^\ell})[q^k_{\vec{v}} : \vec{v} \in \z_{\geq 0}^{\times p^{n-k}}, \ m \leq k \leq \ell]/I_\ell,
		\]
		where $q^k_{\vec{v}}$ represents the element $\tr_{C_{p^k}}^{C_{p^\ell}}(\nm_{C_{p^m}}^{C_{p^k}} (x^{\vec{v}}))$. 
		The Weyl group $W_{C_{p^n}} C_{p^\ell} \cong C_{p^{n-\ell}}$ acts via any section $C_{p^{n-\ell}}\to C_{p^{n-k}}$ of the projection, such that $\gamma\in C_{p^{n-\ell}}$ sends the generator $q^k_{\vec{v}}$ to $q^k_{\gamma \vec{v}}$. 
		The submodule $I_\ell$ is generated by the following relations:

		\begin{enumerate}[(i)]

			\item $q^k_{\vec{0}} = \tr_{C_{p^k}}^{C_{p^\ell}}(\nm_{C_{p^m}}^{C_{p^k}} (1))\in \uA(C_{p^n}/C_{p^\ell})$, where $1\in \uA(C_{p^n}/C_{p^m})$;

			\item $q^k_{\vec{v}} = q^k_{\gamma \vec{v}}$, where $\gamma \in C_{p^{\ell-k}}$ acts on $\mathbb{Z}_{\geq 0}^{\times p^{n-k}}$ via the inclusion $C_{p^{\ell-k}} \to C_{p^{n-k}}$. 





			\item $q^k_{\vec{v}} q^{k'}_{\vec{w}} = \sum_{g \in C_{p^{\ell- k'}}} q^k_{\vec{v} + \sum_{\gamma  \in C_{p^{k'-k}}} g\gamma \vec{w}}$, where $k\leq k'$. 

		\end{enumerate}

		\item For $0 \leq \ell < \ell' \leq m$, we have
		\[
			\res_{C_{p^\ell}}^{C_{p^{\ell'}}}(x^{(i)}) = x^{(i)}, \quad \tr_{C_{p^\ell}}^{C_{p^{\ell'}}}(x^{(i)}) = [C_{p^{\ell'}}/C_{p^\ell}] x^{(i)}.
		\]

		\item For $m < k < \ell \leq n$, we have
		\[
			\res_{C_{p^m}}^{C_{p^\ell}}(q^k_{\vec{v}}) = \sum_{g \in C_{p^{\ell-k}}} \left( \prod_{\gamma \in C_{p^{k-m}}} x^{g\gamma \vec{v}} \right), \quad \tr_{C_{p^k}}^{C_{p^\ell}}(\nm_{C_{p^m}}^{C_{p^k}} (x^{\vec{v}})) = q^k_{\vec{v}}.
		\]
		Other transfers, norms and restrictions of the generators are determined by Mackey decomposition and exponential formulas.
	\end{enumerate}
\end{corollary}

\begin{remark}
	For $m+1 \leq \ell \leq n$, $\uR(C_{p^n}/C_{p^\ell})$ admits a smaller presentation. For instance, the infinite collection of generators
	\[
		\{q^\ell_{\vec{w}} : \vec{w} \in \mathbb{Z}_{\geq 0}^{\times p^{n-\ell}}\}
	\]
	is spanned under the Weyl group action and multiplication by
	\[
		\{n_{\vec{e_0}}, n_{\vec{e_1}}, \ldots, n_{\vec{e_{p^{n-\ell}-1}}} \}
	\]
	in view of (2c), where we abbreviated $q^\ell_{\vec{v}} = n_{\vec{v}}$.

	Moreover, on a summand indexed by $q^k_{\vec{v}}$, the quotients enforced through $I_\ell$ reduce the Burnside ring $\uA(C_{p^n}/C_{p^\ell})$ on this generator to $\uA(C_{p^n}/C_{p^k})$.
\end{remark}

\begin{example}[cf. {\cite[Lemma 3.2]{BH2019}}]
\label{example:freeC3TambaraUnderlyingGen}
	The Lewis diagram for the free $C_3$-Tambara functor $\uR = \uA^{\cO^\top}[x_{C_3/e}]$ on an underlying generator is as below:
	\[ \begin{tikzcd}[row sep=large]
		\bbZ[n, t_{ijk} \colon i,j,k \geq 0] / F
			\ar[d, bend right=50, "\res"{left}]\\
		\bbZ[x^{(0)},x^{(1)},x^{(2)}]
				\ar[u,"\nm" description] \ar[u, bend right=50, "\tr"{right}]
				\arrow[out=240,in=300,loop,swap,looseness=4, "C_3"]
	\end{tikzcd} \]
	In $\uR(C_3/C_3)$, the submodule $F$ is generated by the relations 
	\begin{align*} 
		t_{000}^2 &= 3 t_{000}, \\
		t_{ijk} &= t_{kij} = t_{jki},\\
		n t_{ijk} & = t_{i+1,j+1,k+1},\\
		t_{ijk} t_{pqr} &= t_{i+p,j+q,k+r} + t_{i+r,j+p,k+q} + t_{i+q,j+r,k+p},
	\end{align*}
	and the restriction, norm, and transfer are determined by the following: 
	\begin{align*}
		\res_e^{C_3}(n) &=  x^{(0)}x^{(1)}x^{(2)},\\
		\res_e^{C_3}(t_{ijk}) &= (x^{(0)})^i (x^{(1)})^j (x^{(2)})^k + (x^{(0)})^k (x^{(1)})^i (x^{(2)})^j + (x^{(0)})^j (x^{(1)})^k (x^{(2)})^i,\\[1em]
		\tr_e^{C_3}((x^{(0)})^i (x^{(1)})^j (x^{(2)})^k) & = t_{ijk},\\[1em]
		\nm_e^{C_3}(x^{(0)}) &= \nm_e^{C_3}(x^{(1)}) = \nm_e^{C_3}(x^{(2)}) = n.
	\end{align*}
	The Weyl action on the underlying level is given by $\gamma \cdot x^{(i)}= x^{(i+1)}$, with indices taken mod $3$ so $\gamma \cdot x^{(2)} = x^{(0)}$. 
\end{example}

\begin{example}[cf. {\cite[Lemma 3.6]{BH2019}}]
	\label{example:FreeCpTambaraFixedGen}
	The free $C_p$-Tambara functor $\uR = \uA^{\cO^\top}[x_{C_p/C_p}]$ on a fixed generator is as below:
	\[ 
	\begin{tikzcd}[row sep=large]
		\bbZ[t,n,x] / (t^2 - pt, tx^p - tn)
			\ar[d, bend right=50, "\res"{left}]\\
		\bbZ[x]
				\ar[u,"\nm" description] \ar[u, bend right=50, "\tr"{right}]
				\arrow[out=240,in=300,loop,swap,looseness=4, "C_p"]
	\end{tikzcd} 
	\]
	The restriction, norm, and transfer are determined by 
	\begin{align*}
		\res^{C_p}_e(x) &= x\\
		\res^{C_p}_e(n) &= x^p\\
		\tr_e^{C_p}(f) & = tf \\
		\nm_e^{C_p}(x) & = n\\
		\nm_e^{C_p}(a) & = a + \left(\tfrac{a^p - a}{p}\right) t
	\end{align*}
	Note that this example is not described by the results in this section since the underlying Mackey functor of this Tambara functor is not free.
\end{example}

\section{Tor, Hochschild homology, and K{\"a}hler differentials}
\label{section:HHisTOR}

We now understand our main players, the free incomplete Tambara functors over cyclic groups of prime power order. 
In the sequel, we will provide tools for computing Mackey-functor-valued $\uTor$ over these, but before that we explain how computing $\uTor$ is related to some other objects from equivariant algebra. 

\subsection{Tor, Hochschild homology, and K{\"a}hler differentials in classical algebra}

To begin, recall that in classical algebra, if $R$ is a $k$-algebra that is flat as a module over $k$, then the Hochschild homology of $R$ can be computed using a resolution of $R$ by $(R,R)$-bimodules \cite[1.1.13]{Lod1998}:
\begin{equation}
\label{HHasTor}
	\HH_*(R/k) \cong \Tor_*^{R \otimes_k R^\text{op}}(R, R).
\end{equation}
In fact, we can simplify further. 
The Hochschild homology can be computed using a Koszul complex -- a resolution of the trivial $R$-module $k$ by $R$-modules -- and then tensoring with $R$ (cf. \cite[Proof of Prop. 2.1]{MS93}): 
\begin{equation}
\label{HHTor2}
	\HH_*(R/k) \cong \Tor_*^{R \otimes_k R^\text{op}}(R, R) \cong R \otimes_k \Tor_*^{R}(k,k). 
\end{equation}

Moreover, if $R$ is smooth as a $k$-algebra, the Hochschild--Kostant--Rosenberg theorem \cite{HKR1962, Lod1998} implies that 
\begin{equation}
\label{HKR}
	\HH_n(R/k) \cong \Omega^n_{R/k},
\end{equation}
where $\Omega^n_{R/k}$ is the module of differential $n$-forms defined by
\[
	\Omega^n_{R/k} := \bigwedge^n \Omega^1_{R/k},
\]
with $\Omega^1_{R/k}$ the module of K{\"a}hler differentials. 

In this section, we will demonstrate the isomorphisms \eqref{HHasTor} and \eqref{HHTor2} when $\uR$ is an incomplete Tambara functor and $k$ is replaced by the Burnside functor $\uA$, as well as \eqref{HKR} when $n=1$ and $\uR$ is (only) a Green functor. 
We largely follow \cite{CE56}, modifying the arguments to work with equivariant algebra where necessary. 
We begin with a discussion of Hochschild homology for incomplete Tambara functors. 

\subsection{Hochschild homology of incomplete Tambara functors}

Let $\uR$ be an incomplete Tambara functor and let $\uM$ be an $\uR$-bimodule. 

\begin{definition}
	The \emph{cyclic nerve of $\uR$ with coefficients in $\uM$} is the simplicial $\uR$-module with $k$-simplices 
	\[ 
		[k] \mapsto \uM \boxtimes \uR^{\boxtimes k}.
	\]
	The face map $d_0$ is given by the right action of $\uR$ on $\uM$. 
	The faces $d_i$ for $1 \leq i \leq k-1$ are given by multiplication between the $i$-th and $(i+1)$-st factors. 
	The face $d_k$ is given by wrapping the last factor around to the front and acting on $\uM$ on the left. 
	
	The degeneracies $s_i$ are given by inserting the unit in the $(i+1)$-st copy of $\uR$. 
\end{definition}

\begin{definition}	
	The \emph{Hochschild homology of $\uR$ with coefficients in $\uM$} is the $\uR$-module-valued homology of the cyclic nerve. We denote the $n$-th homology $\uR$-module by $\HH_n(\uR, \uM)$. 
	
	When $\uM = \uR$ with $\uR$-bimodule structure by left and right multiplications, we write $\HH_n(\uR) := \HH_n(\uR,\uR)$. 
\end{definition}

\begin{remark}
	The Hochschild homology we study here is different than that studied in \cite{BGHL2019}, where a twisting is applied to the left module structure of $\uM$ before taking the cyclic nerve. 
	This in particular requires that $G$ is a finite subgroup of $S^1$, while there is no such restriction in the definition above (although all of our calculations below are for cyclic $p$-groups). 
	Nevertheless, many of the arguments that apply to twisted Hochschild homology apply here as well. 
\end{remark}

\begin{remark}
	The Hochschild homology of an incomplete Tambara functor depends only on its underlying Green functor; the norms play no role in the Hochschild complex. 
	Therefore, the first Hochschild homology of an incomplete Tambara functor as defined above will not agree with its genuine equivariant K\"ahler differentials \cite[Definition 5.4]{Hil2017} unless it is itself a Green functor. 
	This is in contrast to the case of ordinary algebra, where the first Hochschild homology of an algebra is isomorphic to its K\"ahler differentials. 
\end{remark}

\subsection{Mapping theorem}

In this section, we prove a technical theorem which will be used to prove the isomorphism \eqref{HHTor2}. 
Our discussion is a straightforward adaptation of the analogous nonequivariant discussion of the mapping theorem in \cite{CE56}. 
Let $\uR$, $\uS$, $\uQ_{\uR}$, and $\uQ_{\uS}$ be Green functors. 
Let $\epsilon_{\uR} : \uR \to \uQ_{\uR}$ and $\epsilon_{\uS} : \uS \to \uQ_{\uS}$ be maps of Green functors. 
Let $\phi : \uR \to \uS$ and $\psi: \uQ_{\uR} \to \uQ_{\uS}$ be maps of Green functors such that we have a commutative diagram
\[
	\begin{tikzcd}
		\uR 
			\arrow{r}{\epsilon_{\uR}} 
			\arrow{d}{\phi} 
			& 
		\uQ_{\uR} 
			\arrow{d}{\psi} 
		\\
		\uS 
			\arrow{r}{\epsilon_{\uS}} 
			& 
		\uQ_{\uS}.
	\end{tikzcd}
\]
For any right $\uS$-module $\uM$, we construct a map
\[
	F^\phi: \uTor^{\uR}(\uM,\uQ_{\uR}) \to \uTor^{\uS}(\uM,\uQ_{\uS})
\]
as follows. 
Let $\uF_{\uR}$ be an $\uR$-projective resolution of $\uQ_{\uR}$ and let $\uF_{\uS}$ be an $\uS$-projective resolution of $\uQ_{\uS}$. 
Define
\[
	g: \uS \boxtimes_{\uR} \uQ_{\uR} \to \uQ_{\uS}
\]
by $s \otimes x \mapsto s \psi(x)$. 
Then $\uS \boxtimes_{\uR} \uF_{\uR}$ is an $\uS$-projective complex of $\uS \boxtimes_{\uR} \uQ_{\uR}$. 
Applying the lifting criterion for projective Mackey functors (cf. \cite[Prop. 2.2.11]{Lee2019}), the usual argument shows that there is a map
\[
	\tilde{g} : \uS \boxtimes_{\uR} \uF_{\uR} \to \uF_{\uS}
\]
over $g$, unique up to homotopy, which yields the desired comparison map
\[
	F^{\phi}: H_*(\uM \boxtimes_{\uR} \uF_{\uR}) 
		\cong H_*(\uM \boxtimes_{\uS} (\uS \boxtimes_{\uR} \uF_{\uR})) 
		\to H_*(\uM \boxtimes_{\uR} \uF_{\uS}).
\]

\begin{theorem}[Mapping theorem for Green functors]
	The map $F^\phi$ is an isomorphism of graded Mackey functors for all right $\uS$-modules $\uM$ if and only if
	\begin{enumerate}[(a)]
		\item $g: \uS \boxtimes_{\uR} \uQ_{\uR} \to \uQ_{\uS}$ is an isomorphism, and
		\item $\uTor_n^{\uR}(\uS,\uQ_{\uR}) = 0$ for all $n >0$. 
	\end{enumerate}
	Moreover, if these conditions are satisfied, then for any projective resolution $\uF_{\uR}$ of $\uQ_{\uR}$, the complex $\uS \boxtimes_{\uR} \uF_{\uR}$ with augmentation $\uS \boxtimes_{\uR} \uF_{\uR} \to \uS \boxtimes_{\uR} \uQ_{\uR} \cong \uQ_{\uS}$ is an $\uS$-projective resolution of $\uQ_{\uS}$. 
\end{theorem}

\begin{proof}
	The proof is identical to the proof of the analogous classical result (\cite[Thm. VIII.3.1]{CE56}). 

	If $F^\phi$ is an isomorphism, then taking $\uM = \uS$ proves $(1)$ and $(2)$. 

	Assume (a) and (b) hold. 
	Let $\uF_{\uS}$ be an $\uS$-projective resolution of $\uQ_{\uR}$, so 
	\[
		H_n(\uS \boxtimes_{\uR} \uF_{\uR}) \cong \uTor_n^{\uR}(\uS,\uQ_{\uR}) = 0
	\]
	for $n >0$. 
	Then (a) and (b) express the fact that
	\[
		\uS \boxtimes_{\uR} \uF_{\uR} 
		\to \uS \boxtimes_{\uR} \uQ_{\uR} \cong \uQ_{\uR}
	\]
	is an $\uS$-projective resolution of $\uQ_{\uS}$. 
	Taking $\uF_{\uS} = \uS \boxtimes_{\uR} \uF_{\uR}$, we can take $\tilde{g} \colon \uS \boxtimes_{\uR} \uF_{\uR} \to \uF_{\uS}$ to be the identity, so $F^\phi$ is an isomorphism. 
\end{proof}

\subsection{Identifying Hochschild Homology and Tor}

We will now apply the mapping theorem to identify $\uHH$ with known $\uTor$ groups under favorable circumstances; in particular, we will produce the isomorphisms \eqref{HHasTor} and \eqref{HHTor2}. 

Suppose $\uK$ is a Green functor and $\uR$ is an augmented $\uK$-algebra. Suppose further that we are given a map $E : \uR \to \uR^e := \uR \boxtimes_{\uK} \uR^{op}$ so that there is a commutative diagram of Green functors
\begin{equation}
	\label{Eqn:ReRDiagram}
	\begin{tikzcd}
	\uR 
		\arrow{r}{\epsilon} 
		\arrow{d}{E} 
		& 
	\uK 
		\arrow{d}{\eta} 
		\\
	\uR^e 
		\arrow{r}{\rho} 
		& 
	\uR,
	\end{tikzcd}
\end{equation}
where $\epsilon$ is the augmentation of $\uR$, $\eta$ is part of the $\uK$-algebra structure of $\uR$, and $\rho$ is the augmentation of $\uR^e$ as an $\uR$-algebra. 
Here, commutativity implies that there is an inclusion $E\uI \subseteq \uJ$, where $\uI$ and $\uJ$ are the augmentation ideals of $\epsilon$ and $\rho$, respectively. 
Since $(E, \eta)$ is a map of augmented Green functors, we obtain a map
\[
	F^E: \uTor_n^{\uR}(\uM_E, \uK) \to \uTor^{\uR^e}_n({}_\epsilon\uM, \uR),
\]
for any $\uR$-bimodule $\uM$, where $\uM_E$ is the right $\uR$-module obtained by regarding $\uM$ as a right $\uR^e$-module and then defining the $\uR$-module structure via $E$, and ${}_\epsilon \uM$ is $\uM$ regarded as a left $\uR$-module via $ra = (\epsilon(r))a$. 

\begin{theorem}
	\label{Thm:TorEnv}
	Assume 
	\begin{enumerate}[(E.1)]
		\item $J = \uR^e_E \uI$;
		\item $\uR^e_E$ is projective as a right $\uR$-module. 
	\end{enumerate}
	Then $F^E$ is an isomorphism, and for any $\uR$-projective resolution $\uF \to \uK$, 
	\[
		\uR^e_E \boxtimes_{\uR} \uF \to \uR
	\]
	is a $\uR^e$-projective resolution. 
\end{theorem}

\begin{proof}
	This follows from the mapping theorem by the same proof used to prove \cite[Thm. X.6.1]{CE56}. We check that (a) and (b) hold. 

	The exact sequence
	\[
		\uI \to \uR \to \uK \to 0
	\]
	yields an exact sequence
	\[
		\uR^e_E \boxtimes_{\uR} \uI \to \uR^e_E \to \uR_E^e \boxtimes_{\uR} \uK \to 0
	\]
	which implies
	\[
		\uR^e_E \boxtimes_{\uR} \uK \cong \coker(\uR^e_E \boxtimes_{\uR} \uI \to \uR_E^e) \cong \uR^e_E / \uR^e_E \uI.
	\]
	But by (E.1) and the fact that $\uJ$ is the augmentation ideal for $\rho$, we have
	\[
		\uR^e_E / \uR^e_E \uI \cong \uR^e_E / \uJ \cong \uR,
	\]
	so 
	\[
		\uR^e_E \boxtimes_{\uR} \uK \cong \uR,
	\]
which proves (a). 

	For (b), it follows from (E.2) that 
	\[
		\uTor^{\uR}_n(\uR^e_E, \uK) = 0
	\]
	for $n>0$. 

	Therefore both conditions of the mapping theorem hold, so $F^E$ is an isomorphism and the claim about resolutions holds. 
\end{proof}

\begin{proposition}
	Let $\uR$ be a free incomplete Tambara functor whose underlying Mackey functor is projective. Then there is an isomorphism of graded Mackey functors
	\[
		\uTor_*^{\uR \boxtimes \uR}(\uR,\uR) \cong \uR \boxtimes \uTor_*^{\uR}(\uA,\uA).
	\]
\end{proposition}

\begin{proof}
	We adapt the proof from \cite[Prop. 2.1]{MS93} of the analogous fact from ordinary algebra. 
	Let $\uR = \uA^{\cO}[x_{G/H}]$ with augmentation $\epsilon: \uA^{\cO}[x_{G/H}] \to \uA$, $\eta: \uA \to \uA^{\cO}[x_{G/H}]$ the Green structure map, and $\rho: \uR^e \to \uR$ the multiplication.
	We define
	\[
		E : \uR \to \uR^e
	\]
	to be the map (on underlying Green functors) corresponding to $x \otimes 1 - 1 \otimes x$ under the isomorphism 
	\[
		\Tamb^{\cO}(\uA^{\cO}[x_{G/H}], \uA^{\cO}[x_{G/H}]^e) \cong \uA^{\cO}[x_{G/H}]^e(G/H).
	\]
	Then the diagram \eqref{Eqn:ReRDiagram} commutes, so $(E,\eta)$ is a map of augmented Green functors which induces a map
	\[
		F^E: \uTor_n^{\uR}(\uM_E,\uA) 
			\to \uTor^{\uR^e}_n({}_{\epsilon}\uM,\uR)
	\]
for any $\uR$-bimodule $\uM$. 
	By construction (cf. the equivalence between (E.1) and (E.1') in \cite[Sec. X.6]{CE56}), the augmentation ideals of $\epsilon$ and $\rho$ satisfy $\uJ = \uR^e_E \uI$, so (E.1) holds. 
	Since $\uR$ is underlying projective, $\uR^e_E$ is projective as a right $\uR$-module, so (E.2) holds. 
	Therefore $F^E$ is an isomorphism by \cref{Thm:TorEnv}.

	Finally, we take $\uM = \uR$. 
	Since $\uR$ is commutative, $\uR_E$ is $\uR$ equipped with the trivial $\uR$-module structure. 
	Therefore we have an isomorphism
	\[
		\uTor^{\uR}(\uR_E,\uA) \cong \uR \boxtimes \uTor^{\uR}(\uA,\uA).
	\]
	Composing this with the isomorphism $F^E$ proves the proposition. 
\end{proof}

\begin{corollary}\label{corollary:HHTorIdenfitication}
	Let $\uA^{\cO}[x_{G/H}]$ be a free incomplete Tambara functor whose underlying Mackey functor is projective. 
	Then there is an isomorphism of graded Mackey functors
	\[
		\uHH_*(\uA^{\cO}[x_{G/H}]) \cong \uA^\cO[x_{G/H}] \boxtimes \uTor_*^{\uA^{\cO}[x_{G/H}]}(\uA,\uA).
	\]
\end{corollary}

\subsection{A simple case of the HKR isomorphism}

We will now prove a simple case of the Hochschild--Kostant--Rosenberg theorem in equivariant algebra. 

To begin, we recall the module of genuine K{\"a}hler differentials over an incomplete Tambara functor. 
For details, we refer the reader to Hill's paper \cite{Hil2017} for the (complete) Tambara case and Leeman's thesis \cite{Lee2019} for the more general incomplete Tambara case. 

\begin{definition}[{\cite[Def. 3.2.1]{Lee2019}}]
	\label{DefDerivations}
	Let $\uR$ be an $\cO$-Tambara functor, $\uS$ an $\uR$-algebra, and $\uM$ an $\uS$-module. 
	A map $d: \uS \to \uM$ is an \emph{$\cO$-genuine $\uR$-derivation} if
	\begin{enumerate}[(a)]
		\item it is a map of Mackey functors,
		\item the composite $\uR \to \uS \to \uM$ is the zero map, and
		\item the map turns all admissible norms and products into transfers and sums in the following sense. 
		Let $f : X \to Y$ be a map of $G$-sets representing an admissible norm map in $\cO$, including products. 
		Let
		\[
			\pi_1, \pi_2 : X \times_Y X - \Delta(X) \to X
		\]
		be the projections. 
		Then for $s \in \uS(X)$, we require
		\[
			d(N_f(s)) = T_f(N_{\pi_2}R_{\pi_1}(s) \cdot d(s)).
		\]
	\end{enumerate}
\end{definition} 

\begin{definition}[{\cite[Sec. 3.3]{Lee2019}}]
	\label{DefKaehlerDiffs}
	Let $\uR$ be an $\cO$-Tambara functor and $\uS$ an $\uR$-algebra. 
	Let $\uI$ denote the kernel of the multiplication
	\[
		\mu: \uS \boxtimes_{\uR} \uS \to \uS.
	\]
	Let $\uI^{>1}$ be the smallest ideal of $\uS$ containing every element $N_f(i)$, where $i \in \uI(X)$ and $f: X \to Y$ is a $2$-surjective $\cO$-admissible map, i.e., an $\cO$-admissible map such that each fiber $f^{-1}(y)$ has cardinality at least $2$. 
	The \emph{module of genuine K{\"a}hler differentials} is defined by
	\[
		\Omega^{1,G}_{\uS/\uR} := \uI/\uI^{>1}.
	\]
\end{definition}

\begin{remark}
	As in classical algebra, the characterizing property of $\Omega^{1,G}_{\uS/\uR}$ is that it corepresents derivations: there is an isomorphism 
	\[
		\Hom_{\uS}(\Omega^{1,G}, \uM) \cong \Der_{\uR}(\uS,\uM)
	\]
	for any $\uS$-module $\uM$ by \cite[Prop. 3.3.2]{Lee2019}. 
	This isomorphism is induced by composition with the \emph{universal derivation}
	\[
		d\colon \uS \to \Omega^{1, G}_{\uS/\uR}, 
	\]
	defined as
	\[
		d(s) = [s\otimes 1 - 1\otimes s].
	\]
\end{remark}


\begin{lemma}[{\cite[Proposition 5.5]{Hil2017}}]
	\label{KaehlerGeneratedByD}
	The $\uS$-module $\Omega^{1,G}_{\uS/\uR}$ is generated by the image of $d \colon \uS \to \Omega^{1,G}_{\uS/\uR}$. 
\end{lemma}

\begin{lemma}
	\label{KahlerGensAndRelations}
	Let $\uR$ and $\uS$ be Green functors. 
	Then $\Omega^{1,G}_{\uS/\uR}$ is the free $\uS$-module generated by symbols $da$ with $a \in \uS(G/H)$ for all $H \leq G$, subject to the relations 
	\begin{align*}
		d(r \cdot s_1 + s_2) &= r d(s_1) + d(s_2)\\
		d(s_1 \cdot s_2) &= s_1 \cdot d(s_2) + s_2 \cdot d(s_1) 
	\end{align*}
	for all $s_1, s_2 \in \uS(G/H)$ and $r \in \uR(G/H)$. 
\end{lemma}

Note that the previous lemma only holds for Green functors; for general incomplete Tambara functors, we must add the relation 
\begin{equation}
	\label{TambaraExtraRelation}
	d(N_f(s)) = T_f(N_{\pi_2}R_{\pi_1}(s) \cdot d(s)),
\end{equation}
as in \cref{DefDerivations}.

\begin{theorem}
	\label{HochschildHomologyAsKaehlerDiffs} 
	Let $\uR$ be a Green functor and $\uS$ be an $\uR$-algebra. 
	There is an isomorphism of $\uS$-modules
	\(
		\HH_1(\uS/\uR) \cong \Omega^{1,G}_{\uS/\uR}.
	\)
\end{theorem}

\begin{proof}
	Consider the morphism of $\uS$-modules 
	\[ 
		\varphi \colon \uS \boxtimes_{\uR} \uS \xto{1 \boxtimes d} \uS \boxtimes_{\uR} \Omega^{1,G}_{\uS/\uR} \to \Omega^{1,G}_{\uS/\uR} 
	\]
	where the second homomorphism is the $\uS$-module structure on the K\"ahler differentials. 
	Because $\uS \boxtimes_{\uR} \uS$ is generated as a Mackey functor by $\uS(G/H) \otimes_{\uR(G/H)} \uS(G/H)$ as $H$ ranges over subgroups of $G$, we may describe this map by its effect on generators. 
	For  
	\[ 
		s_1 \otimes s_2 \in \uS(G/H) \otimes_{\uR(G/H)} \uS(G/H), 
	\]
	we have 
	\[
		\varphi(s_1 \otimes s_2) = s_1 \cdot d(s_2).
	\]
	Note that this map is surjective by \cref{KaehlerGeneratedByD}.

	We claim that $\varphi$ descends to a homomorphism 
	\[ 
		\bar \varphi \colon \HH_1(\uS/\uR) \to \Omega^{1,G}_{\uS/\uR}. 
	\]
	To verify this, we must check that $\varphi$ respects the relations in 
	\[ 
		\HH_1(\uS/\uR) \cong \uS \boxtimes_{\uR} \uS / \im(b), 	
	\]
	where 
	\[ 
		b \colon \uS \boxtimes_{\uR} \uS \boxtimes_{\uR} \uS \to \uS \boxtimes_{\uR} \uS 
	\]
	is the Hochschild differential. 
	Note that $\uS \boxtimes_{\uR} \uS \boxtimes_{\uR} \uS$ is generated as a Mackey functor by all triple tensor products 
	\[ 
		\uS(G/H) \otimes_{\uR(G/H)} \uS(G/H) \otimes_{\uR(G/H)} \uS(G/H), 
	\]
	and $b$ is determined by 
	\[ 
		b(s_1 \otimes s_2 \otimes s_3) = s_1s_2 \otimes s_3 - s_1 \otimes s_2 s_3 + s_3 s_1 \otimes s_2. 
	\]
	Hence, $\im(b)$ is generated by elements of this form. 
	The homomorphism $\bar \varphi$ is well-defined if $\varphi$ sends elements of this form to zero. 
	We check this: 
	\begin{align*}
		\varphi(s_1s_2 \otimes s_3 - s_1 \otimes s_2 s_3 + s_3 s_1 \otimes s_2) 
			&= s_1s_2\cdot d(s_3) - s_1 \cdot d(s_2 s_3) + s_3 s_1 \cdot d(s_2) 
		\\
			& = s_1s_2\cdot d(s_3) - s_1  s_2 \cdot d(s_3) - s_1 s_3 \cdot d(s_2) + s_3 s_1 \cdot d(s_2) 
		\\
			& = 0. 
	\end{align*}
	Hence, this is well-defined. 

	Let $F$ be the free $\uS$-module on the symbols $d(a)$. 
	$\Omega^{1,G}_{\uS/\uR}$ is the quotient of $F$ by the relations in \cref{KahlerGensAndRelations}. 
	We define
	\[
		\psi \colon F \to \HH_1(\uS/\uR) 
	\]
	by 
	\[
		s_1 \cdot d(s_2) \mapsto s_1 \otimes s_2 + \im(b) \in \uS(G/H) \otimes_{\uR(G/H)} \uS(G/H) / \im(b) 
	\]
	for $s_1, s_2 \in \uS(G/H)$. 
	Since the levelwise tensors $\uS(G/H) \otimes_{\uR(G/H)} \uS(G/H)$ generate $\uS \boxtimes_{\uR} \uS$ as a Mackey functor, this yields an element of the box product, and by passing to the quotient, an element of $\HH_1(\uS/\uR)$. 
	To get a homomorphism 
	\[ 
		\bar \psi \colon \Omega^{1,G}_{\uS/\uR} \to \HH_1(\uS/\uR), 
	\]
	we must check that $\psi$ sends the relations of $\Omega^{1,G}_{\uS/\uR}$ (\cref{KahlerGensAndRelations}) to zero.
	\begin{align*}
		\psi( s_1 \cdot d(s_2) + s_2 \cdot d(s_1)- d(s_1 \cdot s_2)) 
			&= s_1 \otimes s_2 + s_2 \otimes s_1 - 1 \otimes s_1s_2  + \im(b).
	\end{align*}
	Note that $b(1 \otimes s_1 \otimes s_2)$ is exactly the right hand side, so we get zero in the quotient. 
	Similarly, 
	\begin{align*}
		\psi(r d(s_1) + d(s_2) - d(r s_1 + s_2)) 
			&= r \otimes s_1 + 1 \otimes s_2 - 1 \otimes (rs_1 + s_2) 
		\\
			& = 1 \otimes rs_1 + 1 \otimes s_2 - 1 \otimes (rs_1 + s_2) 
		\\
			& = 0
	\end{align*}
	Therefore, $\bar\psi$ is a well-defined homomorphism. 

	We have defined homomorphisms of $\uS$-modules
	\begin{align*}
		\bar\psi \colon \Omega^{1,G}_{\uS/\uR} 
			\to \HH_1(\uS/\uR),
			\quad s_1 \cdot d(s_2) 
			\mapsto s_1 \otimes s_1 + \im(b).
			\\
		\bar\varphi \colon \HH_1(\uS/\uR) 
			\to \Omega^{1,G}_{\uS/\uR},
			\quad s_1 \otimes s_2 + \im(b) 
			\mapsto s_1 \cdot d(s_2)
	\end{align*}
	They are clearly inverse to each other, and moreover send generators to generators. 
\end{proof}

\begin{remark}
	\label{HKRnaiveKaehlerDiffs}
	In the general incomplete Tambara case, where $\Omega^{1,G}_{\uS/\uR}$ has an extra relation \eqref{TambaraExtraRelation} that is not present in $\HH_1(\uS/\uR)$, we only get a surjection from Hochschild homology onto the K\"ahler differentials. 
	This is because the Hochschild homology of an incomplete Tambara functor only depends on its underlying Green functor, while the K\"ahler differentials take into account all of the norms.
	
	The same proof as above shows, however, that the Hochschild homology is always isomorphic to a \emph{naive} module of K\"ahler differentials defined as in usual commutative algebra as $\uI/\uI^2$, where $\uI$ is the kernel of the multiplication map
	\[
	 \mu \colon \uS \boxtimes_{\uR} \uS \to \uS. 
	\]
	In comparison to the definition of the genuine K\"ahler differentials in \cref{DefKaehlerDiffs}, for the naive K\"ahler differentials we take the quotient by the square of the ideal $\uI$, which is generated by elements of the form $i \cdot j$ for $i, j \in \uI(G/H)$ but does not contain non-trivial norms of elements in $\uI$.
\end{remark}

\subsection{Free genuine modules}

As noted in \cref{HKRnaiveKaehlerDiffs}, we can identify Hochschild homology over an incomplete Tambara functor with a naive module of K\"ahler differentials. 
Since we take a smaller quotient compared to the genuine K\"ahler differentials, this module carries additional norms. 
In this short section, we show that these norms give the module of naive K\"ahler differentials the structure of a genuine $\uR$-module. 
Moreover, we explicitly describe some free genuine modules over the Burnside functor $\uA$ and over the constant Tambara functor $\ubbZ$. 
The results of this section will not be used anywhere else in the paper.

\begin{definition}[{\cite[Definition 14.3]{Str2012}}]	
	\label{DefGenuineModules}
	Let $\uR$ be an incomplete Tambara functor. 
	Then the category of augmented $\uR$-algebras has products given by the fiber product over $\uR$.
	We define the category of genuine $\uR$-modules as the category of abelian group objects in augmented $\uR$-algebras.
\end{definition}

For a Green functor, the categories of naive and genuine modules are equivalent via the augmentation ideal. 
This generalizes the fact from commutative algebra that the category of modules is the abelianization of the category of augmented algebras, leading to the definition of Andr\'e-Quillen homology \cite{Qui70}. 
For an explanation see \cite[Proposition 14.7f]{Str2012}. 
However, for (incomplete) Tambara functors, there is a difference. 
To see examples of non-trivial genuine modules, we consider modules of naive K\"ahler differentials. 
In the following, we also call any augmentation ideal of an abelian group object in augmented $\uR$-algebras a genuine $\uR$-module.

\begin{proposition}
	Let $\uR$ be an incomplete Tambara functor and $\uS$ be an $\uR$-algebra. 
	Then the module of naive K\"ahler differentials $\Omega_{\uS/\uR}^{1}$ is (the augmentation ideal of) a genuine $\uS$-module.
\end{proposition}

\begin{proof}
	The module of naive K\"ahler differentials by definition is the augmentation ideal of the induced multiplication map
	\[ 
		\mu\colon \uS \boxtimes_{\uR} \uS/I^2 \to \uS, 
	\]
	where $I$ is the kernel of the multiplication $\uS\boxtimes_{\uR} \uS\to \uS$ itself. 
	Since $I^2$ is an incomplete Tambara ideal, $\uT = \uS \boxtimes_{\uR} \uS/I^2$ inherits the structure of an incomplete Tambara functor from $\uS\boxtimes_{\uR} \uS$. 
	The left unit $id \boxtimes 1\colon \uS \to \uT$ is a split of the multiplication and thus makes $\uT$ into an augmented $\uS$-algebra.
	
Thus, we now need to show that $\uT$ carries the structure of an abelian group object in augmented algebras. 
	The left unit splits $\uT \cong \uS\oplus \Omega^1_{\uS/\uR}$ as Mackey functors.
	We now define the group structure using the addition and neutral element in $\Omega_{\uS/\uR}^1$. 
	Explicitly, the group operation is given as 
	\[ 
		\alpha \colon 
			(\uS\oplus \Omega^1_{\uS/\uR}) \times_{\uS} (\uS\oplus \Omega^1_{\uS/\uR}) 
			\to \uS\oplus \Omega^1_{\uS/\uR}, 
			((s, x), (s, y))\mapsto (s, x+y). 
	\]
	We need to check that this map is a map of incomplete Tambara functors. 
	We show that it is compatible with the norm maps. 
	The remaining details to show that $\uT$ becomes an abelian group object with this operation are straight-forward to check.

	Let $K\leq H \leq G$ be subgroups, $s\in \uS(G/K)$ and $x, y\in \Omega^1_{\uS/\uR}(G/K)$. 
	We need to check that 
	\[ 
		\alpha(N_K^H(s, x), N_K^H(s,y)) = N_K^H(s, x+y). 
	\]
	For this we first study in general how a norm behaves on a sum $a+b$ for elements $a,b \in \uT(G/K)$. 
	The norm of a sum has been studied by Mazur \cite{HM2019}; we include some details below. 
	This norm of a sum is calculated by considering the exponential diagram
	\[
		\begin{tikzcd}
			G/K 
				\arrow[d, "f", swap] 
				& 
			G/K \amalg G/K 
				\arrow[l, "\nabla", swap] 
				& 
			X = \{ (gK, \rho\colon f^{-1}(gH) \to \{1,2\}) 
					\mid gK \in G/K \} 
				\arrow[l, "\eval", swap] 
				\arrow[d, "\pr"]
				\\
			G/H 
				&
				& 
			Y \cong \{ (gH, \rho\colon f^{-1}(gH) \to \{1,2\}) 
					\mid gH \in G/H \}. 
			\arrow[ll, "\pr"]
		\end{tikzcd} 
	\]
	Here $f\colon G/K \to G/H$ is the projection map inducing the norm $N_K^H$, and the map 
	\[
		\rho\colon f^{-1}(gH) \to \{1,2\}
	\] 
	describes a section of the fold map $\nabla$ over the fiber of $gH$. 
	In the $G$-set $Y$, choosing the two constant maps as $\rho$ describes an orbit isomorphic to $G/H$ each, and on these the right vertical projection is a copy of $f\colon G/K \to G/H$. 
	If the map $\rho$ is not constant, suppose that $gK$ and $g'K$ are two elements in $f^{-1}(gH)$ such that $\rho(gK) \neq \rho(g'K)$. 
	Then $\eval(gK, \rho)$ and $\eval(g'K, \rho)$ lie in different copies of the disjoint union $G/K \amalg G/K$ and thus the elements $(gK, \rho)$ and $(g'K, \rho)$ lie in different $G$-orbits in the pre-image of $(gH, \rho)$ under the right vertical map. 
	Thus the norm on this part of the vertical projection decomposes as a product.
	
	In total, we obtain that the norm on a sum has the form
	\[ 
		N_K^H(a+b) = N_K^H(a)+ N_K^H(b) + r, 
	\]
	where the remainder term $r$ is a sum of terms containing (transfers of) products of norms of both $a$ and $b$. 
	Applying this to a norm of $s\boxtimes 1 + x$ with $s\in \uS(G/K)$ and $x\in \Omega_{\uS/\uR}^1(G/K)$, we conclude that the remainder term is contained in $\Omega_{\uS/\uR}^1(G/K)$, since it contains factors from this ideal in $\uT$. 
	Thus the $\uS$-part of both $\alpha(N_K^H(s, x), N_K^H(s,y))$ and $N_K^H(s, x+y)$ is given as $N_K^H(s)$. 
	For the part contained in $\Omega_{\uS/\uR}^1$, we consider the norm of $x+y$ for $x, y\in \Omega_{\uS/\uR}^1(G/K')$ for any subgroup $K' \leq G$. 
	Here, the rest term vanishes as it consists of products of elements in $I$ and we take the quotient by $I^2$. 
	Hence $N_{K'}^{H'}(x+y) = N_{K'}^{H'}(x)+ N_{K'}^{H'}(y)$ for any $K'\leq H'\leq G$. 
	Combining these observations we obtain that indeed $\alpha(N_K^H(s, x), N_K^H(s,y)) = N_K^H(s, x+y)$ and thus $\uT$ has the structure of an abelian group object in augmented $\uS$-algebras and $\Omega_{\uS/\uR}^{1}$ is a genuine $\uS$-module.
\end{proof}

\begin{remark}
	The above proof shows that in fact the condition that an augmented algebra has the structure of an abelian group object is equivalent to all products on its augmentation ideal vanishing. 
	This is analogous to the same statement for commutative rings, and the calculation of the norms in the above proof make it clear that the vanishing of the products is enough to make the structure map compatible with the additional norms. 
	The condition to be a naive module over an incomplete Tambara functor additionally asks for all norms on the augmentation ideal to vanish. 
	As shown by Hill \cite{Hil2017}, the augmented algebra then even has the structure of a Mackey functor object.
\end{remark}

Classically, we have the calculation $\Omega^1_{R[x]/R} \cong R[x]\{dx\}$ of the Kähler differentials on a polynomial algebra as a free module. 
	Using the various isomorphism linking K\"ahler differentials, derivations and square-zero extensions also in the case of Tambara functors, as exhibited in \cite{Hil2017} and \cite[Corollary 3.3.9]{Lee2019}, we obtain the isomorphism $\Omega^{1, G}_{\uR[x_{G/H}]/\uR} \cong \uR[x_{G/H}]\{dx_{G/H}\}$. 
	This identifies the genuine K\"ahler differentials with a free naive module. 
	Generalizing this, the naive K\"ahler differentials of a polynomial algebra are also free as genuine modules.

\begin{proposition}
	\label{prop:naiveKahlerDiffsFree}
	Let $\uR$ be an incomplete Tambara functor and $\uR[x_{G/H}]$ be a free $\uR$-algebra on a generator at level $G/H$. 
	Then the module of naive K\"ahler differentials $\Omega_{\uR[x_{G/H}]/\uR}^{1}$ is free on a generator at level $G/H$ as a genuine $\uR[x_{G/H}]$-module. 
	Explicitly, evaluation on the element $dx_{G/H}$ defines an isomorphism
	\[ 
		\Hom_{\gen}(\Omega_{\uR[x_{G/H}]/\uR}^{1}, \uM) \xto{\cong} \uM(G/H)
	\]
	for any genuine $\uR[x_{G/H}]$-module $\uM$.
\end{proposition}

\begin{proof}
	We denote by $\uM$ the augmentation ideal of an abelian group object $\uS$ in the category of augmented $\uR[x_{G/H}]$-algebras.
	The module of naive K\"ahler differentials is the augmentation ideal of the augmented algebra $\uT = \uR[x_{G/H}] \boxtimes_{\uR} \uR[x_{G/H}]/I^2$, where $I$ is the kernel of the multiplication map $\uR[x_{G/H}] \boxtimes_{\uR} \uR[x_{G/H}] \to \uR[x_{G/H}]$. 
	We notice that $\uR[x_{G/H}] \boxtimes_{\uR} \uR[x_{G/H}]$ is the free $\uR[x_{G/H}]$-algebra on a single generator in degree $G/H$, given by the second polynomial generator. 
	To highlight the different tensor factors, we denote the generator of the second tensor factor by $y_{G/H}$. 
	For the augmented $\uR[x_{G/H}]$-algebra $\uS$, evaluation at $y$ thus defines an isomorphism
	\[ 
		\Hom_{\uR[x_{G/H}]}(\uR[x_{G/H}] \boxtimes_{\uR} \uR[y_{G/H}], \uS) \cong \uS(G/H), 
		\]
where the $\Hom$-sets are taken as $\uR[x_{G/H}]$-algebras. 
	Equivalently, such a morphism is uniquely determined by where it sends the element $x-y$.
	Moreover, this element lies in the augmentation ideal of the augmented algebra $\uR[x_{G/H}] \boxtimes_{\uR} \uR[y_{G/H}]$, and a morphism is a morphism of augmented algebras if and only if it maps $x-y$ into the augmentation ideal $\uM$ of $\uS$.

	Finally, we notice that since $\uS$ is an abelian group object in augmented algebras, all products on its augmentation ideal vanish. 
	This implies that a morphism $\uR[x_{G/H}] \boxtimes_{\uR} \uR[y_{G/H}]\to \uS$ of augmented algebras automatically factors through $\uT = \uR[x_{G/H}] \boxtimes_{\uR} \uR[y_{G/H}] / I^2$.
	It is a straight-forward argument that this induced morphism is then automatically a morphism of abelian group objects. 
	In total, we have thus shown that evaluation on the element $x-y = dx$ induces an isomorphism 
	\[ 
		\Hom_{\gen}(\Omega_{\uR[x_{G/H}]/\uR}^{1}, \uM) \xto{\cong} \uM(G/H).\qedhere
	\]
\end{proof}

We can now explicitly calculate the free genuine modules over the Burnside functor $\uA$ using the description as K\"ahler differentials of polynomial algebras.

\begin{example}
	\label{example:free_genuine_module}
	We consider the case $G=H=C_p$, and denote $\uS=\uA[x_{C_p/C_p}]$. 
	Then, we have the definition 
	\[ 
		\Omega^1_{\uS/\uA} = I/I^2, \quad I = \ker (\uS\boxtimes \uS\xto{\mu} \uS). 
	\]
	We already described $\uS$ in \cref{example:FreeCpTambaraFixedGen}. 
	Using \cite[Definition 3.1]{HM2019}, we can compute $\uS\boxtimes \uS$ from this as
	\[
		(\uS\boxtimes \uS)(C_p/C_p) 
			= \z[t,n_x,x, n_y, y]/(t^2=pt,tx^p=tn_x, ty^p=tn_y),	\]
	\[
		(\uS\boxtimes \uS)(C_p/e) = \z[x, y].
	\]
	The kernel of the augmentation ideal $\uI$ can then be described as
	\[ 
		\uI(C_p/C_p)= \langle y-x, n_y-n_x\rangle, 
	\]
	\[ 
		\uI(C_p/e)= \langle y-x\rangle 
	\]
	where at each level, we described $\uI(C_p/H)$ as an ideal in $(\uS\boxtimes\uS)(C_p/H)$. 
	Moreover, we can calculate $\nm(y-x) \in n_y -n_x + \uI(C_p/C_p)$, and hence we can replace the generator $n_y-n_x$ above by $\nm(y-x)$. 
	Calculating now $\uI/\uI^2$ as an $\uS$-module gives $\uI/\uI^2\cong \uS\{y-x, \nm(y-x)\}/(\res(\nm(y-x))=0)$, where both generators live at level $C_p/C_p$.

	From this, by taking $\uA\boxtimes_{\uS} (\_)$ we obtain a description of the free genuine$\uA$ module on a generator at $C_p/C_p$, which we denote as $\uA\{y^\gen_{C_p/C_p} \}$. 
	The final description is:
	\[ 
		\uA\{y^\gen_{C_p/C_p} \}(C_p/C_p) 
			= \bbZ[t]/(t^2=pt) \{ y_{C_p/C_p} \} 
				\oplus \bbZ \{N(y_{C_p/C_p})\} 
	\]
	\[ 
		\uA\{y^\gen_{C_p/C_p} \}(C_p/e) 
			= \bbZ\{ R(y_{C_p/C_p})\}. 
	\] 
	Here, $\res(y_{C_p/C_p})=R(y_{C_p/C_p})$, $\tr(R(y_{C_p/C_p}))=ty_{C_p/C_p}$, $\nm(R(y_{C_p/C_p}))=N(y_{C_p/C_p})$ and $\res(N(y_{C_p/C_p}))=0$.
	Taking $\uA\boxtimes_{\uS} (\_)$ yields the desired description for $\uA\{y^\gen_{C_p/C_p} \}$.

	By a similar calculation, we can describe the free genuine $\uA$-module $\uA\{y^\gen_{C_p/e} \}$ on a generator at level $C_p/e$ as
	\[ 
		\uA\{y^\gen_{C_p/e} \}(C_p/C_p) = \bbZ \{ t(y_{C_p/e}), N(y_{C_p/e})\} 
	\]
	\[ 
		\uA\{y^\gen_{C_p/e} \}(C_p/e) = \bbZ\{ y_{C_p/C_p}, \gamma\cdot y_{C_p/e}, \ldots, \gamma^{p-1}\cdot y_{C_p/e} \}. 
	\] 
	Here, $\gamma$ is a generator of the Weyl group $C_p$ of $e$ in $C_p$, and $\cdot$ signifies the Weyl group action on $\uA\{y^\gen_{C_p/e} \}(C_p/e)$. 
	The restriction, transfer and norm are given by $\tr(\gamma^i\cdot y_{C_p/C_p})=t(y_{C_p/C_p})$, $\nm(\gamma^i\cdot y_{C_p/C_p})=N(y_{C_p/C_p})$, $\res(t(y_{C_p/e}))= \sum_{i=0}^{p-1} \gamma^i\cdot y_{C_p/e}$ and $\res(N(y_{C_p/C_p}))=0$.
\end{example}

Note that in these examples, there are indeed non-trivial norms on the generators, so these genuine modules cannot be described as naive modules.

\begin{example}
	\label{example:comparison_to_Stricklands_genuine_module}
	We also calculate $\ubbZ\{y^\gen_{C_p/C_p}\} = \ubbZ\boxtimes \uA\{y^\gen_{C_p/C_p}\}$, which is the free genuine $\ubbZ$-module on a fixed generator. 
	This is given as
	\[ 
		\ubbZ\{y^\gen_{C_p/C_p} \}(C_p/C_p) 
			= \bbZ \{ y_{C_p/C_p} \} 
				\oplus \bbZ/2 \{N(y_{C_p/C_p})\} 
	\]
	\[ 
		\ubbZ\{y^\gen_{C_p/C_p} \}(C_p/e) 
			= \bbZ\{ R(y_{C_p/C_p})\}. 
	\] 
	Here, $\res(y_{C_p/C_p})=R(y_{C_p/C_p})$, $\tr(R(y_{C_p/C_p}))=2y_{C_p/C_p}$, $\nm(R(y_{C_p/C_p}))=N(y_{C_p/C_p})$ and $\res(N(y_{C_p/C_p}))=0$.
\end{example}

We note that the genuine module in \Cref{example:comparison_to_Stricklands_genuine_module} is almost the genuine module given by Strickland in \cite[7.8 and 14.13]{Str2012}, up to the fact that a factor 2 was moved from the restriction to the transfer. 
Our formulation of this genuine module shows that it is a free genuine module. 
Moreover, in \cite[Section 1.3.a]{Stahlhauer22} it is shown that the first genuine module in \Cref{example:free_genuine_module} and that in \Cref{example:comparison_to_Stricklands_genuine_module} also extend to genuine modules over the Burnside-ring global power functor $\uA$ and the constant global power functor $\ubbZ$.

\begin{remark}
	As alluded to above, the HKR isomorphism does not hold if we work with genuine K{\"a}hler differentials but Hochschild homology only using the underlying Green functor of the input.
	We suspect that some version of Hochschild homology incorporating this genuine structure would be needed to obtain an HKR isomorphism involving genuine K{\"a}hler differentials, but have not been able to carry out this program yet. 
\end{remark}

\begin{remark}
	In \Cref{HochschildHomologyAsKaehlerDiffs} and \Cref{HKRnaiveKaehlerDiffs}, we identified the first Hochschild homology of an incomplete Tambara functor with the module of naive K{\"a}hler differentials as a genuine module. 
	It seems reasonable to expect an HKR theorem in this setting of genuine modules over incomplete Tambara functors as well. 
	To formulate such an HKR theorem, one would need a notion of exterior algebras of genuine modules. 
	Such an equivariant notion should also incorporate norms of generators in degree 1, which are then contained in higher degrees. An instructive example is the calculation in \Cref{Prop:TambaraUnderlying}, where we may interpret the final Tor Mackey functor $\uTor^{\uR}_p (\uA, \uA) \cong \ubbZ$ as generated on a norm element at level $C_p/C_p$ of the generator $x_{C_p/e}$ in degree 1. 
	This seems an interesting direction for future research.
\end{remark}

\section{Koszul resolutions for cyclic groups of prime order}
\label{SS:KoszulCp}

In \cite{MQS24b}, we presented some examples where $\uTor$ is badly behaved. 
In this paper, we study some situations where $\uTor$ is well-behaved. 
In \cref{SS:FreeGreenFixed}, we compute $\uTor$ over the free Green functor on a fixed generator over any finite group; this case is essentially classical homological algebra. 
In contrast, we introduce some new techniques to compute $\uTor$ over the free Tambara functor on an underlying generator over $C_p$ in \cref{SS:CpTambaraFixed}. 
Our approach is modeled on the classical Koszul resolution for a free $C_p$-algebra on one generator, but because of the presence of norms, we actually obtain a resolution by taking the total complex of a bicomplex of Koszul-type resolutions. 

\subsection{$\uTor$ over the free Green functor on a fixed generator}\label{SS:FreeGreenFixed}

Throughout this subsection, $G$ is an arbitrary finite group. 
The free Green functor on a fixed generator, $\uA^{\cO^\bot}[x_{G/G}]$, is particularly simple: we have
\[
	\uA^{\cO^\bot}[x_{G/G}](G/H) = \uA(G/H)[\res_H^G(x_{G/G})]
\]
for all $H \leq G$, i.e., it is levelwise a polynomial ring on one generator obtained by restricting the generator from the top level. 
Consequently, Mackey functor-valued Tor over this Green functor is particularly simple:

\begin{proposition}\label{Prop:GreenFixed}
We have
\[
	\uTor^{\uA^{\cO^\bot}[x_{G/G}]}_*(\uA, \uA) \cong 
		\begin{cases}
			\uA \quad & \text{ if } * =0, 1, \\
			0 \quad & \text{ if } *>1.
		\end{cases}
\]
\end{proposition}

\begin{proof}
	By \cite[Cor. 2.11]{BH2019}, there is an isomorphism of Green functors
	\[
		\uA^{\cO^\bot}[x_{G/G}] \cong \z[x] \otimes \uA.
	\]
	The free $\z[x]$-resolution of $\z$
	\[
		0 \to \z[x] \xrightarrow{\cdot x} \z[x] \xrightarrow{x \mapsto 0} \z
	\]
	can be promoted to a free $\uA^{\cO^\bot}[x_{G/G}]$-resolution of $\uA$
	\[
		0 \to \uA^{\cO^\bot}[x_{G/G}] \xrightarrow{ \cdot x} \uA^{\cO^\bot}[x_{G/G}] \xrightarrow{x \mapsto 0} \uA.
	\]
	Applying $\uA \boxtimes_{\uA^{\cO^\bot}[x_{G/G}]} -$ and computing homology gives the claimed answer. 
\end{proof}

\subsection{$\uTor$ over the free Tambara functor on an underlying generator}\label{SS:CpTambaraFixed}

Let $\uR := \uA^{\cO^{\top}}[x_{C_p/e}]$ be the free $C_p$-Tambara functor on an underlying generator, where $p$ is an odd prime. 
In this subsection, we compute $\uTor_*^{\uR}(\uA,\uA)$, where $\uA$ becomes an $\uR$ module with trivial action of $x_{C_p/e}$. We begin by recalling the concrete description of the Tambara functor $\uR$ obtained in \cref{thm:FreeIncompleteTambaraDescription}. 

\begin{proposition}
	\label{prop:freeCpTambaraUnderlyingGen}
	We have
	\[
		\uR(G/G) = \z[n,t_{\vec{v}}: \vec{v} \in \z_{\geq 0}^{\times p}]/(t_{\vec{0}}^2-pt_{\vec{0}}, t_{\vec{v}}-\gamma t_{\vec{v}}, t_{\vec{v}} t_{\vec{w}} - \sum_{\gamma \in C_p} t_{\vec{v}+\gamma \vec{w}}, n t_{\vec{v}} - t_{\vec{v}+\vec{1}}),
	\]
	\[
		\uR(G/e) = \z[x^{(i)} : 0 \leq i \leq p-1],
	\]
	where 
	\[
		\res_e^{C_p}(t_{\vec{v}}) = \sum_{\gamma \in C_p} x^{\gamma \vec{v}}, \quad \res_e^{C_p}(n) = \prod_{i=0}^{p-1} x^{(i)},
	\]
	\[
		\tr_e^{C_p}(x^{\vec{v}}) = t_{\vec{v}}.
	\]
	The norms are determined by 
	\[ 
		\nm_e^{C_p}(x^{(i)}) = n
	\]
	for all $i$ and the formula for the norm of a sum \cite[Theorem 3.5]{HM2019}. 
	The Weyl action of $W_{C_p}(e) = C_p$ on the underlying level is given by 
	\[ 
		\gamma \cdot x^{(i)} = x^{(i+1)},
	\]
	with the indices taken mod $p$ so that $\gamma \cdot x^{(p-1)} = x^{(0)}$. 
\end{proposition}

The $p = 2$ instance of the above proposition is \cite[Lemma 3.7]{BH2019}. 
The $p = 3$ instance of \cref{prop:freeCpTambaraUnderlyingGen} has already appeared as \cref{example:freeC3TambaraUnderlyingGen}. 

The goal of this section is to prove the following theorem. 

\begin{theorem}\label{Prop:TambaraUnderlying}
	Let $\uR=\uA^{\cO^\top}[x_{C_p/e}]$. We write $\uA_e := \uA\{x_{C_p/e}\}$ as shorthand for the free Mackey functor on an underlying generator. We have for $p\geq 3$ prime
	\[ 
		\uTor_k^{\uR}(\uA, \uA) \cong 
		\begin{cases}
			\uA
				& 
			\textup{ for } k=0,
				\\
			\uA_e \oplus \uI & \textup{ for } k=1,
				\\
			\bigoplus_{{p\choose k}/p} \uA_e 
				& 
			\textup{ for } 2\leq k\leq p-1,
				\\
			\ubbZ 
				& 
			\textup{ for } k=p,
		\end{cases}
	\]
	where $\uI$ is the kernel of the augmentation $\uA \to \ubbZ$. 
\end{theorem}

The definition of our free $\uR$-module resolution of $\uA$ requires two steps. 
Our first step is to lift the ordinary Koszul resolution of $\z[x^{(0)},\ldots,x^{(p-1)}]$ to a complex of $\uR$-modules (\cref{Constr:KoddpFirst}). 
This new complex will have some nontrivial homology because of the presence of norms (\cref{Lem:HKoddpFirst}). 
In the second step (\cref{Constr:KoddpSecond}), we modify this complex by taking an appropriate mapping cone to obtain the desired $\uR$-module resolution of $\uA$ (\cref{Lem:HKoddpSecond}). 

To make the computations more tangible, we will trace through the $p = 3$ example in this section.

\begin{construction}\label{Constr:KoddpFirst}
	Let $p$ be an odd prime. We construct a chain complex of Mackey functors 
	\[ 
		\begin{tikzcd}
			0 \ar[r] & 
			\uK_p \ar[r, "\partial_p"] & 
			\uK_{p-1} \ar[r, "\partial_{p-1}"] & 
			\cdots \ar[r, "\partial_2"] & 
			\uK_1 \ar[r, "\partial_1"] &
			\uK_0 = \uR
		\end{tikzcd}
	\]
	that lifts the Koszul resolution of $\uR(C_p/e) = \bbZ[x^{(0)}, x^{(1)}, \ldots, x^{(p-1)}]$ to a complex of $\uR$-modules in the sense that the underlying level of this chain complex is the Koszul complex resolving the underlying level of $\uR$. 

	Let $\uK_0= \uR\{y_{C_p/C_p}\}$. 
	For $1\leq k\leq p-1$, we define 
	\[ 
		\uK_k= \uR\{(z^{(i_1)}\wedge \ldots \wedge z^{(i_k)})_{C_p/e}\}_{I_k}, 
	\]
	a free $\uR$-module on generators at the underlying level indexed by alternating tensors of Weyl conjugates of $z$. 
	Here, $I_k$ is a set of representatives $\{i_1, \ldots, i_k\}$ of $k$-element subsets of the Weyl group $W_{C_p}(e)\cong C_p$ under the diagonal action of $W_{C_p}(e)$. 
	In particular, the cardinality of $I_k$ is ${p\choose k}/p$, and we might take for example 
	\[ 
		I_1=\{\{0\}\},
		\quad 
		I_2 = \{ \{0,1\}, \ldots, \{0, (p-1)/2\} \} 
		\quad \textup{and} \quad 
		I_{p-1} = \{ \{0, \ldots, p-2\} \}. 
	\]
	Note that by the definition of alternating tensors and by the Weyl action, an $\uR(C_p/e)$-basis of $\uK_k(C_p/e)$ is given by alternating tensors $z^{(i_1)}\wedge \ldots \wedge z^{(i_k)}$ for $i_1\leq \ldots \leq i_k$ and \emph{all} $k$-element subsets $\{i_1, \ldots, i_k\} \subset W_{C_p}(e)$. 
	Using this, we define the usual Koszul differential
	\[ 
		\partial_k\colon \uK_k\to \uK_{k-1}, \quad 
			z^{(i_1)}\wedge \ldots \wedge z^{(i_k)} 
			\mapsto 
			\sum_{j=1}^k (-1)^{j-1} x^{(i_j)} 
				\cdot z^{(i_1)}\wedge 
					\ldots \wedge \widehat{z^{(i_j)}} \wedge 						\ldots \wedge z^{(i_k)}.
	\]
	Finally, we define $\uK_p= \uR\{ N(z^{(0)})_{C_p/C_p} \}$, a free $\uR$-module generated at level $C_p/C_p$ by one formal element $N(z^{(0)})$. 
	We denote its restriction to $C_p/e$ by $z^{(0)}\wedge \ldots \wedge z^{(p-1)}$. 
	On this element, we define the differential as 
	\[ 
		\partial_p( N(z^{(0)}))= x^{(p-1)} \tr (z^{(0)}\wedge \ldots \wedge z^{(p-2)}).
	\]
	Restricting this differential to $C_p/e$, we obtain the usual Koszul differential description for $\partial_p(z^{(0)}\wedge \ldots \wedge z^{(p-1)})$.\footnote{This is where we need to assume $p$ is odd: commuting wedge factors introduces the sign in the Koszul differential.}  
	For $k>p$, we set $\uK_k=0$.
\end{construction}

\begin{remark}
	\label{remark:KoszulDiffOfNorm}
	The differential of the generator $N(z^{(0)})$ defined in the above construction might seem surprising at a first glance, but its formula can be justified by using the formula for the value of a (genuine) derivation on a norm in a Tambara functor as described in \Cref{DefDerivations}(c). 
	In fact, the Koszul differential on the element $N(z^{(0)})$ is of the form described by this formula, since the $C_p$-set $C_p/e \times C_p/e \setminus \Delta$ decomposes into $p-1$ copies of $C_p/e$. 
	The composition of restriction and norm thus is a simple product, and represented here as the wedge product $z^{(0)}\wedge \ldots \wedge z^{(p-2)}$. 
	This term is then multiplied with the differential of $z^{(p-1)}$ and transferred up to $C_p/C_p$ to obtain the formula above.
\end{remark}

\begin{example}
	\label{spiderExample} 
	For $p = 3$, this chain complex is the following: 
	\begin{equation}
	\label{smallKoszulC3}
		\begin{tikzcd}
			\uR\{(N(z^{(0)})_{C_3/C_3}\} 
				\ar[r, "\partial_3"] 
				&
			\uR\{(z^{(0)} \wedge z^{(1)})_{C_3/e}\}
				\ar[r, "\partial_2"] 
				& 
			\uR\{(z^{(0)})_{C_3/e}\}
				\ar[r, "\partial_1"]
				& 
			\uR.
		\end{tikzcd}
	\end{equation}
	Recall the explicit description of $\uR = \uA^{\cO^\top}[x_{C_3/e}]$ from \cref{example:freeC3TambaraUnderlyingGen}. 
	The $\uR$-modules in the above chain complex are the following: 
	\[ 
		\uR\{(z^{(0)})_{C_3/e}\} = 
		\begin{tikzcd}
			\uR(C_3/e)\{\tr(z)\}
				\ar[d, bend right=50, "\res"{left}]
			\\
			\uR(C_3/e)\{z^{(0)},z^{(1)},z^{(2)}\}
				\ar[u, bend right=50, "\tr"{right}]	
				\arrow[out=240,in=300,loop,swap,looseness=4, "C_3"]
		\end{tikzcd}
	\]
	The transfer in this $\uR$-module sends $z^{(i)}$ to $\tr(z)$ and the Weyl action permutes the module generators. Restriction sends $\tr(z)$ to the sum of its Weyl conjugates: $\res(\tr(z)) = z^{(0)} + z^{(1)} + z^{(2)}$. 
	\[ 
		\uR\{(z^{(0)} \wedge z^{(1)})_{C_3/e}\} = 
		\begin{tikzcd}
			\uR(C_3/e)\{\tr(z^{(0)} \wedge z^{(1)})\}
				\ar[d, bend right=50, "\res"{left}]
			\\
			\uR(C_3/e)\left\{
					z^{(0)} \wedge z^{(1)},
					z^{(1)} \wedge z^{(2)},
					z^{(2)} \wedge z^{(0)}
				\right\}
				\ar[u, bend right=50, "\tr"{right}]
				\arrow[out=240,in=300,loop,swap,looseness=4, "C_3"]
		\end{tikzcd}
	\]
	The transfer in this $\uR$-module sends $z^{(i)}$ to $\tr(z)$ and the Weyl action permutes the module generators by distributing over $\wedge$:  
	\[
		\gamma \cdot z^{(0)} \wedge z^{(1)} = z^{(1)} \wedge z^{(2)} \qquad
		\gamma \cdot z^{(1)} \wedge z^{(2)} = z^{(2)} \wedge z^{(0)} \qquad
		\gamma \cdot z^{(2)} \wedge z^{(0)} = z^{(0)} \wedge z^{(1)} 
	\]
	Restriction sends $\tr(z^{(0)} \wedge z^{(1)})$ to the sum of its Weyl conjugates. 
	
	The $\uR$-module $\uR\{N(z^{(0)})_{C_3/C_3}\}$ is isomorphic to $\uR$; we have simply named the generator at the top level. 
	
	The differential $\partial_3$ is determined by where it sends the generator on the top level. 
	To describe its behavior on the underlying level, we take a restriction of $x^{(2)} \tr(z^{(0)} \wedge z^{(1)})$ by summing over Weyl conjugates of $x^{(2)} z^{(0)} \wedge z^{(1)}$. 
	We have:
	\begin{align}
		\partial_3(C_3/C_3) &\colon F \cdot 
			N(z^{(0)}) 
			\mapsto 
			\res(F) x^{(2)}\tr(z^{(0)} \wedge z^{(1)}),
			\label{differential3top}
		\\
		\partial_3(C_3/e) &\colon 
			\res\left( N(z^{(0)})\right)
			\mapsto 
			x^{(2)} z^{(0)} \wedge z^{(1)} 
				+ x^{(0)}z^{(1)} \wedge z^{(2)} 
				+ x^{(1)} z^{(2)} \wedge z^{(0)}
				\label{differential3bottom}
	\end{align}
	
	The differential $\partial_2$ is determined on the underlying level. 
	To find what it does to $\tr(z^{(0)} \wedge z^{(1)})$ on the top level, we simply take a transfer of $x^{(0)} z^{(1)} - x^{(1)} z^{(0)}$. 
	Since the transfer is equivariant for the Weyl actions, we can rewrite 
	\[
		\tr(x^{(0)} z^{(1)} - x^{(1)} z^{(0)}) 
			= \tr(x^{(0)} z^{(1)}) - \tr(x^{(1)} z^{(0)}) 
			= \tr(x^{(2)} z^{(0)}) - \tr(x^{(1)} z^{(0)}) 
			= (x^{(2)} - x^{(1)}) \tr(z^{(0)}),
	\]
	so the second differential in this complex is determined by: 
	\begin{align}
		\partial_2(C_3/C_3) &\colon 
			\tr(z^{(0)} \wedge z^{(1)}) 
			\mapsto 
			(x^{(2)} - x^{(1)}) \tr(z^{(0)}),
			\label{differential2top}
		\\
		\partial_2(C_3/e) &\colon 
			z^{(i)} \wedge z^{(j)} 
			\mapsto 
			x^{(i)} z^{(j)} - x^{(j)} z^{(i)} 
			\label{differential2bottom}
	\end{align}
	
	Finally, $\partial_1$ is determined by what it does on the underlying level, and on the top level by transferring from the underlying:
	\begin{align*}
		\partial_1(C_3/C_3) &\colon \tr(z^{(0)}) \mapsto t_{100}, \\
		\partial_1(C_3/e) &\colon z^{(i)} \mapsto x^{(i)}
	\end{align*}
	
	The fact that this is a complex at the underlying level is straightforward, so we explain why this defines a complex at the fixed level. 
	In degree $2$, the effect of $\partial_2$ on $x^{(2)} \tr(z^{(0)} \wedge z^{(1)})$ is determined by what happens on the underlying level:
	\[
		\begin{tikzcd}[column sep=huge]
			x^{(2)} z^{(0)} \wedge z^{(1)} 
			\ar[r, mapsto, "\partial_2(C_3/e)"]
			&
			x^{(2)} x^{(0)} z^{(1)} - x^{(2)} x^{(1)} z^{(0)}
		\end{tikzcd}
	\]
	After transferring back up, we get  
	\[
		\tr(x^{(2)} x^{(0)} z^{(1)} - x^{(2)} x^{(1)} z^{(0)}) = 
		\tr(x^{(2)} x^{(0)} z^{(1)}) - \tr(x^{(2)} x^{(1)} z^{(0)}) = 
		\tr(x^{(1)} x^{(2)} z^{(0)}) - \tr(x^{(2)} x^{(1)} z^{(0)}) = 
		0
	\]
	using Weyl invariance of transfers. 
	
	In degree $1$, something similar happens: 
	\[ 
		\begin{tikzcd}[column sep=huge]
			x^{(2)} z^{(0)} - x^{(1)} z^{(0)} 
				\ar[r, mapsto, "\partial_1(C_3/e)"] 
			& 
			x^{(2)} x^{(0)} - x^{(1)} x^{(0)},
		\end{tikzcd}
	\]
	so after transferring, we get 
	\[ 
		\begin{tikzcd}[column sep=huge]
			(x^{(2)} - x^{(1)}) \tr(z^{(0)}) 
				\ar[r, mapsto, "\partial_1(C_3/C_3)"] 
			& 
			\tr(x^{(2)} x^{(0)}) - \tr(x^{(1)} x^{(0)}) 
			= 
			t_{101} - t_{110} 
			= 0,
		\end{tikzcd}
	\]	
	remembering that transfers are Weyl invariant and $x^{(2)} x^{(0)} = \gamma^2 \cdot x^{(1)} x^{(0)}$. 
	
\end{example}

\begin{lemma}
	\label{Lem:HKoddpFirst}
	With $\uK_\bullet$ as in \cref{Constr:KoddpFirst}, we have
	\[ 
		H_\ast(\uK_\bullet) = 
		\begin{cases}
			\uA\oplus \uI\{n^i_{C_p} \mid i\geq 1\} 
				& 
			\textup{ for } \ast =0,
				\\
			0 
				& 
			\textup{ for } 1\leq \ast\leq p-1,
				\\
			\uI_{\uR}\{t-p\} 
				& 
			\textup{ for } \ast= p.\\
		\end{cases} 
	\]
\end{lemma}

Recall $t$ is the class of the finite $C_p$-set $C_p/e$ in the top level of the Burnside Tambara functor $\uA$ which injects into $\uR = \uA^{\cO^\top}[x_{C_3/e}]$. 
	Then $\langle t - p \rangle$ is the kernel of the restriction from the top level to the underlying level in $\uA$, and similarly in $\uR$. 
	Here, we consider the ideal $\uI \subseteq \uA$ generated by the class $t - p$ in degree $p$, and the corresponding ideal $\uI_{\uR} \subseteq \uR$.

\begin{proof}
	At level $C_p/e$, we observe that $\uK_\bullet$ is a classical Koszul complex and hence a resolution of $\mathbb Z$. At level $C_p/C_p$, we observe that the image of $\partial_1\colon \uK_1\to \uK_0$ is exactly the image of the transfer of the ideal in $\uR(C_p/e)$ generated by $x$ and all its Weyl conjugates. 
	A straightforward calculation shows that the complex $\uK_\bullet$ is exact for $1\leq k\leq p-2$.

	We now describe the final map $\partial_{p}\colon \uK_{p}\to \uK_{p-1}$. 
	We may describe the differential $\partial_p$ on general elements of $\uK_p(C_p/C_p)$ as
	\[ 
		\partial_p (f\cdot N(z^{(0)})) 
		= (\res(f)\cdot x^{(p-1)})\tr (z^{(0)}\wedge \ldots \wedge z^{(p-2)}) 
	\] 
by $\uR$-linearity. 
Moreover, a straight-forward calculation shows that the kernel of $\partial_{p-1}$ at level $C_p/C_p$ is of the form
\[ 
	\{ F\cdot x^{(p-1)} \tr(z^{(0)}\wedge \ldots \wedge z^{(p-2)}) \mid F\in\uR(C_p/e)^{W_{C_p}(e)} \}. 
\]
We observe that the map $\res\colon \uR(C_p/C_p)\to \uR(C_p/e)^{W_{C_p}e}$ is surjective.
In particular, the norm element $n=\nm_e^{C_p}(x)$ maps to $x^{(0)}\ldots x^{(p-1)}$. 
Hence, the Koszul complex is also exact at $p-1$. 
Finally, since the element $x^{(p-1)} \tr(z^{(0)}\wedge \ldots \wedge z^{(p-2)})\in \uK_{p-1}(C_p/e)$ is annihilated only by $0$, we observe that $H_p(\uK_\bullet)(C_p/C_p)=\ker(\res\colon \uR(C_p/C_p)\to \uR(C_p/e))\cong \uI_{\uR}$ by \cref{lemma:kernels}.
\end{proof}

\begin{example}
	\label{spiderExampleCohomology}
	Continuing the $C_3$ example from \cref{spiderExample}, we compute the homology of the complex of \cref{smallKoszulC3}. Recall $\uR = \uA^{\cO^\top}[x_{C_3/e}]$ from \cref{example:freeC3TambaraUnderlyingGen}. 
	
	In degree zero, we get the cokernel of the first differential. 
	\[  
		\begin{tikzcd}
			\uR(C_3/e)\{\tr(z^{(0)})\}
				\ar[d, bend right=50, "\res"{left}]
				\ar[r, "\tr(z) \mapsto t_{100}"]
			& 
			\uR(C_3/C_3)
				\ar[d, bend right=50, "\res"{left}]
				\ar[r]
			& 
			\bbZ[n, t_{000}]/(t_{000}^2 - 3t_{000}, t_{000} n)
				\ar[d, bend right=50, 
					"t_{000} \mapsto 3"{left,near start}, 
					"n \mapsto 0"{left, near end}]
			\\
			\uR(C_3/e)\{z^{(0)},z^{(1)},z^{(2)}\}
				\ar[u, bend right=50, "\tr"{right}]	
				\arrow[out=240,in=300,loop,swap,looseness=4, "C_3"]
				\ar[r, "z^{(i)} \mapsto x^{(i)}"]
			& 
			\uR(C_3/e)
				\ar[u, bend right=50, "\tr"{right}]	
				\arrow[out=240,in=300,loop,swap,looseness=4, "C_3"]
				\ar[r]
			&
			\bbZ
				\ar[u, bend right=50, "\cdot t_{000}"{right}]	
			\\
			\\
			\uR\{z^{(0)}_{C_3/e}\}
				\ar[r, "\partial_1"]
			& 
			\uR 
				\ar[r]
			& 
			\coker(\partial_1)
		\end{tikzcd}
	\]
	To understand the image of $\partial_1$ on the top level, note that every element of $\uR\{z^{(0)}_{C_3/e}\}(C_3/C_3)$ is a transfer $f \tr(z) = \tr(fz^{(0)})$ for $f \in \uR(C_3/e) = \bbZ[x^{(0)}, x^{(1)}, x^{(2)}]$. 
	Therefore, an element $f \tr(z) = \tr(fz^{(0)}) \in \uR\{z^{(0)}_{C_3/e}\}$ is sent to $\tr(\partial_1(fz^{(0)})) = \tr(x^{(0)} f) \in \uR(C_3/C_3)$. 
	In particular, since every polynomial in the $x^{(i)}$ is hit on the underlying level, we hit every $t_{ijk}$ at the top level, so long as $i,j,k$ are not all zero. 
	This leaves us with the Mackey functor displayed above, which is a sum of a copy of the Burnside functor $\uA$, generated by $1$ at the top level, and one copy of the augmentation ideal $\uI$ for each nonzero power of $n$. 
	
	In degree one, the underlying level is a classical Koszul complex, so it has zero homology. 
	The fixed level is entirely determined by the underlying one: the kernel of $\partial_1$ at the fixed level is $(x^{(i)} - x^{(j)}) \tr(z^{(0)})$, which is in the image of $\partial_2$ at that level.
	\[  
		\hspace*{-2cm}
		\begin{tikzcd}[column sep=huge]
			\uR(C_3/e)\{\tr(z^{(0)} \wedge z^{(1)})\}
				\ar[d, bend right=50, "\res"{left}]
				\ar[r, "\cref{differential2top}"{below}, 
					"\tr(z^{(0)} \wedge z^{(1)}) \mapsto (x^{(2)} - x^{(1)}) \tr(z^{(0)})"]
			&
			\uR(C_3/e)\{\tr(z^{(0)})\}
				\ar[d, bend right=50, "\res"{left}]
				\ar[r, "\tr(z) \mapsto t_{100}"]
			& 
			\uR(C_3/C_3)
				\ar[d, bend right=50, "\res"{left}]
			\\
			\uR(C_3/e)\left\{
					z^{(0)} \wedge z^{(1)},
					z^{(1)} \wedge z^{(2)},
					z^{(2)} \wedge z^{(0)}
				\right\}
				\ar[u, bend right=50, "\tr"{right}]
				\arrow[out=240,in=300,loop,swap,looseness=4, "C_3"]
				\ar[r, "\cref{differential2bottom}"]
			&
			\uR(C_3/e)\{z^{(0)},z^{(1)},z^{(2)}\}
				\ar[u, bend right=50, "\tr"{right}]	
				\arrow[out=240,in=300,loop,swap,looseness=4, "C_3"]
				\ar[r, "z^{(i)} \mapsto x^{(i)}"]
			& 
			\uR(C_3/e)
				\ar[u, bend right=50, "\tr"{right}]	
				\arrow[out=240,in=300,loop,swap,looseness=4, "C_3"]
			\\
			\\
			\uR\{(z^{(0)} \wedge z^{(1)})_{C_3/e}\} 
				\ar[r, "\partial_2"]
			&
			\uR\{z^{(0)}_{C_3/e}\}
				\ar[r, "\partial_1"]
			& 
			\uR 
		\end{tikzcd}
	\]

	In degree two, we have: 
	\[	
		\hspace*{-3cm}
		\begin{tikzcd}[column sep=7em]
			\uR(C_3/C_3)\left\{N(z^{(0)})\right\}
				\ar[d, bend right=50, "\res"{left}]
				\ar[r, "G\cdot N(z^{(0)}) \mapsto \res(G)x^{(2)} \tr(z^{(0)} \wedge z^{(1)})", "\cref{differential3top}"{below}]
			& 
			\uR(C_3/e)\{\tr(z^{(0)} \wedge z^{(1)})\}
				\ar[d, bend right=50, "\res"{left}]
				\ar[r, "\cref{differential2top}"{below}, "\tr(z^{(0)} \wedge z^{(1)}) \mapsto (x^{(2)} - x^{(1)}) \tr(z^{(0)})"]
			& 
			\uR(C_3/e)\{\tr(z^{(0)})\}
				\ar[d, bend right=50, "\res"{left}]
			\\
			\uR(C_3/e)\{\res(N(z^{(0)}))\}
				\ar[u, bend right=50, "\tr"{right}]
				\ar[r, "\cref{differential3bottom}"]
				\arrow[out=240,in=300,loop,swap,looseness=4, "C_3"]
			&
			\uR(C_3/e)\left\{
					z^{(0)} \wedge z^{(1)},
					z^{(1)} \wedge z^{(2)},
					z^{(2)} \wedge z^{(0)}
				\right\}
				\ar[u, bend right=50, "\tr"{right}]
				\arrow[out=240,in=300,loop,swap,looseness=4, "C_3"]
				\ar[r, "\cref{differential2bottom}"]
			& 
			\uR(C_3/e)\{z^{(0)}, z^{(1)}, z^{(2)}\}
				\ar[u, bend right=50, "\tr"{right}]
				\arrow[out=240,in=300,loop,swap,looseness=4, "C_3"]
			\\
			\\
			\uR\big\{\left(N(z^{(0)})\right)_{C_3/e} \big\}
				\ar[r, "\partial_3"]
			&
			\uR\{(z^{(0)} \wedge z^{(1)})_{C_3/e}\} 
				\ar[r, "\partial_2"]
			& 
			\uR\{z^{(0)}_{C_3/e}\}
		\end{tikzcd}
	\]
	On the underlying level, we have a classical Koszul complex with zero homology. 
	At the top level, the kernel of $\partial_2$ is elements of the form $F x^{(2)} \tr(z^{(0)} \wedge z^{(1)})$, where $F$ is a Weyl-fixed element of $\uR(C_3/e)$. 
	Such an element is necessarily in the image of $\partial_3$ because any restriction is Weyl-fixed. 	
	
	In degree three, we have the kernel of the differential $\partial_3$:
	\[ 
		\begin{tikzcd}
			\langle t_{000} - 3 \rangle
				\ar[d, bend right=50]
				\ar[r]
			&
			\uR(C_3/C_3)\left\{N(z^{(0)})\right\}
				\ar[d, bend right=50, "\res"{left}]
				\ar[r, "\cref{differential3top}"]
			& 
			\uR(C_3/e)\{\tr(z^{(0)} \wedge z^{(1)})\}
				\ar[d, bend right=50, "\res"{left}]
			\\
			0 
				\ar[u, bend right=50]
				\ar[r] 
			& 
			\uR(C_3/e)\{\res(N(z^{(0)}))\}
				\ar[u, bend right=50, "\tr"{right}]
				\ar[r, "\cref{differential3bottom}"]
				\arrow[out=240,in=300,loop,swap,looseness=4, "C_3"]
			&
			\uR(C_3/e)\left\{
					z^{(0)} \wedge z^{(1)},
					z^{(1)} \wedge z^{(2)},
					z^{(2)} \wedge z^{(0)}
				\right\}
				\ar[u, bend right=50, "\tr"{right}]
				\arrow[out=240,in=300,loop,swap,looseness=4, "C_3"]
			\\
			\\
			\ker(\partial_3)
				\ar[r]
			& 
			\uR\big\{\left(N(z^{(0)})\right)_{C_3/C_3} \big\}
				\ar[r, "\partial_3"]
			&
			\uR\{(z^{(0)} \wedge z^{(1)})_{C_3/e}\} 
		\end{tikzcd}
	\]
	By \cref{differential3bottom}, we see that the differential $\partial_3$ is injective on the underlying level. To understand the kernel on the top level, recall from \cref{differential3top} that 
	\[
		\partial_3(C_3/C_3)(f N(z^{(0)})) = \res(f)x^{(2)} \tr(z^{(0)} \wedge z^{(1)})
	\]
	for any $f \in \uR(C_3/C_3)$. In particular, as in \cref{lemma:kernels}, the kernel at this level is the kernel of restriction, which is generated by $t_{000} - 3$ (recall that $t_{000} = t \in \uA(C_3/C_3) \subseteq \uR(C_3/C_3)$). 
\end{example}

With that example out of the way, we see that the complex $\uK_{\bullet}$ is not quite a resolution of $\uA$; the degree zero homology contains generators $n^i$ for $i \geq 1$, and the top homology does not vanish. 
In order to obtain a resolution for $\uA$, we need to kill the elements $n^i$ in the 0-th homology. 
We do this by forming a mapping cone for a map of complexes describing multiplication by $n$ in the 0-th component:

\begin{construction}
	\label{Constr:KoddpSecond}
	We define an ``$n$-divided Koszul complex'' $\ndivkoszul$ by $\ndivkoszul_0= \uR\{u_{C_p}\}$, 
	\[ 
		\ndivkoszul_{k} 
		= \uR\left\{ \left(\frac{z^{(i_1)}\wedge \ldots\wedge z^{(i_k)} \wedge R(u)}{x^{(i_1)}\cdot \ldots \cdot x^{(i_k)}} \right)_e \right\}_{I_k} 
	\]
for $1\leq k\leq p-1 $, and $\ndivkoszul_p= \uR\{ (N(z^{(0)}) \wedge u/\nm_e^{C_p}(x^{(0)}) )_{C_p} \}$.\footnote{The generating elements are treated as formal elements, not as actual quotients. Also here, $I_k$ is the same set of representatives of $k$-element subsets of $\{0, \ldots, p-1\}$ under the Weyl group action.} 
	The differentials on this complex are given by
	\[ 
		\ndivpartial_k \left(
			\frac{z^{(i_1)}\wedge \ldots\wedge z^{(i_k)} \wedge R(u)}{x^{(i_1)}\cdot \ldots \cdot x^{(i_k)}}
			\right)
			= 
			\sum_{j=1}^k (-1)^{j-1} 
				\frac{z^{(i_1)}\wedge \ldots \wedge \widehat{z^{(i_j)}} \wedge \ldots \wedge z^{(i_k)} \wedge R(u)}{x^{(i_1)}\cdot \ldots \cdot \widehat{x^{(i_j)}}\cdot\ldots \cdot x^{(i_k)}} 
	\]
	and 
	\[ 
	\ndivpartial_p \left(
		\frac{N(z^{(0)}) \wedge u}{\nm_e^{C_p}(x^{(0)})} 
		\right)
		= 
		\tr_e^{C_p} \left( 
			\frac{z^{(0)}\wedge \ldots\wedge z^{(p-2)} \wedge R(u)}{x^{(0)}\cdot \ldots \cdot x^{(p-2)}}
			\right). 
		\]
\end{construction}

This again is a Koszul complex with an added formal variable $u$, where instead of taking the (regular) sequence $(x^{(0)}, \ldots, x^{(p-1)})$ for the differentials, we instead take $(1, \ldots, 1)$. 

\begin{example}
	\label{C3dividedKoszul}
	Continuing \cref{spiderExampleCohomology}. 
	The $n$-divided Koszul complex for $C_3$ is 
	\[	
	\begin{tikzcd}
		\uR\{\left(\frac{N(z^{(0)}) \wedge u}{n}\right)_{C_3/C_3}\}
			\ar[r, "\ndivpartial_3"]
		& 
		\uR\{\left(\frac{z^{(0)} \wedge z^{(1)} \wedge \res(u)}{x^{(0)}x^{(1)}}\right)_{C_3/e}\}
			\ar[r, "\ndivpartial_2"]
		&
		\uR\{\left(\frac{z^{(0)} \wedge \res(u)}{x^{(0)}}\right)_{C_3/e}\}
			\ar[r, "\ndivpartial_1"]
		&
		\uR\{u_{C_3/C_3}\}
	\end{tikzcd}
	\]
	At each level, there is a single $\uR$-module generator; at level $C_3/C_3$ in degrees 0 and 3, and at level $C_3/e$ in degrees 1 and 2. 
	We write 
	\[ 
		\frac{z^{(1)} \wedge \res(u)}{x^{(1)}} 
		\quad \text{ and }\quad  
		\frac{z^{(2)} \wedge \res(u)}{x^{(2)}}
	\]
	for the Weyl conjugates of $\frac{z^{(0)} \wedge \res(u)}{x^{(0)}}$ in degree $1$, and 
	\[
		\frac{z^{(1)} \wedge z^{(2)} \wedge \res(u)}{x^{(1)}x^{(2)}}
		\quad \text{ and } \quad
		\frac{z^{(2)} \wedge z^{(0)} \wedge \res(u)}{x^{(2)}x^{(0)}}
	\]
	for the Weyl conjugates of $\frac{z^{(0)} \wedge z^{(1)} \wedge \res(u)}{x^{(0)}x^{(1)}}$ in degree $2$. 

	The differentials are:
	\[
	\begin{tikzcd}
		\displaystyle \frac{z^{(0)} \wedge \res(u)}{x^{(0)}} 
			\ar[r, mapsto, "\ndivpartial_1"] 
		& 
		\res(u) 
		\\
		\displaystyle \frac{z^{(0)} \wedge z^{(1)} \wedge \res(u)}{x^{(0)}x^{(1)}} 
			\ar[r, mapsto, "\ndivpartial_2"] 
		& 
		\displaystyle 
		\frac{z^{(1)} \wedge \res(u)}{x^{(1)}} - \frac{z^{(0)}\wedge \res(u)}{x^{(0)}}
		\\
		\displaystyle
		\frac{N(z^{(0)}) \wedge u}{n}
			\ar[r, mapsto, "\ndivpartial_3"]
		& 
		\displaystyle
		\tr\left(\frac{z^{(0)} \wedge z^{(1)} \wedge \res(u)}{x^{(0)}x^{(1)}}\right)
	\end{tikzcd}
	\]
\end{example}

A similar calculation to the one in \cref{Lem:HKoddpFirst} for the complex $\uK_\bullet$ yields: 

\begin{lemma}
\[ H_\ast(\ndivkoszul_\bullet) = \begin{cases}
\uI\{n^i_{C_p} \mid i\geq 0\} & \textup{ for } \ast =0,\\
0 & \textup{ for } 1\leq \ast\leq p-1,\\
\uI_{\uR}\{t-p\} & \textup{ for } \ast= p.\\
\end{cases} \]
\end{lemma}

\begin{example}
Comparing the homology calculation in \cref{spiderExampleCohomology} and the complex in \cref{C3dividedKoszul}, we see that the only real difference in the homology is in degree zero, where the differential is now surjective on the underlying level. Therefore, every transfer is hit by $\partial_1$ at the fixed level, and we are left with $\bbZ[n]$. As an $\uR$-module, this is a copy of $\uI$ for each power of $n$. 
\end{example}

\begin{construction}\label{Constr:KoddpMap}
Let $f\colon \ndivkoszul_\bullet \to \uK_\bullet$ be the morphism induced by $u\mapsto n$. On generators in higher degrees, we define
\[ f\left(\frac{z^{(i_1)}\wedge \ldots\wedge z^{(i_k)} \wedge R(u)}{x^{(i_1)}\cdot \ldots \cdot x^{(i_k)}} \right) = \frac{\res(n)}{x^{(i_1)}\cdot \ldots \cdot x^{(i_k)}} \cdot z^{(i_1)}\wedge \ldots\wedge z^{(i_k)},\]
and
\[ f\left(\frac{N(z^{(0)}) \wedge u}{\nm_e^{C_p}(x^{(0)})} \right) = N(z^{(0)}).\]
\end{construction}

\begin{example}
	\label{C3map}
	Continuing the \cref{C3dividedKoszul,spiderExample}, the morphism $f \colon \ndivkoszul \to \uK_\bullet$ induced by $u \mapsto n$ is displayed below. Vertical arrows are labelled with the image of the generator.
	\[
	\hspace*{-2cm}
	\begin{tikzcd}[row sep=huge]
		\ndivkoszul_\bullet 
			\ar[d, "f"] 
		& 
		\uR\{(\tfrac{N(z^{(0)}) \wedge u}{n})_{C_3/C_3}\}
			\ar[r, "\ndivpartial_3"]
			\ar[d, "N(z^{(0)})"]
		& 
		\uR\{(\tfrac{z^{(0)} \wedge z^{(1)} \wedge \res(u)}{x^{(0)}x^{(1)}})_{C_3/e}\}
			\ar[r, "\ndivpartial_2"]
			\ar[d, "x^{(2)} z^{(0)} \wedge z^{(1)}"]
		&
		\uR\{(\tfrac{z^{(0)} \wedge \res(u)}{x^{(0)}})_{C_3/e}\}
			\ar[r, "\ndivpartial_1"]
			\ar[d, "x^{(1)}x^{(2)} z^{(0)}"]
		&
		\uR\{u_{C_3/C_3}\}
			\ar[d, "u \mapsto n"] 
		\\
		\uK_\bullet
			& 
		\uR\{(N(z^{(0)}))_{C_3/C_3}\} 
			\ar[r, "\partial_3"] 
			&
		\uR\{(z^{(0)} \wedge z^{(1)})_{C_3/e}\}
			\ar[r, "\partial_2"] 
			& 
		\uR\{(z^{(0)})_{C_3/e}\}
			\ar[r, "\partial_1"]
			& 
		\uR.
	\end{tikzcd}
	\]
\end{example}

\begin{lemma}\label{Lem:HKoddpSecond}
The map $f: \ndivkoszul_\bullet \to \uK_{\bullet}$ defined in \cref{Constr:KoddpMap} is a map of complexes, and on homology it induces multiplication by $n$ in degree 0 and the identity in degree $p$. Hence, the mapping cone $\uF_\bullet=\Cone(f) \to \uA$ is a resolution of $\uA$ by free $\uR$-modules. 
\end{lemma}

\begin{example}
	The mapping cone of \cref{C3map} yields the following resolution of $\uA$ by free $\uR$-modules: 
	\[ 
		\hspace*{-2cm}
		\begin{tikzcd}[ampersand replacement=\&]
			\uR\{(\tfrac{N(z^{(0)}) \wedge u}{n})_{C_3/C_3}\}
				\ar[r] 
			\&
			\begin{matrix}
				\uR\{(N(z^{(0)})_{C_3/C_3}\}\\
				\bigoplus\\
				\uR\{(\tfrac{z^{(01)}\wedge \res(u)}{x^{(0)}x^{(1)}})_{C_3/e}\}	
			\end{matrix}
				\ar[r]
			\& 
			\begin{matrix}
				\uR\{(z^{(01)})_{C_3/e}\} \\
				\bigoplus \\
				\uR\{(\tfrac{z^{(0)} \wedge \res(u)}{x^{(0)}})_{C_3/e}\}
			\end{matrix}
				\ar[r]
			\&
			\begin{matrix}
				\uR\{(z^{(0)})_{C_3/e}\}\\
				\bigoplus\\
				\uR\{u_{C_3/C_3}\}
			\end{matrix}
				\ar[r]
			\&
			\uR 
				\ar[d] 
			\\
			\&
			\&
			\&
			\&
			\uA
		\end{tikzcd}
	\]
	Above, we have adopted the shorthand notation $z^{(ij)} = z^{(i)} \wedge z^{(j)}$. We can use this to calculate $\uTor_*^{\uR}(\uA,\uA)$. Taking a box product with $\uA$ kills each of the $x^{(i)}$. So the differentials in $\uK_\bullet$ vanish and we are left with the complex below. The effect of the differentials on generators is displayed below the complex.
	\[ 
		\hspace*{-2cm}
		\begin{tikzcd}[ampersand replacement=\&]
			\uA\{(\tfrac{N(z^{(0)}) \wedge u}{n})_{C_3/C_3}\}
				\ar[r] 
			\&
			\begin{matrix}
				\uA\{(N(z^{(0)}))_{C_3/C_3}\}\\
				\bigoplus\\
				\uA\{(\tfrac{z^{(01)}\wedge \res(u)}{x^{(0)}x^{(1)}})_{C_3/e}\}	
			\end{matrix}
				\ar[r]
			\& 
			\begin{matrix}
				\uA\{(z^{(01)})_{C_3/e}\} \\
				\bigoplus \\
				\uA\{(\tfrac{z^{(0)} \wedge \res(u)}{x^{(0)}})_{C_3/e}\}
			\end{matrix}
				\ar[r]
			\&
			\begin{matrix}
				\uA\{(z^{(0)})_{C_3/e}\}\\
				\bigoplus\\
				\uA\{u_{C_3/C_3}\}
			\end{matrix}
				\ar[r, "0"]
			\&
			\uA 
			\\
			\tfrac{N(z^{(0)}) \wedge u}{n} 
				\ar[r, mapsto] 
			\&
			\tr\left(\tfrac{z^{(01)}\wedge \res(u)}{x^{(0)}x^{(1)}}\right) 
			+ N(z^{(0)})
			\& 
			\begin{matrix}
				z^{(01)} \\
				\tfrac{z^{(0)} \wedge \res(u)}{x^{(0)}} 
			\end{matrix}
			\ar[r,mapsto, shift left=4]
			\ar[r,mapsto, shift right=4]
			\&
			\begin{matrix}
				0 \\
				\res(u)
			\end{matrix}
			\\
			\&
			\begin{matrix}
				\tfrac{z^{(01)} \wedge \res(u)}{x^{(0)}x^{(1)}}\\
				N(z^{(0)})
			\end{matrix}
			\ar[r,mapsto, shift left=4]
			\ar[r,mapsto, shift right=4]
			\&
			\begin{matrix}
				\frac{z^{(1)} \wedge \res(u)}{x^{(1)}} - \frac{z^{(0)} \wedge \res(u)}{x^{(0)}}\\
				0 \hphantom{xxxxxxxxxxxxxxxx}
			\end{matrix}
			\&
			\begin{matrix}
				u \\
				z^{(0)}
			\end{matrix}
			\ar[r,mapsto, shift left=4]
			\ar[r,mapsto, shift right=4]
			\&
			\begin{matrix}
				0 \\
				0
			\end{matrix}
		\end{tikzcd}
	\]
	From the differentials above, we see that the homology in degree zero is $\uA$, and the homology in degree $1$ is $\uA\{z^{(0)}\} \oplus \uA\{u\}/\res(u) \cong \uA\{z^{(0)}\} \oplus \uI\{u\}$. In degree 2, the term from the $n$-divided Koszul complex $\ndivkoszul_\bullet$ contributes nothing to homology, while the term from $\uK_\bullet$ emerges unscathed. So 
	\[ 
		\Tor_2^{\uR}(\uA,\uA) = \uA\{z^{(0)} \wedge z^{(1)}\}.
	\]
	In degree 4, we observe that the differential followed by projection onto $\uA\{N(z^{(0)})\}$ is an isomorphism, so it is injective. In particular, $\Tor_4^{\uR}(\uA, \uA) = 0$. 
	
	Thus it remains to calculate the homology in degree 3. The kernel of the third differential contains both $\tr\left(\tfrac{z^{(01)}\wedge \res(u)}{x^{(0)}x^{(1)}}\right)$ (which is sent to the transfer of a difference of Weyl conjugates) and $N(z^{(0)})$. In fact, these two elements generate the kernel, which is 
	\[ 
		\uA\{N(z^{(0)})\} \oplus \ubbZ\{\tr\left(\tfrac{z^{(01)}\wedge \res(u)}{x^{(0)}x^{(1)}}\right)\}.
	\]
	Note that the submodule of $\uA\{(\tfrac{z^{(01)}\wedge \res(u)}{x^{(0)}x^{(1)}})_{C_3/e}\}$ generated by the transfer is a copy of the constant module $\ubbZ$; this happens because the restriction of $\tr(g)$ in $\uA\{g_{C_3/e}\}$ is a sum of Weyl conjugates of $g$, and transfer of this sum is $3 \tr(g)$. 

	The homology is then the quotient of the above by
	\[ 
		\tr\left(\tfrac{z^{(01)}\wedge \res(u)}{x^{(0)}x^{(1)}}\right) 
		+ N(z^{(0)}),
	\]
	which is isomorphic to the cokernel of
	\[
		\begin{tikzcd}
			\bbZ[t]/t^2-3t
				\ar[d, bend right=50, "t \mapsto 3"'] 
				\ar[r, "1 \mapsto {(1,1)}"]
			& 
			\bbZ[t]/t^2-3t
				\ar[d, bend right=50, "t \mapsto 3"'] 
				\ar[r, phantom, "\bigoplus" description]
			&
			\bbZ
				\ar[d, bend right=50, "1"'] 
			\\
			\bbZ
				\ar[u, bend right=50, "\cdot t"']
				\ar[r, "1 \mapsto {(1,1)}"']
			&
			\bbZ
				\ar[u, bend right=50, "\cdot t"']
				\ar[r, phantom, "\bigoplus" description]			
			&
			\bbZ
				\ar[u, bend right=50, "3"']
			\\
			\uA \ar[r, "1 \mapsto {(1,1)}"]
			& 
			\uA \oplus \ubbZ
		\end{tikzcd}
	\]
	On the underlying level, this identifies the two copies of $\bbZ$. Then on the top level, this identifies $\ubbZ(C_3/C_3)$ with the trivial $G$-sets in $\uA(C_3)$. Finally, because the transfer of $1 \in \uA(C_3/e)$ is $t$ but the transfer of $1 \in \ubbZ(C_3/e)$ is $3$, we see that $t = 3$ in the top level of cokernel, and we are left with a copy of $\ubbZ$. So we have:
	\[ 
		\Tor^{\uR}_3(\uA,\uA) = \ubbZ\{(N(z^{(0)}))_{C_3/C_3}\}.
	\]
\end{example}

To calculate $\uTor$ in general, we now take $H_*(\uF_\bullet \boxtimes_{\uR} \uA)$. Since all differentials of $\uK_\bullet$ contain multiplications by $x$, the differentials in $\uK_\bullet\boxtimes_{\uR} \uA$ vanish. On the other hand, the differentials in $\ndivkoszul_\bullet$ only contain coefficients in $\uA$, hence we obtain $H_*(\ndivkoszul_\bullet \boxtimes_{\uR}\uA)\cong H_*(\ndivkoszul_\bullet)\boxtimes_{\uR} \uA$. Finally, the morphism $f_k\colon \ndivkoszul_k\to \uK_k$ contains factors of $x$ or $n$ except for $k=p$, where it is the identity. Hence, the long exact sequence of the mapping cone resulting from applying $- \boxtimes_{\uR}\uA$ shows that 
\[ \uTor_k^{\uR}(\uA, \uA) \cong \begin{cases}
\uA & \textup{ for } k=0,\\
\uA\{z_e\} \oplus \uI\{u_{C_p} \} & \textup{ for } k=1,\\
\bigoplus_{I_k} \uA\{(z^{(i_1)} \wedge \ldots\wedge z^{(i_k)})_e\} & \textup{ for } 2\leq k\leq p-1,\\
\ubbZ\{N(z^{(0)}) \} & \textup{ for } k=p, \textup{ and }\\
0 & \textup{ otherwise. }
\end{cases}\]
This finishes the proof of \cref{Prop:TambaraUnderlying} for $p\geq 3$.

\subsection{Koszul Complexes and Monomorphic Restriction Property}

The previous calculations simplify on a well-studied class of Tambara functors, namely those satisfying the monomorphic restriction property \cite[Definition 4.19]{Nak2012}:

\begin{definition}\label{def:MRP}
	A Mackey functor $\uM$ satisfies the monomorphic restriction property if for any subgroup inclusion $K\leq H\leq G$, the restriction 
	\[ \res_{K}^{H} \colon \uM(G/H)\to \uM(G/K)\]
	is injective.
\end{definition}

Tambara functors satisfying this monomorphic restriction property can be characterized as subfunctors of fixed-point Tambara functors, i.e., $\uT(G/H) = R^H$ for some commutative ring $R$ with a $G$-action, by \cite[Proposition 4.21]{Nak2012}. This is also a necessary condition for a Mackey functor to be a zero-slice of an equivariant spectrum by \cite[Proposition 4.50]{HHR2016}. In the context of equivariant algebra, an interesting class of Tambara functors satisfying the monomorphic restriction property are the field-like Tambara functors \cite[Definition 4.28, Theorem 4.32]{Nak2012}. These are Tambara functors which have no non-trivial Tambara ideals, mirroring the classical notion of a field.

We now consider how the previously studied Koszul resolutions behave over Tambara functors satisfying the monomorphic restriction property.

\begin{lemma}\label{lemma:MRPforFreeModules}
	Let $\uR$ be a $G$-Tambara functor satisfying the monomorphic restriction property.
	\begin{enumerate}
		\item For any morphism $f\colon U\to V$ of $G$-sets, the restriction $f^\ast\colon \uR(V)\to \uR(U)$ is injective.
		\item If $\uM$ is a free $\uR$-module, then $\uM$ satisfies the monomorphic restriction property.
		\item If $\uS$ is a free Tambara functor over $R$ which is free as an $\uR$-module, then $\uS$ satisfies the monomorphic restriction property.
	\end{enumerate}
\end{lemma}
\begin{proof}
	\begin{enumerate}
		\item We decompose the $G$-sets $U$ and $V$ into $G$-orbits. Then, the map $f$ decomposes into a disjoint union of compositions of fold maps with maps between transitive $G$-sets. Restriction along fold maps gives a diagonal map, which is injective, and restrictions along maps between transitive $G$-sets are injective by the monomorphic restriction property. Finally, disjoint unions are taken to direct sums, so the restriction $f^\ast$ is indeed injective.
		\item We may assume that $\uM= \uR\{x_{G/H}\}$ is a free module on a single generator, since the monomorphic restriction property is preserved under direct sum. We observe that the restriction $\uM(G/K)\to \uM(G/L)$ then agrees with the restriction $\uR(G/K\times G/H)\to \uR(G/L\times G/H)$, which is injective by the previous part, since $\uR$ satisfies the monomorphic restriction property.
		\item This follows from the previous assertion, since the monomorphic restriction property only depends on the Mackey functor structure. \qedhere
	\end{enumerate}
\end{proof}

\begin{proposition}\label{Prop:KoszulResolutionForMRP}
	Let $\uT$ be a Tambara functor satisfying the monomorphic restriction property, and $p>2$. 
	Then the Koszul complex $\uK_\bullet$, defined in the proof of \Cref{Prop:TambaraUnderlying}, defines a free resolution $\uK_\bullet\boxtimes \uT$ of $\uT^{\cO^\top}[x_{C_p/e}]/\langle x\rangle$, where $\langle x \rangle$ is the Green ideal generated by $x$. 
\end{proposition}

\begin{proof}
	We study the base-changed Koszul complex $\uK_\bullet\boxtimes \uT$. 
	Its zeroth homology is the desired $\uT^{\cO^\top}[x_{C_p/e}]/\langle x\rangle$ by right exactness of $\boxtimes$. 
	For all $1\leq i\leq p-1$, we observe that the homology $H_i(\uK_\bullet\boxtimes \uT)=0$ by the same arguments as in the proof of \Cref{Prop:TambaraUnderlying}, since we again analyse a Koszul complex, now over $\uT[x^{(0)}, \ldots, x^{(p-1)}]$. 
	Finally, the $p$-th homology is given as the kernel of the restriction. 
	But since $\uT$ satisfies the monomorphic restriction property and this transfers to the free algebra by \Cref{lemma:MRPforFreeModules}, this homology is also trivial. 
	Hence $\uK_\bullet \boxtimes \uT$ indeed is a free resolution.
\end{proof}

\section{A Koszul resolution for $C_9$}
\label{section:c9}

Computing $\uTor$ over $\uA^{\cO^\top}[x_{C_9/e}]$ is more involved than the $C_3$-case due to the presence of the additional norm $\nm_e^{C_9}(x)$ at the fixed level. 
A similar strategy to the $C_3$-case produces a bicomplex whose total complex is not quite a resolution -- the homology has terms generated by $\nm_e^{C_9}(x)$. 
Hence, we form a tricomplex with an additional layer designed to kill these norms. 
We will see that the total complex of this tricomplex is then a resolution. 
In general, the Koszul resolution for the free Tambara functor on $C_{p^n}$ will be the total complex of an $(n+1)$-dimensional complex (see \cref{section:KoszulCpn}).  

Throughout this section, let $\uR = \uA^{\cO^\top}[x_{C_9/e}]$ be the free Tambara functor on an underlying generator. 

\begin{example}
	Recall the free Tambara functor $\uR = \uA^{\cO^\top}[x_{C_9/e}]$ from \cref{thm:FreeIncompleteTambaraDescription}:
	\[
		\begin{tikzcd}[row sep=large]
			\bbZ[t'_{\vec v} \mid \vec v \in \bbZ^9_{\geq 0}][s_{\vec w} \mid \vec w \in \bbZ^3_{\geq 0}][n']/\sim
				\ar[d, bend right=50, "\res"{left}] \\
			\bbZ[t_{\vec v} \mid v \in \bbZ^9_{\geq 0}][n^{(0)}, n^{(1)}, n^{(2)}]/\sim
				\ar[u,"\nm" description] 
				\ar[u, bend right=50, "\tr"{right}]
				\ar[d, bend right=50, "\res"{left}] \\
			\bbZ[x^{(0)}, \ldots, x^{(8)}] 
				\ar[u,"\nm" description] 
				\ar[u, bend right=50, "\tr"{right}]
		\end{tikzcd}
	\]
	Elements $t_{\vec v} = \tr_e^{C_3}(x^{\vec v})$ and $t'_{\vec v} = \tr_e^{C_9}(x^{\vec v})$ are transfers of monomials in the underlying generators. 
	In particular, $t_{\vec 0} = \tr_e^{C_3}(1) = [C_3/e]$ and $t'_{\vec 0} = \tr_e^{C_9}(1) = [C_9/e]$. 
	The classes $n^{(i)} = \nm_e^{C_3}(x^{(i)})$ for $i = 0,1,2$ are the three distinct norms of the underlying generators, and $n' = \nm_e^{C_9}(x^{(0)})$ is the unique norm of any single underlying generator. 
	The classes $s_{\vec w} = \tr_{C_3}^{C_9}(n^{\vec w})$ are transfers of norms. 
	In particular, $s_{\vec 0} = \tr_{C_3}^{C_9}(\nm_e^{C_3}(1)) = [C_9/C_3]$. 
	The relations are as described in \cref{thm:FreeIncompleteTambaraDescription}. 
\end{example}

\subsection{Lifting the classical Koszul complex}

To begin, we first define a chain complex of $\uR$-modules that lifts the classical Koszul resolution of $\uR(C_9/e) = \bbZ[x^{(0)}, \ldots, x^{(8)}]$, as in \cref{Constr:KoddpFirst}. 

\begin{notation}\label{notation:wedge_power_set}
	Consider the $C_9$-set $Z = \{z^{(0)}, z^{(1)}, \ldots, z^{(8)}\}$, which is a single orbit of the form $C_9/e$. Write
	\[
		\bigwedge\nolimits^{\!\!k} Z := Z^k \setminus \{(z_1, \ldots, z_k) \in Z^k \mid z_i = z_j \text{ for some $i, j$} \}.
	\]
	We write elements of $\bigwedge^k Z$ as $z^{(i_1)} \wedge \ldots \wedge z^{(i_k)}$. 
	This inherits a $C_9$-action from the diagonal action on $Z^k$, which commutes with the $\Sigma_k$-action that permutes the tuples. 
\end{notation}
	
	Classically, the $k$-th stage of the Koszul resolution of $\bbZ[x^{(i)}]$ arises from a free module on $\bigwedge^k Z$ by identifying the permutation $\Sigma_k$-action with the action by the sign of the permutation. 
	Equivariantly, however, we want to only add one $\uR$-module generator for each orbit of the $C_9$-action to avoid double-counting. To do so, we divide $(\bigwedge^k Z)/\Sigma_k$ into orbits:\footnote{Here, since $\bigwedge^k Z$ is only a $\Sigma_k$-set (or a $C_9 \times \Sigma_k$-set) and not a $\Sigma_k$-module, by $(\bigwedge^k Z)/\Sigma_k$ we mean the coinvariants of the $\Sigma_k$-action on $\bigwedge^k Z$. Of course, this misses out on the signs that are important in the Koszul complex, but this still gives the correct count of generators. After all, $x$ and $-x$ generate the same $\uR$-module. } 
	\begin{equation}
		\label{C9gens}
		\left(\bigwedge\nolimits^{\!\!k} Z\right)/\Sigma_k \cong 
		\begin{cases}
			C_9/e & (k = 1) \\
			4(C_9/e) & (k = 2) \\
			9(C_9/e) + C_9/C_3 & (k = 3)\\
			14(C_9/e) & (k = 4) \\
			14(C_9/e) & (k = 5) \\
			9(C_9/e) + C_9/C_3 & (k = 6)\\
			4(C_9/e) & (k = 7) \\
			C_9/e & (k = 8) \\
			C_9/C_9 & (k = 9) \\
			\emptyset & (k > 9)
		\end{cases}
	\end{equation}
To build our lift of the classical Koszul complex, we add generators at the levels prescribed by the orbits above in each degree. 
	
\begin{construction}
	\label{c9KoszulLift}
	We construct a chain complex $\uK_\bullet$ of $\uR$-modules.
	Let $\uK_0 = \uR$, and for $0 < k \leq 9$, let 
	\[
		\uK_r = \uR\{(\textstyle\bigwedge^{\!r} Z)/\Sigma_r\}.
	\]
	Note $\uK_r = 0$ for $r >9$ because $\textstyle\bigwedge^{\!r} Z/\Sigma_r = \emptyset$. 
	It will be convenient to have names for the generators of $\uK_r$. 
	To do so, we pick a representative for each $C_9$-orbit in $(\textstyle\bigwedge^{\!r} Z)/\Sigma_r$ with the convention that the indices are in increasing order, beginning with $0$. 
	For example, we pick 
	\[
		z^{(0)} \wedge z^{(1)}, \quad 
		z^{(0)} \wedge z^{(2)}, \quad 
		z^{(0)} \wedge z^{(3)}, \quad
		z^{(0)} \wedge z^{(4)}
	\]
	for representatives of the four $C_9$-orbits of $(\textstyle\bigwedge^{\!2} Z)/\Sigma_2$. 
	If a representative $z^{(i_1)} \wedge \ldots \wedge z^{(i_r)}$ generates an orbit of the form $C_9/e$, we take this as the name of the $\uR$-module generator of $\uK_r$. 
	
	There are two cases where $z^{(i_1)} \wedge \ldots \wedge z^{(i_r)}$ generates an orbit of the form $C_9/C_3$: 
	\[
		z^{(0)} \wedge z^{(3)} \wedge z^{(6)} \qquad 
		\text{ and } \qquad 
		z^{(0)} \wedge z^{(3)} \wedge z^{(6)} \wedge z^{(1)} \wedge z^{(4)} \wedge z^{(7)}
	\] 
	in degrees 3 and 6, respectively. 
	These arise from the action of $C_3 \subseteq C_9$ on $z^{(0)}$ and $z^{(0)} \wedge z^{(1)}$. Therefore, we write  
	\[ 
		N_e^{C_3}(z^{(0)}) 
		\qquad \text{ and } \qquad
		N_e^{C_3}(z^{(0)} \wedge z^{(1)}), 
	\]
	for the corresponding generators in $\uK_3$ and $\uK_6$. The letter $N$ is supposed to indicate a formal norm of a variable. 
	Similarly, we write the generator of $\uK_9 = \uR\{C_9/C_9\}$ as
	\[
		N_e^{C_9}(z^{(0)})
	\]
	
	We take the convention that  
	\[
		z^{(i_{\sigma(1)})} \wedge \ldots \wedge z^{(i_{\sigma(r)})} = 
		\sgn(\sigma) z^{(i_1)} \wedge \ldots \wedge z^{(i_r)},
	\]
	for a permutation $\sigma \in \Sigma_r$, 
	and similarly 
	\[
		N_e^{C_{3^m}}( z^{(i_{\sigma(1)})} \wedge \ldots\wedge z^{(i_{\sigma(r)})}) 
		= \sgn(\sigma) N_e^{C_{3^m}}( z^{(i_1)}\wedge \ldots\wedge z^{(i_r)} ).
	\]

	The differentials are determined by 
	\[
		\partial_r (z^{(i_1)} \wedge\ldots \wedge z^{(i_r)}) 
			= 
		\sum_{j=1}^r (-1)^{j-1} x^{(i_j)}
		(z^{(i_1)} \wedge \ldots \wedge \widehat{z^{(i_j)}} \wedge \ldots \wedge z^{(i_r)})
	\]
	on the underlying generators, which is exactly the usual Koszul differential.
	For generators at other levels, 
	\begin{align*}
		\partial_2 ( N_e^{C_3} (z^{(0)}))
		& = 
		x^{(0)}\tr_e^{C_3}(z^{(3)} \wedge z^{(6)}).\\
		\partial_5 ( N_e^{C_3}(z^{(0)} \wedge z^{(1)}) )
		& = 
		x^{(0)} \tr_e^{C_3}(z^{(1)} \wedge z^{(3)} \wedge z^{(4)} \wedge z^{(6)} \wedge z^{(7)}) 
		- x^{(1)} \tr_e^{C_3}(z^{(0)} \wedge z^{(3)} \wedge z^{(4)} \wedge z^{(6)} \wedge z^{(7)})\\
		\partial_8 ( N_e^{C_9}(z^{(0)}) )
		& = 
		x^{(0)} \tr_e^{C_9} (z^{(1)} \wedge z^{(2)} \wedge \ldots \wedge z^{(8)}).
	\end{align*}
	Note that the above differentials are again compatible with the classical Koszul differential after restriction to the underlying level. 
	As in \Cref{remark:KoszulDiffOfNorm}, these differentials satisfy the Leibniz rule for the derivation of a norm.
\end{construction}

This construction gives the following chain complex of $\uR$-modules (shown without generator names), augmented to $\uA$ via the map $\varepsilon \colon \uK_0 \cong \uR \to \uA$ that sends $x^{(i)}$ to zero. 
\[
	\hspace*{-2cm}
	\begin{tikzcd}
			0
			\ar[r] 
				&
			\uR\{\sfrac{C_9}{C_9}\}
			\ar[r, "\partial_8"]
				&
			\uR\{\sfrac{C_9}{e}\}
			\ar[r, "\partial_7"]
				&
			\uR\{4(\sfrac{C_9}{e})\}
			\ar[r, "\partial_6"]
				&
			\uR\{9 (\sfrac{C_9}{e}) + \sfrac{C_9}{C_3}\}
			\ar[r, "\partial_5"]
				&
			\uR\{14 (\sfrac{C_9}{e})\}
			\ar[r, "\partial_4"]
				&
			\cdots
				\\
			\cdots 
			\ar[r, "\partial_4"]
				& 
			\uR\{14 (\sfrac{C_9}{e})\}
			\ar[r, "\partial_3"]
				&
			\uR\{9 (\sfrac{C_9}{e}) + \sfrac{C_9}{C_3}\}
			\ar[r, "\partial_2"] 
				& 
			\uR\{4(\sfrac{C_9}{e})\}
			\ar[r, "\partial_1"]
				& 
			\uR\{\sfrac{C_9}{e}\}
			\ar[r, "\partial_0"]
				& 
			\uR\{\sfrac{C_9}{C_9}\}
			\ar[r, "\varepsilon"]
				& 
			\uA
	\end{tikzcd}
\]

This is not yet a resolution of $\uA$ since it has nontrivial homology:

\begin{proposition}
	\label{C9HomologyOfKoszulLift}
	The complex $\uK_\bullet$ constructed above has homology:
	\[
		H_t(\uK_\bullet) \cong
		\begin{cases}
			\uR /\langle x \rangle^\bot \quad & \text{ if } t=0, \\
			\ker(\res^{C_3}_e) \quad & \text{ if } t=3, \\
			\ker(\res^{C_3}_e) \quad & \text{ if } t=6, \\
			\ker(\res^{C_9}_e) \quad & \text{ if } t=9, \\
			0 \quad & \text{ otherwise,}
		\end{cases}
	\]
	where 
	\begin{itemize}
		\item $\langle x \rangle^\bot$ denotes the Green ideal generated by the variable $x^{(0)}$
		\item $\ker(\res^H_e)$ denotes the sub-$\uR$-module of $\uK_t$ generated by the kernel of the restriction in $\uK_t$. 
	\end{itemize}
\end{proposition}

\begin{proof}
	The complex $\uK_\bullet$ agrees with the classical Koszul complex on the underlying level, so it is exact at that level. 
	
	In degree zero, the homology is $\uR$ modulo the image of the differential $\partial_0$, determined by  $z^{(i)} \mapsto x^{(i)}$. 
	Thus, the image of $\partial_0$ is exactly the sub-$\uR$-module of $\uR$ generated by sums, products, transfers, and restrictions of the $x^{(i)}$. 
	This is succinctly described as the Green ideal generated by these variables. 
	
	In positive degree $k > 0$, the kernel of $\partial_k \colon \uK_{k+1} \to \uK_k$ is exactly the image of $\partial_{k+1}$, just as in a classical Koszul complex, unless the domain of $\partial_k$ is generated at a different level than the codomain. This occurs in degrees 3, 6, and 9. 
	In degree 3, for example, the differential is given by 
	\[
		f N_e^{C_3}(z^{(0)} \wedge z^{(3)} \wedge z^{(6)}) 
		\overset{\partial_2}{\longmapsto}
	 	\res^{C_3}_e(f) x^{(0)} \tr_e^{C_3}(z^{(3)} \wedge z^{(6)})
	\]
	and so the kernel of $\partial_2$ is the sub-$\uR$-module generated by the kernel of $\res^{C_3}_e$ by \cref{lemma:kernels}. 
	These kernels are never in the image of $\partial_3$, however, since $\uK_4$ only has generators at the underlying level, so $\partial_3$ contains only transfers in $\uK_3(C_9/C_3)$. 
	A transfer is never in the kernel of restriction in these free $\uR$-modules. 
	Hence, $H_3(\uK_\bullet) \cong \ker(\res^{C_3}_e)$. 
	
	A similar calculation shows that the homology in degree 6 is also given by a kernel of restriction, and in degree 9 this is another application of \cref{lemma:kernels}. 
\end{proof}

\begin{remark}
	$\uR$ is generated as an $\uR$-module by sum, products, restrictions, transfers, and norms of the variables $x^{(i)}$ and $1 \in \uR(C_9/C_9)$. 
	The quotient by the Green ideal $\langle x \rangle^\bot$ leaves only the norms of the variables, and the copy of $\uA$ generated by $1 \in \uR(C_9/C_9)$.
	\[
		H_0(\uK_\bullet) 
			\cong \uA\{1\} 
			\oplus \bigoplus_{\text{monomials } f} 
				\uA\{\nm_e^{C_3}(f)\} 
			\oplus \bigoplus_{\text{monomials } g} 
				\uA\{\nm_e^{C_9}(g)\}
	\]    
\end{remark}

\subsection{Building a resolution}

As in the $C_p$-case, we eliminate these norms by adding dimensions to our resolution. 
This time, however, we must form a tricomplex: we add a dimension to kill norms at level $C_9/C_3$, and another to eliminate norms at level $C_9/C_9$. 

\begin{construction}
	\label{C9BigKoszulModules}
	Define a collection $\bigkoszul_{r,s,t}$ of free $\uR$-modules indexed by $0 \leq r \leq 9$, $0 \leq s \leq 3$, and $0 \leq t \leq 1$ as follows. 
	In \cref{C9BigKoszulDifferentials} below, we add differentials to turn this into a tricomplex of free $\uR$-modules. 
	First, we set
	\[
		\bigkoszul_{r,0,0} := \uK_r. 
	\] 
	In particular, $\bigkoszul_{0,0,0} = \uK_0 = \uR$. 
	Next, consider the $C_9$-set $U = \{u^{(0)}, u^{(1)}, u^{(2)}\} \cong C_9/C_3$, and let 
	\[
		\bigkoszul_{0,s,0} := \uR\{(\textstyle \bigwedge^s U)/\Sigma_s\} 
	\]
	for $0 < s \leq 3$. 
	Let 
	\[
		\bigkoszul_{0,0,1} := \uR\{V\} = \uR\{v_{C_9/C_9}\},
	\]
	where $V$ is the trivial $C_9$-set $V = \{v\} \cong C_9/C_9$.
	In general, define 
	\[
		\bigkoszul_{r,s,t} := \bigkoszul_{r,0,0} \boxtimes \bigkoszul_{0,s,0} \boxtimes \bigkoszul_{0,0,t} 
	\]
\end{construction}
	
	\medskip
	
	Recall that for finite $G$-sets $A$, $B$, and $C$, 
	\[
		\uR\{A\} \boxtimes \uR\{B\} \boxtimes \uR\{C\} \cong \uR\{A \times B \times C\}. 
	\]
	Therefore, $\bigkoszul_{r,s,t}$ is the free $\uR$-module generated by the product of generating $G$-sets of $\bigkoszul_{r,0,0}$, $\bigkoszul_{0,s,0}$, and $\bigkoszul_{0,0,t}$. Since the generating set of $\bigkoszul_{0,0,t}$ is always $C_9/C_9$, the generating set of $\bigkoszul_{r,s,t}$ in effect only depends on $\bigkoszul_{r,0,0}$ and $\bigkoszul_{0,s,0}$. 
	
\begin{example}
	In tridegrees $(3,1,0)$ and $(3,1,1)$, $\bigkoszul_{3,1,0}$ and $\bigkoszul_{3,1,1}$ are free $\uR$-modules generated by the $C_9$-set  
	\[	
		\big(9(\sfrac{C_9}{e}) + \sfrac{C_9}{C_3} \big) \times \sfrac{C_9}{C_3} \times \sfrac{C_9}{C_9} \cong 27(\sfrac{C_9}{e}) + 3(\sfrac{C_9}{C_3}). 
	\]
	The finite $G$-sets generating $\bigkoszul_{r,s,t}$ can be seen in \cref{C9tricomplex}. 
\end{example}
	
\begin{notation}
	\label{notation:generator names for big koszul C9}
		As before, it will be convenient to name the generators of $\bigkoszul_{r,s,t}$. 
		We have already named the generators of $\bigkoszul_{r,0,0} = \uK_0$ in \cref{c9KoszulLift}, and $\bigkoszul_{0,0,1} = \uR\{v_{C_9/C_9}\}$ has generator $v$.		 
		Next, we label the generators of $\bigkoszul_{0,s,0}$ as follows. 
		\begin{align*}
			\bigkoszul_{0,1,0} &= \uR\{C_9/C_3\} = \uR\{u^{(0)}_{C_9/C_3}\} \\
			\bigkoszul_{0,2,0} &= \uR\{C_9/C_3\} = \uR\{(u^{(0)} \wedge u^{(1)})_{C_9/C_3}\} \\
			\bigkoszul_{0,3,0} &= \uR\{C_9/C_9\} = \uR\{N_{C_3}^{C_9}(u^{(0)})_{C_9/C_9}\}
		\end{align*}
	
		With these labels in hand, we produce labels for the remaining generators of $\bigkoszul_{r,s,t}$. First, we introduce some convenient shorthand. 
		For $I = \{i_1, \ldots, i_k\} \subseteq \{0,\ldots,8\} \cong C_9/e$, we write 
		\[
			z^{(I)} = z^{(i_1)} \wedge z^{(i_2)} \wedge \ldots \wedge z^{(i_k)}, 
		\]
		and similarly write $u^{(J)}$ for $J \subseteq \{0,1,2\} \cong C_9/C_3$. 
		For the underlying Tambara generators $x^{(i)}$ of $\uR$, we write 
		\[
			x^{(I)} = x^{(i_1)}x^{(i_2)} \cdots x^{(i_k)} = \prod_{i \in I} x^{(i)}. 
		\]
		
		Given $I \subseteq C_9/e$ and $J \subseteq C_9/C_3$, let 
		\[
			\Inc(I,J) = \{i \in I \mid i + C_3 \in J\}.
		\]
		We call this the \emph{incidence set} of $I$ and $J$. 
		
		A pair $(I,J)$ with $|I| = r$, $I \subseteq C_9/e$ and $|J| = s$, $J \subseteq C_9/C_3$ determines a pair of generators of $z^{(I)}$ of $\bigkoszul_{r,0,0}$ and $u^{(J)}$ of $\bigkoszul_{0,s,0}$. 
		If this pair of generators represents an orbit of the form $C_9/e$, then we denote the corresponding generators of $\bigkoszul_{r,s,0}$ and $\bigkoszul_{r,s,1}$ by 
		\[
			\frac{z^{(I)} \wedge R^{C_3}_e(u^{(J)})}{x^{(\Inc(I,J))}} 
			\qquad \text{ and } \qquad 
			\frac{z^{(I)} \wedge R^{C_3}_e(u^{(J)}) \wedge R^{C_9}_{e}(v)}{x^{(\Inc(I,J))} \res_e^{C_3}\nm_e^{C_3}(x^{(J)})},
		\]	
		respectively. 
		Note that $\res^{C_3}_e(\nm_e^{C_3}(x^{(j)})) = x^{(j)} x^{(j+3)} x^{(j+6)}$, so the denominators are both monomials in the variables $x^{(i)}$.

		If the pair $(I,J)$ represents an orbit of the form $C_9/C_3$, then we know that $I$ has stabilizer $C_3 \subseteq C_9$. 
		Write $\overline{I} \subseteq C_9/e$ for a set of representatives of the $C_9$-orbits of $I$. 
		Note that $\Inc(\overline{I},J)$ is a set of representatives for the $C_9$-orbits of $\Inc(I,J)$. 
		We denote the corresponding generators of $\bigkoszul_{r,s,0}$ and $\bigkoszul_{r,s,1}$ by 
		\[
			\frac{N_e^{C_3}(z^{(\overline I)}) \wedge u^{(J)}}{\nm_e^{C_3}(x^{(\Inc(\overline I,J))})} 	
			\qquad
			\text{ and } 
			\qquad
			\frac{N_e^{C_3}(z^{(\overline I)}) \wedge u^{(J)} \wedge R^{C_9}_{C_3}(v)}{\nm_e^{C_3}(x^{(\Inc(\overline I,J))}) \nm_e^{C_3}(x^{(J)})},
		\]
		respectively. 
		
		Generators of the form $C_9/C_9$ occur in the corners: tridegrees $(0,0,0)$, $(9,0,0)$, $(0,3,0)$, $(9,3,0)$, $(0,0,1)$ $(9,0,1)$, $(0,3,1)$, and $(9,3,1)$. 
		These generators are: 
		\[
		\begin{array}{ll}
			1 \in \bigkoszul_{0,0,0} & 
			v \in \bigkoszul_{0,0,1}\\[1em]
			N_e^{C_9}(z^{(0)}) \in \bigkoszul_{9,0,0} & 
			\displaystyle \frac{N_e^{C_9}(z^{(0)}) \wedge v}{\nm_e^{C_9}(x^{(0)})} \in \bigkoszul_{9,0,1}\\[1em]
			N_{C_3}^{C_9}(u^{(0)}) \in \bigkoszul_{0,3,0} & 
			\displaystyle \frac{N_{C_3}^{C_9}(u^{(0)}) \wedge v}{\nm_e^{C_9}(x^{(0)})} \in \bigkoszul_{0,3,1}\\[1em]
			\displaystyle \frac{N_e^{C_9}(z^{(0)}) \wedge N_{C_3}^{C_9}(u^{(0)})}{\nm_e^{C_9}(x^{(0)})} \in \bigkoszul_{9,3,0} & 
			\displaystyle \frac{N_e^{C_9}(z^{(0)}) \wedge N_{C_3}^{C_9}(u^{(0)}) \wedge v}{\nm_e^{C_9}(x^{(0)})\nm_e^{C_9}(x^{(0)})}\in \bigkoszul_{9,3,1} 
		\end{array}
		\]
\end{notation}

\begin{remark}
	Although the sane choice of generator names $\bigkoszul_{r,s,t}$ appears to be triples of generators, one from each of $\bigkoszul_{r,0,0}$, $\bigkoszul_{0,s,0}$, and $\bigkoszul_{0,0,t}$, this doesn't actually give the right number of generators because the generators of each of the box product factors may live in different orbits. 
	Moreover, the differentials in this tricomplex differ slightly from the obvious choices of $\partial_{1,0,0} = \partial \boxtimes 1 \boxtimes 1$, etc., and the generator labels given above help us remember exactly how these differentials behave. 
\end{remark}
	
\begin{remark}
	The denominators of the generators only ever contain monomials in the $x^{(i)}$ or their norms. 
	These are designed so that, after performing a Koszul differential on $z^{(I)}$ or $u^{(J)}$, these denominators cancel with some of the coefficients of the Koszul differential. 
\end{remark}
	
\begin{remark}
	We distinguish between formal norms of module generators $z$, $u$, and $v$ and norms of elements of the Tambara functor by writing (for example) $N_e^{C_3}(z^{(I)})$ for a formal norm and $\nm_e^{C_3}(x^{(I)})$ for an element of $\uR$. 
	Similarly, we write $R^{C_3}_e(u^{(J)})$ for restrictions of module generators, mostly to save space in the generator names.  
\end{remark}

\begin{construction}
	\label{C9BigKoszulDifferentials}
	We add differentials to the collection $\bigkoszul_{r,s,t}$ of free $\uR$-modules defined in \cref{C9BigKoszulModules} to make it into a tricomplex. 
	
	The subcomplex $\bigkoszul_{0,s,0}$ will again look like a Koszul complex, but this time at level $C_9/C_3$. We define differentials $\partial_{0,1,0} \colon \bigkoszul_{0,s,0} \to \bigkoszul_{0,s-1,0}$ by 
	\begin{align*}
		\partial_{0,1,0}(u^{(0)}) &= \nm_e^{C_3}(x^{(0)}), \\
		\partial_{0,1,0}(u^{(0)}\wedge u^{(1)}) &= \nm_e^{C_3}(x^{(0)}) u^{(1)} - \nm_e^{C_3}(x^{(1)}) u^{(0)},  \\
		\partial_{0,1,0}(N_{C_3}^{C_9}(u^{(0)})) &= \nm_e^{C_3}(x^{(0)}) \tr_{C_3}^{C_9}(u^{(1)} \wedge u^{(2)}).
	\end{align*}
	We also define a differential $\partial_{0,0,1} \colon \bigkoszul_{0,0,1} \to \bigkoszul_{0,0,0}$ by 
	\[
		\partial_{0,0,1}(v) = \nm_e^{C_9}(x^{(0)}).
	\]	
	A general differential will be determined by performing these Koszul-type differentials on $z^{(I)}$, $u^{(J)}$, or $v$, and then cancelling as many coefficients as possible with terms in the denominators of the generators. 
	Explicitly, we define differentials $\partial_{1,0,0}$, $\partial_{0,1,0}$ and $\partial_{0,0,1}$ of tridegrees $(-1,0,0)$, $(0,-1,0)$ and $(0,0,-1)$ as follows. 
	
	If the pair $(I,J)$ determines an orbit $C_9/e$, then  
	\begin{align*}
		\partial_{1,0,0}\left(
				\frac{z^{(I)} \wedge R^{C_3}_e(u^{(J)})}{x^{(\Inc(I,J))}}
			\right) 
		&= 
		\sum_{k=1}^{|I|}
			(-1)^{k-1} 
			x^{(i_k)} 
			\frac{z^{(I\setminus \{i_k\})}\wedge R^{C_3}_e(u^{(J)})}
				{x^{(\Inc(I, J))}},
		\\
		\partial_{0,1,0}\left(
				\frac{z^{(I)} \wedge R^{C_3}_e(u^{(J)})}{x^{(\Inc(I,J))}}
			\right) 
		&= 
		\sum_{k=1}^{|J|}
			(-1)^{k-1} 
			\res_e^{C_3}(\nm_e^{C_3}(x^{(j_k)})) 
			\frac{z^{(I)}\wedge R^{C_3}_e(u^{(J\setminus \{j_k\})})}
				{x^{(\Inc(I, J) )} }.
		\\
	\end{align*}
	The terms $x^{i_k}$ and $\res_e^{C_3} (\nm_e^{C_3} (x^{(j_k)}))$ pulled out by the Koszul differential are understood to cancel against terms in the denominator, if possible. Note also that $\res^{C_3}_e(\nm_e^{C_3}(x^{(j)})) = x^{(j)} x^{(j+3)} x^{(j+6)}$.
	
	\begin{align*}
		\partial_{1,0,0}\left(
			\frac{z^{(I)} \wedge R^{C_3}_e(u^{(J)}) \wedge R^{C_9}_{e}(v)}
				{x^{(\Inc(I,J))} \res_e^{C_3}\nm_e^{C_3}(x^{(J)})}
			\right)
		&= 
		\sum_{k=1}^{|I|}
			(-1)^{k-1} 
			x^{(i_k)}
			\frac{z^{(I \setminus \{i_k\})} 
					\wedge R^{C_3}_e(u^{(J)}) 
					\wedge R^{C_9}_{e}(v)}
				{x^{(\Inc(I,J))} \res_e^{C_3}\nm_e^{C_3}(x^{(J)})}
		\\
		\partial_{0,1,0}\left(
			\frac{z^{(I)} \wedge R^{C_3}_e(u^{(J)}) \wedge R^{C_9}_{e}(v)}
				{x^{(\Inc(I,J))} \res_e^{C_3}\nm_e^{C_3}(x^{(J)})}
			\right)
		&= 
		\sum_{k=1}^{|J|}
			(-1)^{k-1} 
			\res_e^{C_3} \nm_e^{C_3}(x^{(j_k)})
			\frac{z^{(I)} 
					\wedge R^{C_3}_e(u^{(J \setminus \{j_k\})}) 
					\wedge R^{C_9}_{e}(v)}
				{x^{(\Inc(I,J))} \res_e^{C_3}\nm_e^{C_3}(x^{(J)})}
		\\
		\partial_{0,0,1}\left(
			\frac{z^{(I)} \wedge R^{C_3}_e(u^{(J)}) \wedge R^{C_9}_{e}(v)}
				{x^{(\Inc(I,J))} \res_e^{C_3}\nm_e^{C_3}(x^{(J)})}
			\right)
		&= 
		\res_e^{C_9} \nm_e^{C_9}(x^{(0)}) 
		\frac{z^{(I)} \wedge R^{C_3}_e(u^{(J)})}{x^{(\Inc(I,J))} \res_e^{C_3}\nm_e^{C_3}(x^{(J)})}
	\end{align*}
	
	If the pair $(I,J)$ determines an orbit $C_9/C_3$, and $\overline I = \{\overline{i}_k\}$ is a set of representatives for the $C_9/C_3$-orbits of $I$, we define
	\begin{align*}
		\partial_{1,0,0} \left( 
				\frac{N_e^{C_3}(z^{(\overline{I})}) \wedge u^{(J)}}
					{\nm_e^{C_3}(x^{(\Inc(\overline I, J))})}
			\right)
		& = 
		\sum_{k=1}^{|\overline I|}
			(-1)^{k-1}
			\tr_e^{C_3} \left(
					x^{(\{\overline{i}_k\})}
					\frac{z^{(I \setminus \{\overline i_k\})} 
						\wedge R^{C_3}_e(u^{(J)})}
						{\res_e^{C_3}\nm_e^{C_3}(x^{(\Inc(\overline I, J))})}
				\right)
		\\
		\partial_{0,1,0} \left(
				\frac{N_e^{C_3}(z^{(\overline{I})}) \wedge u^{(J)}}
					{\nm_e^{C_3}(x^{(\Inc(\overline I, J))})}
			\right)
		& = 
		\sum_{k=1}^{|J|}
			(-1)^{k-1}
			\nm_e^{C_3}(x^{(\{j_k\})})
			\left(
					\frac{N(z^{(\overline I)}) \wedge u^{(J \setminus \{j_k\})}}
					{\nm_e^{C_3}(x^{(\Inc(\overline I,J)}))}
				\right)
	\end{align*}
	Again, the terms pulled out by the Koszul differentials on $z^{(I)}$ or $u^{(J)}$ are understood to cancel with terms in the denominators, if possible. 
	\begin{align*}
		\partial_{1,0,0} \left( 
				\frac{N_e^{C_3}(z^{(\overline I)}) 
					\wedge u^{(J)} 
					\wedge R^{C_9}_{C_3}(v)}
				{\nm_e^{C_3}(x^{(\Inc(\overline I,J))}) \nm_e^{C_3}(x^{(J)})}
			\right)
		&= 
		\sum_{k=1}^{|\overline I|}
			(-1)^{k-1}
			\tr_e^{C_3}\left(
				x^{(\overline i_k)}
				\frac{z^{(I \setminus \{\overline i_k\})} \wedge R_e^{C_3}(u^{(J)}) \wedge R_e^{C_9}(v)}{\res^{C_3}_e\nm_e^{C_3}(x^{(\Inc(\overline I,J))}) \res^{C_3}_e\nm_e^{C_3}(x^{(J)})}
			\right)
		\\
		\partial_{0,1,0} \left( 
				\frac{N_e^{C_3}(z^{(\overline I)}) 
					\wedge u^{(J)} 
					\wedge R^{C_9}_{C_3}(v)}
				{\nm_e^{C_3}(x^{(\Inc(\overline I,J))}) \nm_e^{C_3}(x^{(J)})}
			\right)
		&= 
		\sum_{k=1}^{|J|}
			(-1)^{k-1}
			\nm_e^{C_3}(x^{(j_k)})
			\frac{N_e^{C_3}(z^{(\overline I)}) 
					\wedge u^{(J \setminus \{j_k\})} 
					\wedge R^{C_9}_{C_3}(v)}
				{\nm_e^{C_3}(x^{(\Inc(\overline I,J))}) \nm_e^{C_3}(x^{(J)})}
		\\
		\partial_{0,0,1} \left( 
				\frac{N_e^{C_3}(z^{(\overline I)}) 
					\wedge u^{(J)} 
					\wedge R^{C_9}_{C_3}(v)}
				{\nm_e^{C_3}(x^{(\Inc(\overline I,J))}) \nm_e^{C_3}(x^{(J)})}
			\right)
		&= 
		\res_{C_3}^{C_9}\nm_e^{C_9}(x^{(0)})
		\frac{N_e^{C_3}(z^{(\overline I)}) \wedge u^{(J)}}
			{\nm_e^{C_3}(x^{(\Inc(\overline I,J))}) \nm_e^{C_3}(x^{(J)})}
	\end{align*} 
	In the corner with tridegree $(9,3,0)$, we have 
	\begin{align*}
		\partial_{1,0,0} \left(
				\frac{N_e^{C_9}(z^{(0)}) \wedge N_{C_3}^{C_9}(u^{(0)})}
					{\nm_e^{C_9}(x^{(0)})}
			\right)
		& =
			\tr_e^{C_9} \left(
					\frac{z^{(\{1,2,\ldots,8\})} \wedge R^{C_3}_e(u^{(\{0,1,2\})})}
					{x^{({1,2,\ldots,8})}}
				\right)
		\\
		\partial_{0,1,0} \left(
				\frac{N_e^{C_9}(z^{(0)}) \wedge N_{C_3}^{C_9}(u^{(0)})}
					{\nm_e^{C_9}(x^{(0)})}
			\right)
		& =
			\tr_{C_3}^{C_9}\left(
					\frac{N_e^{C_3}(z^{(\{0,1,2\})}) \wedge u^{(\{1,2\})}}
					{\nm_e^{C_3}(x^{(\{1,2\})})}
				\right)
	\end{align*}
	In the corner with tridegree $(9,0,1)$, we have 
	\begin{align*}
		\partial_{1,0,0} \left(
			\frac{N_e^{C_9}(z^{(0)}) \wedge v}{\nm_e^{C_9}(x^{(0)})}
		\right)
		&= 
			\tr_e^{C_9} \left( \frac{z^{(\{1,2,\ldots,8\})} \wedge R_e^{C_9}(v)}{x^{({1,2,\ldots,8})}}
			\right)
		\\
		\partial_{0,0,1} \left(
			\frac{N_e^{C_9}(z^{(0)}) \wedge v}{\nm_e^{C_9}(x^{(0)})}
		\right)
		&= N_e^{C_9}(z^{(0)})
	\end{align*}
	In the corner with tridegree $(0,3,1)$, we have 
	\begin{align*}
		\partial_{0,1,0} \left(
			\frac{N_{C_3}^{C_9}(u^{(0)}) \wedge v}{\nm_e^{C_9}(x^{(0)})}
		\right)
		&= 
			\tr_{C_3}^{C_9} \left(
				\frac{u^{(\{1,2\})} \wedge R_{C_3}^{C_9}(v)}{\nm_e^{C_3}(x^{(\{1,2\})})}
			\right)
		\\
		\partial_{0,0,1} \left(
			\frac{N_{C_3}^{C_9}(u^{(0)}) \wedge v}{\nm_e^{C_9}(x^{(0)})}
		\right)
		&= N_{C_3}^{C_9}(u^{(0)})
	\end{align*}
	Finally, in the corner with tridegree $(9,3,1)$, we have 
	\begin{align*}
		\partial_{1,0,0} \left(
			\frac{N_e^{C_9}(z^{(0)}) \wedge N_{C_3}^{C_9}(u^{(0)}) \wedge v}
				{\nm_e^{C_9}(x^{(0)})\nm_e^{C_9}(x^{(0)})}
		\right)
		&= 
			\tr_{e}^{C_9} \left(
				\frac{x^{(\{1,2,\ldots,8\})} \wedge R_e^{C_3}(u^{(\{0,1,2\})}) \wedge R_e^{C_9}(v)}{x^{(\{1,2,\ldots,8\})}\ \res_{e}^{C_9} \nm_e^{C_9}(x^{(0)})}
			\right)
		\\
		\partial_{0,1,0} \left(
			\frac{N_e^{C_9}(z^{(0)}) \wedge N_{C_3}^{C_9}(u^{(0)}) \wedge v}
				{\nm_e^{C_9}(x^{(0)})\nm_e^{C_9}(x^{(0)})}
		\right)
		&= 
			\tr_{C_3}^{C_9} \left(
				\frac{N_e^{C_3}(z^{(\{0,1,2\})}) \wedge u^{(\{1,2\})} \wedge R_{C_3}^{C_9}(v)}{\nm_e^{C_3}(x^{(\{1,2\})})\ \res_{C_3}^{C_9} \nm_e^{C_9}(x^{(0)})}
			\right)
		\\
		\partial_{0,0,1} \left(
			\frac{N_e^{C_9}(z^{(0)}) \wedge N_{C_3}^{C_9}(u^{(0)}) \wedge v}	
				{\nm_e^{C_9}(x^{(0)})\nm_e^{C_9}(x^{(0)})}
		\right)
		&= 
			\frac{N_e^{C_9}(z^{(0)}) \wedge N_{C_3}^{C_9}(u^{(0)})}	
				{\nm_e^{C_9}(x^{(0)})}		
	\end{align*}
	To understand how the terms pulled out by Koszul differentials are cancelled with terms in the denominators in the equations above, it's useful to remember that $\res_{C_3}^{C_9} \nm_e^{C_9}(x^{(0)}) = \nm_e^{C_3}(x^{(0)}) \nm_e^{C_3}(x^{(1)}) \nm_e^{C_3}(x^{(2)})$.
\end{construction}
	
\begin{example}
	\label{example:bigkoszul210}
		$\bigkoszul_{2,1,0}$ has twelve generators, all at level $C_9/e$: 
		\[
			\begin{array}{lll}
				\frac{z^{(0)} \wedge z^{(1)} \wedge R^{C_3}_e(u^{(0)})}{x^{(0)}} & 
				\frac{z^{(0)} \wedge z^{(1)} \wedge R^{C_3}_e(u^{(1)})}{x^{(1)}} & 
				\frac{z^{(0)} \wedge z^{(1)} \wedge R^{C_3}_e(u^{(2)})}{1} \\[1em]
				\frac{z^{(0)} \wedge z^{(2)} \wedge R^{C_3}_e(u^{(0)})}{x^{(0)}} & 
				\frac{z^{(0)} \wedge z^{(2)} \wedge R^{C_3}_e(u^{(1)})}{1} & 
				\frac{z^{(0)} \wedge z^{(2)} \wedge R^{C_3}_e(u^{(2)})}{x^{(2)}} \\[1em] 
				\frac{z^{(0)} \wedge z^{(3)} \wedge R^{C_3}_e(u^{(0)})}{x^{(0)}x^{(3)}} & 
				\frac{z^{(0)} \wedge z^{(3)} \wedge R^{C_3}_e(u^{(1)})}{1} & 
				\frac{z^{(0)} \wedge z^{(3)} \wedge R^{C_3}_e(u^{(2)})}{1} \\[1em] 
				\frac{z^{(0)} \wedge z^{(4)} \wedge R^{C_3}_e(u^{(0)})}{x^{(0)}} & 
				\frac{z^{(0)} \wedge z^{(4)} \wedge R^{C_3}_e(u^{(1)})}{x^{(4)}} & 
				\frac{z^{(0)} \wedge z^{(4)} \wedge R^{C_3}_e(u^{(2)})}{1} 	
			\end{array}
		\]
		A few sample differentials: 
		\begin{align*}
			\partial_{1,0,0}\left( \frac{z^{(0)} \wedge z^{(1)} \wedge R^{C_3}_e(u^{(0)})}{x^{(0)}} \right) 
			& = 
			\cancel{x^{(0)}}  \frac{z^{(1)} \wedge R^{C_3}_e(u^{(0)})}{\cancel{x^{(0)}}} - x^{(1)} \frac{z^{(0)} \wedge R^{C_3}_e(u^{(0)})}{x^{(0)}}\\
			& = 
			z^{(1)} \wedge R^{C_3}_e(u^{(0)}) - x^{(1)} \frac{z^{(0)} \wedge R^{C_3}_e(u^{(0)})}{x^{(0)}}
			\\[1em]
			\partial_{0,1,0}\left( \frac{z^{(0)} \wedge z^{(1)} \wedge R^{C_3}_e(u^{(0)})}{x^{(0)}} \right) & = 
			\res_e^{C_3}(\nm_e^{C_3}(x^{(0)}))\frac{z^{(0)} \wedge z^{(1)} }{x^{(0)}} 
			\\
			& = \frac{\cancel{x^{(0)}}x^{(3)}x^{(6)}}{\cancel{x^{(0)}}} z^{(0)} \wedge z^{(1)}\\
			& = x^{(3)}x^{(6)} z^{(0)} \wedge z^{(1)}
		\end{align*}
\end{example}
	
\begin{example}
		\label{example:bigkoszul310}
		$\bigkoszul_{3,1,0}$ has 27 generators at level $C_9/e$: 
		\[
		\begin{array}{lll}
			\frac{z^{(0)} \wedge z^{(1)} \wedge z^{(2)} \wedge R^{C_3}_e(u^{(0)})}
			{x^{(0)}}, & 
			\frac{z^{(0)} \wedge z^{(1)} \wedge z^{(2)} \wedge R^{C_3}_e(u^{(1)})}
			{x^{(1)}}, & 
			\frac{z^{(0)} \wedge z^{(1)} \wedge z^{(2)} \wedge R^{C_3}_e(u^{(2)})}
			{x^{(2)}} 
			\\[1em]
			\frac{z^{(0)} \wedge z^{(1)} \wedge z^{(3)} \wedge R^{C_3}_e(u^{(0)})}
			{x^{(0)} x^{(3)}}, & 
			\frac{z^{(0)} \wedge z^{(1)} \wedge z^{(3)} \wedge R^{C_3}_e(u^{(1)})}
			{x^{(1)}}, & 
			\frac{z^{(0)} \wedge z^{(1)} \wedge z^{(3)} \wedge R^{C_3}_e(u^{(2)})}
			{1} 
			\\[1em]
			\frac{z^{(0)} \wedge z^{(1)} \wedge z^{(4)} \wedge R^{C_3}_e(u^{(0)})}
			{x^{(0)}}, & 
			\frac{z^{(0)} \wedge z^{(1)} \wedge z^{(4)} \wedge R^{C_3}_e(u^{(1)})}
			{x^{(1)} x^{(4)}}, & 
			\frac{z^{(0)} \wedge z^{(1)} \wedge z^{(4)} \wedge R^{C_3}_e(u^{(2)})}
			{1} 
			\\[1em]
			\frac{z^{(0)} \wedge z^{(1)} \wedge z^{(5)} \wedge R^{C_3}_e(u^{(0)})}
			{x^{(0)}}, & 
			\frac{z^{(0)} \wedge z^{(1)} \wedge z^{(5)} \wedge R^{C_3}_e(u^{(1)})}
			{x^{(1)}}, &
			\frac{z^{(0)} \wedge z^{(1)} \wedge z^{(5)} \wedge R^{C_3}_e(u^{(2)})}
			{x^{(5)}} 
			\\[1em]
			\frac{z^{(0)} \wedge z^{(1)} \wedge z^{(6)} \wedge R^{C_3}_e(u^{(0)})}
			{x^{(0)} x^{(6)}}, & 
			\frac{z^{(0)} \wedge z^{(1)} \wedge z^{(6)} \wedge R^{C_3}_e(u^{(1)})}
			{x^{(1)}}, &
			\frac{z^{(0)} \wedge z^{(1)} \wedge z^{(6)} \wedge R^{C_3}_e(u^{(2)})}
			{1} 
			\\[1em]
			\frac{z^{(0)} \wedge z^{(1)} \wedge z^{(7)} \wedge R^{C_3}_e(u^{(0)})}
			{x^{(0)}}, &
			\frac{z^{(0)} \wedge z^{(1)} \wedge z^{(7)} \wedge R^{C_3}_e(u^{(1)})}
			{x^{(1)}}, &
			\frac{z^{(0)} \wedge z^{(1)} \wedge z^{(7)} \wedge R^{C_3}_e(u^{(2)})}
			{x^{(7)}} 
			\\[1em]
			\frac{z^{(0)} \wedge z^{(2)} \wedge z^{(4)} \wedge R^{C_3}_e(u^{(0)})}
			{x^{(0)}}, &
			\frac{z^{(0)} \wedge z^{(2)} \wedge z^{(4)} \wedge R^{C_3}_e(u^{(1)})}
			{x^{(4)}}, &
			\frac{z^{(0)} \wedge z^{(2)} \wedge z^{(4)} \wedge R^{C_3}_e(u^{(2)})}
			{x^{(2)}} 
			\\[1em]
			\frac{z^{(0)} \wedge z^{(2)} \wedge z^{(5)} \wedge R^{C_3}_e(u^{(0)})}
			{x^{(0)}}, &
			\frac{z^{(0)} \wedge z^{(2)} \wedge z^{(5)} \wedge R^{C_3}_e(u^{(1)})}
			{1}, &
			\frac{z^{(0)} \wedge z^{(2)} \wedge z^{(5)} \wedge R^{C_3}_e(u^{(2)})}
			{x^{(2)} x^{(5)}} 
			\\[1em]
			\frac{z^{(0)} \wedge z^{(2)} \wedge z^{(6)} \wedge R^{C_3}_e(u^{(0)})}
			{x^{(0)} x^{(6)}}, &
			\frac{z^{(0)} \wedge z^{(2)} \wedge z^{(6)} \wedge R^{C_3}_e(u^{(1)})}
			{1}, &
			\frac{z^{(0)} \wedge z^{(2)} \wedge z^{(6)} \wedge R^{C_3}_e(u^{(2)})}
			{x^{(2)}} 
		\end{array}
		\]
		and $3$ generators at level $C_9/C_3$
		\[	
			\frac{N_e^{C_3}(z^{(0)}) \wedge u^{(0)}}{\nm_e^{C_3}(x^{(0)})}, \qquad 
			\frac{N_e^{C_3}(z^{(0)}) \wedge u^{(1)}}{1}, \qquad
			\frac{N_e^{C_3}(z^{(0)}) \wedge u^{(2)}}{1}.
		\]
		A few differentials: 
		\begin{align*}
			\partial_{1,0,0}\left(\frac{N_e^{C_3}(z^{(0)}) \wedge u^{(0)}}{\nm_e^{C_3}(x^{(0)})} \right) & = \tr_e^{C_3}\left( \frac{z^{(1)} \wedge z^{(2)} \wedge R^{C_3}_e(u^{(0)})}{x^{(1)} x^{(2)}} \right)  \\
			\partial_{0,1,0} \left(\frac{N_e^{C_3}(z^{(0)}) \wedge u^{(0)}}{\nm_e^{C_3}(x^{(0)})}\right) 
			& = 
			\xcancel{\nm_e^{C_3}(x^{(0)})} \frac{ N_e^{C_3}(z^{(0)})}{\xcancel{\nm_e^{C_3}(x^{(0)})}} = N_e^{C_3}(z^{(0)})
		\end{align*}
		\begin{align*}
			\partial_{1,0,0}\left(\frac{z^{(0)} \wedge z^{(1)} \wedge z^{(3)} \wedge R^{C_3}_e(u^{(0)})}{x^{(0)}  x^{(3)}}\right)
			& = 
			\cancel{x^{(0)}} \frac{ z^{(1)} \wedge z^{(3)} \wedge R^{C_3}_e(u^{(0)})}
			{\cancel{x^{(0)}} x^{(3)}} \\
			& - x^{(1)} \frac{ z^{(0)} \wedge z^{(3)} \wedge R^{C_3}_e(u^{(0)})}{x^{(0)} x^{(3)}}\\
			&+\cancel{x^{(3)}}  \frac{ z^{(0)} \wedge z^{(1)} \wedge R^{C_3}_e(u^{(0)})}{x^{(0)}\cancel{ x^{(3)}}}\\[1em]
			\partial_{0,1,0}\left(\frac{z^{(0)} \wedge z^{(1)} \wedge z^{(3)} \wedge R^{C_3}_e(u^{(0)})}{x^{(0)}  x^{(3)}}\right) &= \res^{C_3}_e(\nm_e^{C_3}(x^{(0)}))
			\frac{z^{(0)} \wedge z^{(1)} \wedge z^{(3)}}{x^{(0)} x^{(3)}} \\
			& = \cancel{x^{(0)}} \cancel{x^{(3)}} x^{(6)} \frac{z^{(0)} \wedge z^{(1)} \wedge z^{(3)}}{\cancel{x^{(0)}} \cancel{x^{(3)}}} \\
			& = x^{(6)} z^{(0)} \wedge z^{(1)} \wedge z^{(3)}
		\end{align*}			
\end{example}
	
	A picture of this tricomplex (without generator names) appears in \cref{C9tricomplex}. 
	
	Let $\overline{\bigkoszul}$ denote the totalization of this tricomplex. 
	
\begin{theorem}
	The total complex $\overline{\bigkoszul}$ of this tricomplex $\bigkoszul_{r,s,t}$ is a resolution of $\uA$ by free $\uR$-modules. 
\end{theorem}
	
\begin{proof}
		To compute the homology of the total complex, we successively compute horizontal homology (in the $r$-direction), then vertical homology (in the $s$-direction), then homology in the third direction (in the $t$-direction). This strategy ensures that the homology vanishes on the $C_9/e$-level after taking horizontal homology (\cref{C9TricomplexHorizontalHomology}); on levels $C_9/e$ and $C_9/C_3$ after taking vertical homology of horizontal homology (\cref{C9TricomplexVerticalHomology}); and on all levels after taking the homology in the third direction. 
		
		We have already computed the horizontal homology of the first row in \cref{C9HomologyOfKoszulLift}. By design, the first row restricts to a Koszul complex on the underlying level, and therefore has zero homology on the underlying level save for $H_0(\uK_\bullet)(C_9/e) \cong \bbZ$. A similar calculation with judicious use of \cref{lemma:kernels} finds the homology of the subsequent rows on the underlying level; the horizontal differentials again behave as Koszul differentials and don't touch the $u$'s and $v$. 
		
		Note that the generators that remain after taking horizontal homology yield Koszul-type complexes in columns $0$, $3$, $6$, and $9$; these look like inflations of $C_3$-Koszul complexes. So the vertical homology of the horizontal homology vanishes unless the generators live at different levels, in which case we get a kernel of restriction, again using \cref{lemma:kernels}. In tridegree $(0,0,0)$, the variable $u \in \bigkoszul_{0,1,0}$ is sent to $\nm_e^{C_3}(x^{(0)})$, so we get the quotient of $\uR/\langle x \rangle$ by the Green ideal generated by $\nm_e^{C_3}(x^{(0)})$. The only tridegrees that survive after taking vertical homology are the corners (\cref{C9TricomplexVerticalHomology}).

		Finally, the differentials in the third direction are isomorphisms except in tridegree $(0,0,0)$, by inspection of \cref{C9BigKoszulDifferentials}. In tridegree $(0,0,0)$, the differential sends $v \in \bigkoszul_{0,0,1}$ to $\nm_e^{C_9}(x^{(0)})$. This leaves $\uR/\langle x^{(0)}, \nm_e^{C_3}(x^{(0)}), \nm_e^{C_9}(x^{(0)}) \rangle$. The quotient by this Green ideal kills all $\uR$-module (Green functor) generators of $\uR$ except $1 \in \uR(C_9/C_9)$, so 
		\[
			\uR/\langle x^{(0)}, \nm_e^{C_3}(x^{(0)}), \nm_e^{C_9}(x^{(0)}) \rangle 
			\cong 
			\uA
		\]	
		Hence, the homology vanishes except at the origin, where it is $\uA$. 
\end{proof}

\subsection{Computing Tor}

We can use this resolution to compute the Mackey-functor valued $\uTor$ of $\uA$ as an $\uR$-module. 

\begin{theorem}
	Let $\uR = \uA^{\cO^\top}[x_{C_9/e}]$ be the free $C_9$-Tambara functor on an underlying generator. 
	Consider $\uA$ as an $\uR$-module via the augmentation $\varepsilon \colon \uR \to \uA$, $x\mapsto 0$. Then
\begin{equation}
	\label{C9CaseFinalCalculation}
	\uTor_n^{\uR} (\uA, \uA) \cong 
	\begin{cases}
		\uA & (n=0)\\
		\uA\{\sfrac{C_9}{e}\} \oplus \uA\{\sfrac{C_9}{C_3}\}/\res_{e} \oplus \uA\{\sfrac{C_9}{C_9}\}/\res_{C_3} & (n=1)\\
\uA\{4(\sfrac{C_9}{e})\} \oplus \uA\{\sfrac{C_9}{C_3}\}/\res_{e} & (n=2)\\
		\uA\{9(\sfrac{C_9}{e})\} \oplus \ubbZ\{\sfrac{C_9}{C_3}\} \oplus \Inf_{C_3}^{C_9}(\ubbZ) & (n=3)\\
		\uA\{14(\sfrac{C_9}{e})\} \oplus \uL^{\oplus 2} & (n=4)\\
		\uA\{14(\sfrac{C_9}{e})\} \oplus \uL & (n=5)\\
		\uA\{9(\sfrac{C_9}{e})\} \oplus \ubbZ\{\sfrac{C_9}{C_3}\} & (n=6)\\
		\uA\{4(\sfrac{C_9}{e})\} \oplus \uL & (n=7)\\
		\uA\{\sfrac{C_9}{e}\} & (n=8)\\
		\ubbZ & (n=9).
	\end{cases}
\end{equation}
where
\begin{itemize}
	\item $\uT\{C_9/H\}$ is the free $\uT$-module generated at level $C_9/H$, for $\uT = \uA$ or $\uT = \ubbZ$,
	\item $\uA\{C_9/H\}/\res^H_K$ is the quotient of $\uA\{C_9/H\}$ generated by the submodule generated by restriction of generators to level $K$,
	\item $\uL$ is the Mackey functor
	\[
	\uL = 
	\begin{tikzcd}
		\bbZ/3\{\tr_{C_3}^{C_9}(g)\} 
			\ar[d, bend right=30]
			\\
		\bbZ/3\{g^{(0)}, g^{(1)}, g^{(2)}\} 
			\ar[u, bend right=30] 
			\ar[d, bend right=30]
			\\
		0
			\ar[u, bend right=30]
	\end{tikzcd}
	\]
\end{itemize}
\end{theorem}

\begin{proof}
	In order to calculate $\uTor^{\uR}(\uA,\uA)$, we take $\bigkoszul_{\bullet} \boxtimes_{\uR} \uA$ and calculate the homology. 
	For this, we use the standard spectral sequence and first calculate the horizontal homology, then the vertical homology of the horizontal homology, and finally compute the homology of the differential between the two layers.

	In the base-changed complex $\bigkoszul_{\bullet} \boxtimes_{\uR} \uA$, observe that all differentials containing a coefficient of $x^{(0)}$ or any of its Weyl conjugates or norms vanish. 
	The calculation of the horizontal homology is carried out in \cref{TorHorizontalHomology}. 
	We highlight one subtlety by way of example. 
	In general, elements which are not in the kernel of the horizontal differential may lie in the kernel after base-change. 
	For example, the generator $z \wedge u^{(1)}$ has nontrivial image $x^{(0)}u^{(1)}$ in the Koszul complex, but after base-change to $\uA$, it maps to zero and thus lies in the kernel of the horizontal differential. 
	However, it does not contribute to homology since it is in the image of $\frac{z^{(0)} \wedge z^{(1)} \wedge u^{(1)}}{x^{(1)}}$. 

	To calculate the vertical homology in columns 3 and 6, we have to calculate the homology of the two-term complex 
	\[
		\ker(\res^{C_3}_e) \to \uA\{\sfrac{C_9}{C_3}\}/\res_e^{C_3},
	\]
	where $\ker(\res^{C_3}_e)$ is the submodule of $\uA\{\sfrac{C_9}{C_3}\}$ generated by the kernel of restriction from $C_3$ to $e$. 
	The homomorphism is induced by the identity on $\uA$. 
	As Mackey functors, this homomorphism is 
	\[
	\begin{tikzcd}[row sep=huge]
		\langle t - 3 \rangle\{\tr_{C_3}^{C_9}(g)\} 
			\ar[r] 
			\ar[d, bend right=30] 
			& 
		A(C_3)/\langle t \rangle \{\tr_{C_3}^{C_9}(g)\} 
			\ar[d, bend right=30] 
			\\
		\langle t - 3 \rangle\{g^{(0)}, g^{(1)}, g^{(2)}\} 
			\ar[r] 
			\ar[u, bend right=30] 
			\ar[d, bend right=30] 
			& 
		A(C_3)/\langle t \rangle \{g^{(0)}, g^{(1)}, g^{(2)}\}
			\ar[u, bend right=30] 
			\ar[d, bend right=30] 
			\\
		0 
			\ar[r] 
			\ar[u, bend right=30] 
			& 
		0 
			\ar[u, bend right=30] 
	\end{tikzcd}
	\]
	where $g = g_{C_9/C_3}$ is a generator for $\uA\{C_9/C_3\}$, $A(C_3) \cong \bbZ[t]/\langle t^2 - 3t \rangle$ is the Burnside ring of $C_3$ with $t$ the class of a free orbit, and $\langle f\rangle$ is the ideal of $A(C_3)$ generated by $f$ for $f\in A(C_3)$. 
	We can see that this homomorphism is injective with cokernel 
	\[
	\uL = 
	\begin{tikzcd}
		\bbZ/3\{\tr_{C_3}^{C_9}(g)\} 
			\ar[d, bend right=30]
			\\
		\bbZ/3\{g^{(0)}, g^{(1)}, g^{(2)}\} 
			\ar[u, bend right=30] 
			\ar[d, bend right=30]
			\\
		0
			\ar[u, bend right=30]
	\end{tikzcd}
	\]
	In the top layer $\bigkoszul_{\bullet,\bullet,1}$ of the complex calculating $\uTor$, we claim that no differentials vanish.	
	Recall from \cref{notation:generator names for big koszul C9} that generators at level $C_9/e$ in tridegrees $(r,s,1)$ have denominators that are monomial in the $x^{(i)}$, and those $x^{(i)}$ that appear are exactly those that any $z^{(i)}$ or $u^{(j)}$ in the numerator are sent to under the horizontal and vertical Koszul differentials. 
	The upshot is that there are no coefficients of $x^{(i)}$ in the image of any vertical or horizontal differentials in the top layer. 
	In the bottom layer, $\bigkoszul_{\bullet, \bullet, 0}$, this is not the case, as the sample differentials in \cref{example:bigkoszul210,example:bigkoszul310} show.   

	Thus, most of the horizontal and vertical homology at this layer vanishes, as in the calculation of the homology of $\bigkoszul$. 
	The calculation is depicted in \cref{TorHorizontalHomology,TorVerticalHomology}.

	The differential in the third direction then only has an effect in the corners, that is, in those degrees where there are generators at level $C_9$ present. 
	\begin{itemize}
		\item At the $(0,0)$-spot, the differential is trivial, and thus we obtain a new homology group in total degree 1 of the form $\uA\{v_{C_9}\}/ \res^{C_9}_{C_3}$.
		\item At the $(0,3)$-spot, the surviving term in the top layer is $\ker(\res^{C_9}_{C_3})$. 
		The differential to the bottom layer is injective and yields as cokernel
		\[ 
			\Inf_{C_3}^{C_9}(\ubbZ)
				\cong
				\begin{tikzcd}
					\mathbb Z 
						\arrow[d, bend right=30, "1", swap] 
						\\
					\mathbb Z 
						\arrow[u, bend right=30, "3", swap] 
						\\
					0.
				\end{tikzcd} 
		\]
		\item At the $(9,0)$-spot, we obtain $\ker(\res^{C_9}_e)/\res^{C_9}_{C_3}$. 
		The differential to the bottom layer is injective, and the cokernel is isomorphic to $\ubbZ$, on the generator $N_e^{C_9} (z^{(0)})$.
		\item At the $(9,3)$-spot, the term in the top layer is $\ker(\res^{C_9}_e)$, and the differential to the bottom layer is the identity. 
		Thus both terms vanish in the homology of the total complex.
	\end{itemize} 
We observe that there is no homology of total degree more than $9$.

	We obtain the final result by summing over terms with the same total degree. 
	For this to yield the correct homology, we must argue that there are no extensions. 
	However, differentials of elements on the second page all vanish, either for degree reasons or because differentials of elements in $\ker(\res^{H}_K)$ are a sum of terms which each contain a restriction into $K$ or a subgroup of $K$. 
	Therefore, these elements also represent homology classes for the total complex. 
	Having a concrete cycle in the total complex for each element in the second page provides lifts that trivialize all extensions. 
	Thus, we obtain the $\uTor$-calculation by adding the terms with same total degree.
\end{proof}

\section{Koszul resolutions for cyclic groups of odd prime power order}
\label{section:KoszulCpn}

	Let $G = C_{p^n}$ with $p$ an odd prime and let $\uR := \uA^{\cO^{\top}}[x_{C_{p^n}/e}]$ be the free Tambara functor on an underlying generator. 
	In this section, we define a free $\uR$-module resolution of the Burnside Mackey functor $\uA$, generalizing the Koszul resolution of $\mathbb{Z}$ by free $\mathbb{Z}[x]$-modules. 
	Our resolution will be obtained by totalizing an $(n+1)$-dimensional complex $\bigkoszul_\bullet$ of free $\uR$-modules; to define this, we need one auxiliary definition:

\begin{definition}
	Let $G=C_{p^n}$ and $H=C_{p^m}$. 
	Let $T= C_{p^n}/C_{p^m}=\{ z^{(0)}, \ldots, z^{(p^{n-m}-1)}\}$ be the $C_{p^n}$-set of Weyl conjugates of a generator $z$ at level $C_{p^n}/C_{p^m}$ and let $S=C_{p^n}/C_{p^{m+1}}=\{ u^{(0)}, \ldots, u^{(p^{n-m-1}-1)}\}$ be the  $C_{p^{n}}$-set of Weyl conjugates of a generator $u$ at level $C_{p^n}/C_{p^{m+1}}$. 
	Let $X_{t,s}=(\bigwedge^t T\times \bigwedge^s S)/C_{p^n}$, graded by the types of $C_{p^n}$-orbits in $X_{t,s}/(\Sigma_{t}\times\Sigma_s)$ (see \Cref{notation:wedge_power_set} for the definition of $\bigwedge$ used here). 

	For an ordered $t$-element subset $I\subset T = C_{p^n}/C_{p^m}$ and an ordered $s$-element subset $J\subset S = C_{p^n}/C_{p^{m+1}}$, we define the \emph{incidence set} by
	\[ 
		\Inc(I,J) = \{ i\in C_{p^n}/C_{p^m} :  i\in I, i+C_{p^{m+1}} \in J\}. 
	\]
\end{definition}

\begin{construction}\label{Constr:General}
	Let $\bigkoszul_\bullet$ be the $(n+1)$-dimensional complex of free $\uR$-modules defined as follows. 
	For $i \in \mathbb{Z}$ and $0 \leq j \leq n$, let
	\[
		\bigkoszul_{i\vec{e}_j} := 
		\begin{cases}
			\uR\left\{
					(\bigwedge^i C_{p^n}/C_{p^j})/\Sigma_i
				\right\} \quad 
				& 
			\text{ if } 0 \leq i \leq p^{n-j}, 
				\\
			0 \quad 
				& 
			\text{ otherwise. }
		\end{cases}
	\]
	For $\vec{v} = \sum_{j=0}^n i_j \vec{e}_j$, let
	\[
		\bigkoszul_{\vec{v}} 
		:= \bigkoszul_{i_0 \vec{e}_0} 
				\boxtimes 
			\bigkoszul_{i_1 \vec{e}_1} 
				\boxtimes 
				\cdots 
				\boxtimes 
			\bigkoszul_{i_n \vec{e}_n}.
	\]
	When $\bigkoszul_{\vec{v}} \neq 0$, its generators have the form
	\[
		\left(\dfrac{N_e^{C_{p^\ell}}(z^{(\overline{I_0})}) \wedge N_{C_p}^{C_{p^\ell}}(u_1^{(\overline{I_1})}) \wedge \cdots \wedge R^{C_{p^n}}_{C_{p^\ell}}(u_n^{(I_n)})}{\nm_e^{C_{p^\ell}}(x^{(\overline{\Inc_0})})\nm_{C_p}^{C_{p^\ell}}(\nm_e^{C_p}(x)^{(\overline{\Inc_1})})\cdots \res^{C_{p^{n-1}}}_{C_{p^\ell}}(\nm_e^{C_{p^{n-1}}}(x)^{(\Inc_{n-1})})}\right)_{C_{p^n}/C_{p^\ell}}
	\]
	where $|I_j| = i_j$, $\Inc_j := \Inc(I_j,I_{j+1})$, $\ell = \min_{0 \leq j \leq n} \{ k : C_{p^k} = \stab(I_j)\}$ and for a subset $J\subset C_{p^n}/C_{p^k}$ with $k< \ell$ and stabilizer at least $C_{p^\ell}$, $\overline{J}$ is a set of representatives of $J$ under the $C_{p^\ell}$-action. 
	Define the differential $d_{\vec{e}_j}: \bigkoszul_{\vec{v}} \to \bigkoszul_{\vec{v}-\vec{e}_j}$ by sending the generator above to 
	\[
		\sum_{k \in I_j} (-1)^k \res_{C_{p^\ell}}^{C_{p^j}}\nm_{e}^{C_{p^j}}(x^{(k)})  \left(\dfrac{N_e^{C_{p^\ell}}(z^{(\overline{I_0})}) \wedge \cdots \wedge R^{C_{p^j}}_{C_{p^\ell}}(u_j^{(I_j)\setminus k}) \wedge \cdots \wedge R^{C_{p^n}}_{C_{p^\ell}}(u_n^{(I_n)})}{\nm_e^{C_{p^\ell}}(x^{(\overline{\Inc_0})})\cdots \res^{C_{p^{n-1}}}_{C_{p^\ell}}(\nm_e^{C_{p^{n-1}}}(x)^{(\Inc_{n-1})})}\right)
	\]
if $j \geq \ell$, and
	\[
		\sum_{k \in \overline{I_j}} (-1)^k \tr_{C_{p^j}}^{C_{p^{\ell}}} \left( \nm_{e}^{C_{p^j}}(x^{(k)})  \dfrac{N_e^{C_{p^j}}(z^{(\overline{I_0})}) \wedge \cdots \wedge u_j^{(I_j)\setminus k} \wedge \cdots \wedge R^{C_{p^n}}_{C_{p^j}}(u_n^{(I_n)})}{\nm_e^{C_{p^j}}(x^{(\overline{\Inc_0})})\cdots \res^{C_{p^{n-1}}}_{C_{p^j}}(\nm_e^{C_{p^{n-1}}}(x)^{(\Inc_{n-1})})}\right)
	\]
	if $j < \ell$. 
	Notice that $\overline{J}$ in the formal fraction now refers to the $C_{p^j}$-action, while in the indexing set for the sum it refers to the $C_{p^\ell}$-action. 
	In both cases, we reduce as much as possible by cancelling common factors in the coefficient and denominator. 
	Since this differential is essentially a Koszul differential, $\bigkoszul_\bullet$ is indeed an $(n+1)$-dimensional complex. 

	Let $\overline{\bigkoszul}_\bullet$ denote the totalization of $\bigkoszul_\bullet$.
\end{construction}

Our goal in this section is to prove the following:

\begin{theorem}
	\label{Thm:Koszulpn}
	Let $p$ be an odd prime, $G = C_{p^n}$, and let $\uR = \uA^{\cO^{\top}}[x_{G/e}]$ be the free Tambara functor on an underlying generator. 
	The complex of $\uR$-modules $\overline{\bigkoszul}_\bullet$ defined in \cref{Constr:General} defines a free $\uR$-module resolution of the Burnside Mackey functor $\overline{\bigkoszul}_\bullet \to \uA$, where the map $\overline{\bigkoszul}_0 \to \uA$ is the quotient map sending $z \mapsto 0$. 
	The length of this resolution is $\sum_{i=0}^n p^i$. 
\end{theorem}

To prove $\overline{\bigkoszul}_\bullet \to \uA$ is a resolution, we will show that iteratively taking homology with respect to $d_{\vec{e}_0}$ through $d_{\vec{e}_k}$ kills homology up to level $C_{p^k}$. 
This will imply that the $E_2$-page of a spectral sequence converging to $H_*(\overline{\bigkoszul}_\bullet \to \uA)$ is zero, and thus we have a resolution. 

We begin by computing the ``horizontal homology" of $\bigkoszul_\bullet$, i.e.,  homology with respect to the differential $d_{\vec{e}_0}: \bigkoszul_{\vec{v}} \to \bigkoszul_{\vec{v}-\vec{e}_0}$. 

\begin{proposition}
	Let $H_{\vec{v}}^0 := H_{\vec{v}}(\bigkoszul_{\bullet},d_{\vec{e}_0})$.

	\begin{enumerate}[(a)]

		\item We have
			\[
				H^0_{\vec{0}} 
					\cong \uR/\langle x \rangle^{\bot},
			\]
		where $\langle x \rangle^\bot$ is the Green ideal generated by $x \in \uR(C_{p^n}/e)$. 

		\item If $\vec{v}_0 = 0$ with $\vec{v} \neq 0$, then
			\[
				H^0_{\vec{v}} \cong 
				\bigoplus_{\ell=1}^n
				\bigoplus_{\stab(I_1,\ldots,I_n)=C_{p^\ell}}
				\uR\{ U_{I_1, \ldots, I_n} \}/
				\langle \res^{C_{p}}_{e}(U_{I_1, \ldots, I_n}) \rangle,
			\]
			where we abbreviate 
			\[ 
				U_{I_1, \ldots, I_n} = 
				\frac{N_{C_p}^{C_{p^\ell}}(u_1^{(\overline{I_1})})\wedge \cdots \wedge R_{C_{p^\ell}}^{C_{p^n}}(u_n^{(I_n)})}{\nm_e^{C_{p^\ell}}(x)^{(\overline{\Inc_1})}\cdots \res^{C_{p^{n-1}}}_{C_{p^\ell}}\nm_e^{C_{p^{n-1}}}(x)^{(\Inc_{n-1})}} 
			\]

		\item If $\vec{v}_0 \neq 0$, then after restriction to $C_{p^\ell}$:
			\[
				H^0_{\vec{v}} \cong 
				\bigoplus_{\ell=1}^n
				\bigoplus_{ \stab(I_0,\ldots,I_n) 
					= C_{p^{\ell}}} 
					\ker(\res^{C_{p^{\ell}}}_e )
					\cdot U_{I_0, \ldots, I_n}.
			\]
	\end{enumerate}

	In particular, this homology vanishes at level $C_{p^n}/e$ in all three cases. 
\end{proposition}

We note that the last case is often zero, e.g., if $\bigkoszul_{\vec{v}} = 0$ or if $\vec{v}_0$ is not divisible by $p$. 

\begin{proof}
	Since $(\bigkoszul_\bullet, d_{\vec{e}_0})$ is a Koszul complex at the underlying level by construction, the only nontrivial homology will appear at levels $C_{p^n}/C_{p^\ell}$ for $\ell > 0$, and it will be nontrivial precisely when we have a generator in $\bigkoszul_\bullet$ with $\stab(I_0) \neq e$.
	The resulting homology can then be computed from the definition of $\bigkoszul_\bullet$ using \cref{lemma:kernels} (cf. the analogous calculations for $n=1$ and $n=2$ in the previous sections). 
\end{proof}

In other words, the horizontal homology is supported strictly above the underlying level, but will be nontrivial at levels $C_{p^n}/C_p$, $C_{p^n}/C_{p^2}$, and so on, all the way up to level $C_{p^n}/C_{p^n}$. We now wish to iteratively take homology. If
\[
	H^0_{\vec{v}} := H_{\vec{v}}(\bigkoszul_{\bullet}, d_{\vec{e}_0})
\]
denotes the horizontal homology groups just computed, we define
\[
	H^k_{\vec{v}} := H_{\vec{v}}(H^{k-1}_{\bullet}, d_{\vec{e}_k})
\]
for $1 \leq k \leq n$. The following proposition says that $H^k_\bullet$ is supported in levels strictly higher than $C_{p^n}/C_{p^k}$. 

\begin{proposition}
	For all $\vec{v} \in \z^{n+1}$ and $0 \leq j \leq k$, 
	\[
		H^k_{\vec{v}}(C_{p^n}/C_{p^j}) = 0.
	\]
\end{proposition}

\begin{proof}
	We proceed by induction on $k$. 
	The base case $k=0$ is handled by the previous proposition. 

	We may suppose inductively that $H^{k-1}_{\vec{v}}(C_{p^n}/C_{p^j}) = 0$ for all $\vec{v} \in \z^{n+1}$ and $0 \leq j \leq k-1$, so we only need to show that $H^k_{\vec{v}}(C_{p^n}/C_{p^k}) = 0$. 
	But by construction and analogous calculations to before, the complex $(H^{k-1}_{\vec{v}}, d_{\vec{e}_k})$ is a Koszul complex at level $C_{p^n}/C_{p^k}$, so $H^k_{\vec{v}}(C_{p^n}/C_{p^k})=0$.  
\end{proof}

\begin{proof}[{Proof of \cref{Thm:Koszulpn}}]
	Filtering $\bigkoszul_\bullet$ by 
	\[
		F_j \bigkoszul_\bullet 
		:= \bigoplus_{\vec{v}: \vec{v}_n < j}
			\bigkoszul_{\vec{v}}
	\]
	gives rise to a spectral sequence with
	\begin{equation}
		\label{Eqn:E2a}
			E^2_{s,t} = 
				H_s( H_t( \Tot (\bigkoszul|_{\vec{v}_n = s}))) \Rightarrow H_{s+t} \overline{\bigkoszul}.
	\end{equation}
	For each $s \in \mathbb{Z}$, we may filter $\bigkoszul|_{\vec{v}_n =s}$ by
	\[
		F_j \bigkoszul|_{\vec{v}_n=s} 
			:= \bigoplus_{\vec{v} : \vec{v}_n=s, \vec{v}_{n-1} < j} \bigkoszul_{\vec{v}}
	\]
	to obtain an analogous spectral sequence
	\[
		E^2_{i,j} = H_i( H_j( \Tot (\bigkoszul|_{\vec{v}_n=s, \vec{v}_{n-1}=j})) \Rightarrow H_{i+j}(\Tot (\bigkoszul|_{\vec{v}_n = s})).
	\]
	Iterating this, we obtain a sequence of $n+1$ spectral sequences which allows us to identify the $E^2$-page in \eqref{Eqn:E2a} with $H^n_\bullet$. 
	Since $H^n_\bullet = 0$ by the previous proposition, the result follows. 
\end{proof}

\begin{remark}
	To simplify notation, we have focused on the free Tambara functors on underlying generators in this section. 
	More generally, the $(n+1)$-dimensional complex defined in \cref{Constr:General} may be modified to produce a resolution of $\uA$ by free $\uA^{\cO}[x_{C_{p^n}/C_{p^m}}]$-modules, where $\cO$ is an indexing category such that the underlying Mackey functor of $\uA^{\cO}[x_{C_{p^n}/C_{p^m}}]$ is free (cf. \cite{HMQ21}). For $i \in \z$ and $m \leq j \leq n$, one defines 
	\[
		\bigkoszul_{i \vec{e}_j} := 
		\begin{cases}
			\uR \left\{ (\bigwedge^i C_{p^n}/C_{p^j})/\Sigma_i\right\} \quad 
				& 
			\text{ if } 0 \leq i \leq p^{n-j}, 
				\\
			0 \quad 
				& 
			\text{ otherwise.}
		\end{cases}
	\]
	This produces an $(n-m+1)$-dimensional complex with Koszul-like differentials $d_{\vec{e}_j}$, $m \leq j \leq n$, defined as in \cref{Constr:General}. 
	The ``horizontal homology" becomes homology with respect to $d_{\vec{e}_m}$ and vanishes in levels $C_{p^n}/C_{p^j}$ for $0 \leq j \leq m$, and as above, each additional direction we take homology kills another level. 
	An analogous spectral sequence argument then implies the total complex is a resolution of $\uA$. 
\end{remark}

\appendix

\section{Figures}

	\begin{landscape}
	\begin{figure}	
		\scalebox{0.8}{
		\begin{tikzcd}[ampersand replacement=\&, column sep=small]
			\uR\{\sfrac{C_9}{C_9}\}
			\ar[r]
				\&
			\uR\{\sfrac{C_9}{e}\}
			\ar[r]
				\&
			\uR\{4(\sfrac{C_9}{e})\}
			\ar[r]
				\&
			\uR\{9 (\sfrac{C_9}{e}) + \sfrac{C_9}{C_3}\}
			\ar[r]
				\&
			\uR\{14 (\sfrac{C_9}{e})\}
			\ar[r]
				\&
			\uR\{14 (\sfrac{C_9}{e})\}
			\ar[r]
				\&
			\uR\{9 (\sfrac{C_9}{e}) + \sfrac{C_9}{C_3}\}
			\ar[r]
				\&
			\uR\{4(\sfrac{C_9}{e})\}
			\ar[r]
				\&
			\uR\{\sfrac{C_9}{e}\}
			\ar[r]
				\&
			\uR\{\sfrac{C_9}{C_9}\}
			\ar[ddddd, bend left=50, "\nm_e^{C_9}(x^{(0)})" description]
				\\
			\uR\{\sfrac{C_9}{C_3}\}
			\ar[r]
			\ar[u]
				\&
			\uR\{3(\sfrac{C_9}{e})\}
			\ar[r]
			\ar[u]
				\&
			\uR\{12(\sfrac{C_9}{e})\}
			\ar[r]
			\ar[u]
				\&
			\uR\{27 (\sfrac{C_9}{e}) + 3(\sfrac{C_9}{C_3})\}
			\ar[r]
			\ar[u]
				\&
			\uR\{42 (\sfrac{C_9}{e})\}
			\ar[r]
			\ar[u]
				\&
			\uR\{42 (\sfrac{C_9}{e})\}
			\ar[r]
			\ar[u]
				\&
			\uR\{27 (\sfrac{C_9}{e}) + 3(\sfrac{C_9}{C_3})\}
			\ar[r]
			\ar[u]
				\&
			\uR\{12(\sfrac{C_9}{e})\}
			\ar[r]
			\ar[u]
				\&
			\uR\{3(\sfrac{C_9}{e})\}
			\ar[r]
			\ar[u]
				\&
			\uR\{\sfrac{C_9}{C_3}\}
			\ar[u]
				\\
			\uR\{\sfrac{C_9}{C_3}\}
			\ar[r]
			\ar[u]
				\&
			\uR\{3(\sfrac{C_9}{e})\}
			\ar[r]
			\ar[u]
				\&
			\uR\{12(\sfrac{C_9}{e})\}
			\ar[r]
			\ar[u]
				\&
			\uR\{27 (\sfrac{C_9}{e}) + 3(\sfrac{C_9}{C_3})\}
			\ar[r]
			\ar[u]
				\&
			\uR\{42 (\sfrac{C_9}{e})\}
			\ar[r]
			\ar[u]
				\&
			\uR\{42 (\sfrac{C_9}{e})\}
			\ar[r]
			\ar[u]
				\&
			\uR\{27 (\sfrac{C_9}{e}) + 3(\sfrac{C_9}{C_3})\}
			\ar[r]
			\ar[u]
				\&
			\uR\{12(\sfrac{C_9}{e})\}
			\ar[r]
			\ar[u]
				\&
			\uR\{3(\sfrac{C_9}{e})\}
			\ar[r]
			\ar[u]
				\&
			\uR\{\sfrac{C_9}{C_3}\}
			\ar[u]
				\\
			\uR\{\sfrac{C_9}{C_9}\}
			\ar[r]
			\ar[u]
				\&
			\uR\{\sfrac{C_9}{e}\}
			\ar[r]
			\ar[u]
				\&
			\uR\{4(\sfrac{C_9}{e})\}
			\ar[r]
			\ar[u]
				\&
			\uR\{9 (\sfrac{C_9}{e}) + \sfrac{C_9}{C_3}\}
			\ar[r]
			\ar[u]
				\&
			\uR\{14 (\sfrac{C_9}{e})\}
			\ar[r]
			\ar[u]
				\&
			\uR\{14 (\sfrac{C_9}{e})\}
			\ar[r]
			\ar[u]
				\&
			\uR\{9 (\sfrac{C_9}{e}) + \sfrac{C_9}{C_3}\}
			\ar[r]
			\ar[u]
				\&
			\uR\{4(\sfrac{C_9}{e})\}
			\ar[r]
			\ar[u]
				\&
			\uR\{\sfrac{C_9}{e}\}
			\ar[r]
			\ar[u]
				\&
			\uR\{\sfrac{C_9}{C_9}\}
			\ar[u]
				\\
				\\
			\uR\{\sfrac{C_9}{C_9}\}
			\ar[r]
				\&
			\uR\{\sfrac{C_9}{e}\}
			\ar[r]
				\&
			\uR\{4(\sfrac{C_9}{e})\}
			\ar[r]
				\&
			\uR\{9 (\sfrac{C_9}{e}) + \sfrac{C_9}{C_3}\}
			\ar[r]
				\&
			\uR\{14 (\sfrac{C_9}{e})\}
			\ar[r]
				\&
			\uR\{14 (\sfrac{C_9}{e})\}
			\ar[r]
				\&
			\uR\{9 (\sfrac{C_9}{e}) + \sfrac{C_9}{C_3}\}
			\ar[r]
				\&
			\uR\{4(\sfrac{C_9}{e})\}
			\ar[r]
				\&
			\uR\{\sfrac{C_9}{e}\}
			\ar[r]
				\&
			\uR\{\sfrac{C_9}{C_9}\}
				\\
			\uR\{\sfrac{C_9}{C_3}\}
			\ar[r]
			\ar[u]
				\&
			\uR\{3(\sfrac{C_9}{e})\}
			\ar[r]
			\ar[u]
				\&
			\uR\{12(\sfrac{C_9}{e})\}
			\ar[r]
			\ar[u]
				\&
			\uR\{27 (\sfrac{C_9}{e}) + 3(\sfrac{C_9}{C_3})\}
			\ar[r]
			\ar[u]
				\&
			\uR\{42 (\sfrac{C_9}{e})\}
			\ar[r]
			\ar[u]
				\&
			\uR\{42 (\sfrac{C_9}{e})\}
			\ar[r]
			\ar[u]
				\&
			\uR\{27 (\sfrac{C_9}{e}) + 3(\sfrac{C_9}{C_3})\}
			\ar[r]
			\ar[u]
				\&
			\uR\{12(\sfrac{C_9}{e})\}
			\ar[r]
			\ar[u]
				\&
			\uR\{3(\sfrac{C_9}{e})\}
			\ar[r]
			\ar[u]
				\&
			\uR\{\sfrac{C_9}{C_3}\}
			\ar[u, "\nm_e^{C_3}(x^{(0)})"]
				\\
			\uR\{\sfrac{C_9}{C_3}\}
			\ar[r]
			\ar[u]
				\&
			\uR\{3(\sfrac{C_9}{e})\}
			\ar[r]
			\ar[u]
				\&
			\uR\{12(\sfrac{C_9}{e})\}
			\ar[r]
			\ar[u]
				\&
			\uR\{27 (\sfrac{C_9}{e}) + 3(\sfrac{C_9}{C_3})\}
			\ar[r]
			\ar[u]
				\&
			\uR\{42 (\sfrac{C_9}{e})\}
			\ar[r]
			\ar[u]
				\&
			\uR\{42 (\sfrac{C_9}{e})\}
			\ar[r]
			\ar[u]
				\&
			\uR\{27 (\sfrac{C_9}{e}) + 3(\sfrac{C_9}{C_3})\}
			\ar[r]
			\ar[u]
				\&
			\uR\{12(\sfrac{C_9}{e})\}
			\ar[r]
			\ar[u]
				\&
			\uR\{3(\sfrac{C_9}{e})\}
			\ar[r]
			\ar[u]
				\&
			\uR\{\sfrac{C_9}{C_3}\}
			\ar[u]
				\\
			\uR\{\sfrac{C_9}{C_9}\}
			\ar[r]
			\ar[u]
				\&
			\uR\{\sfrac{C_9}{e}\}
			\ar[r]
			\ar[u]
				\&
			\uR\{4(\sfrac{C_9}{e})\}
			\ar[r]
			\ar[u]
				\&
			\uR\{9 (\sfrac{C_9}{e}) + \sfrac{C_9}{C_3}\}
			\ar[r]
			\ar[u]
				\&
			\uR\{14 (\sfrac{C_9}{e})\}
			\ar[r]
			\ar[u]
				\&
			\uR\{14 (\sfrac{C_9}{e})\}
			\ar[r]
			\ar[u]
				\&
			\uR\{9 (\sfrac{C_9}{e}) + \sfrac{C_9}{C_3}\}
			\ar[r]
			\ar[u]
				\&
			\uR\{4(\sfrac{C_9}{e})\}
			\ar[r]
			\ar[u]
				\&
			\uR\{\sfrac{C_9}{e}\}
			\ar[r]
			\ar[u]
				\&
			\uR\{\sfrac{C_9}{C_9}\}
			\ar[u]
			\\
			\\
		\end{tikzcd}
		}
		\caption{The tricomplex $\bigkoszul_{r, s, t}$ for $C_9$, without names of generators. The first row in the second block is the complex $\bigkoszul_{\bullet,0,0} = \uK_\bullet$ which lifts the classical Koszul complex. The index $r$ describes the column, $s$ describes the row, $t = 0$ is the bottom $10 \times 4$-block, and $t = 1$ is the top $10 \times 4$-block. Differentials from $t = 1$ to $t = 0$ are omitted except between degrees $(0,0,1)$ and $(0,0,0)$. 
		}\label{C9tricomplex}
	\end{figure}
	\thispagestyle{empty}
	\end{landscape}
	
	\newpage
	\thispagestyle{empty}
	\pagestyle{empty}
	\begin{figure}
		\hspace*{-3cm}
		\begin{tikzcd}[ampersand replacement=\&]
			\ker(\res^{C_9}_e)
				\&
			0
				\&
			0
				\&
			\ker(\res^{C_3}_e)
				\&
			0
				\&
			0
				\&
			\ker(\res^{C_3}_e)
				\&
			0
				\&
			0
				\&
			\uR\{\sfrac{C_9}{C_9}\}/\res^{C_9}_e
				\ar[ddddd, bend left=60]
				\\
			\ker(\res^{C_3}_e)^{\oplus 3}
			\ar[u]
				\&
			0
			\ar[u]
				\&
			0
			\ar[u]
				\&
			\ker(\res^{C_3}_e)^{\oplus 3}
			\ar[u]
				\&
			0
			\ar[u]
				\&
			0
			\ar[u]
				\&
			\ker(\res^{C_3}_e)^{\oplus 3}
			\ar[u]
				\&
			0
			\ar[u]
				\&
			0
			\ar[u]
				\&
			\uR\{\sfrac{C_9}{C_3}\}/\res^{C_3}_e
			\ar[u]
				\\
			\ker(\res^{C_3}_e)^{\oplus 3}
			\ar[u]
				\&
			0
			\ar[u]
				\&
			0
			\ar[u]
				\&
			\ker(\res^{C_3}_e)^{\oplus 3}
			\ar[u]
				\&
			0
			\ar[u]
				\&
			0
			\ar[u]
				\&
			\ker(\res^{C_3}_e)^{\oplus 3}
			\ar[u]
				\&
			0
			\ar[u]
				\&
			0
			\ar[u]
				\&
			\uR\{\sfrac{C_9}{C_3}\}/\res^{C_3}_e
			\ar[u]
				\\
			\ker(\res^{C_9}_e)
			\ar[u]
				\&
			0
			\ar[u]
				\&
			0
			\ar[u]
				\&
			\ker(\res^{C_3}_e)
			\ar[u]
				\&
			0
			\ar[u]
				\&
			0
			\ar[u]
				\&
			\ker(\res^{C_3}_e)
			\ar[u]
				\&
			0
			\ar[u]
				\&
			0
			\ar[u]
				\&
			\uR\{\sfrac{C_9}{C_9}\}/\res^{C_9}_e
			\ar[u]
				\\
				\\
			\ker(\res^{C_9}_e)
				\&
			0 
				\&
			0
				\&
			\ker(\res^{C_3}_e)
				\&
			0
				\&
			0
				\&
			\ker(\res^{C_3}_e)
				\&
			0
				\&
			0
				\&
			\uR/\langle x \rangle^\perp			
				\\
			\ker(\res^{C_3}_e)^{\oplus 3}
			\ar[u]
				\&
			0
			\ar[u]
				\&
			0
			\ar[u]
				\&
			\ker(\res^{C_3}_e)^{\oplus 3}
			\ar[u]
				\&
			0
			\ar[u]
				\&
			0
			\ar[u]
				\&
			\ker(\res^{C_3}_e)^{\oplus 3}
			\ar[u]
				\&
			0
			\ar[u]
				\&
			0
			\ar[u]
				\&
			\uR\{\sfrac{C_9}{C_3}\} / \res^{C_3}_e
			\ar[u]
				\\
			\ker(\res^{C_3}_e)^{\oplus 3}
			\ar[u]
				\&
			0
			\ar[u]
				\&
			0
			\ar[u]
				\&
			\ker(\res^{C_3}_e)^{\oplus 3}
			\ar[u]
				\&
			0
			\ar[u]
				\&
			0
			\ar[u]
				\&
			\ker(\res^{C_3}_e)^{\oplus 3}
			\ar[u]
				\&
			0
			\ar[u]
				\&
			0
			\ar[u]
				\&
			\uR\{\sfrac{C_9}{C_3}\} /\res^{C_3}_e
			\ar[u]
				\\
			\ker(\res^{C_9}_e)
			\ar[u]
				\&
			0
			\ar[u]
				\&
			0
			\ar[u]
				\&
			\ker(\res^{C_3}_e)
			\ar[u]
				\&
			0
			\ar[u]
				\&
			0
			\ar[u]
				\&
			\ker(\res^{C_3}_e)
			\ar[u]
				\&
			0
			\ar[u]
				\&
			0
			\ar[u]
				\&
			\uR\{\sfrac{C_9}{C_9}\}/ \res^{C_9}_e
			\ar[u]
			\\
		\end{tikzcd}
		\caption{The horizontal homology of the tricomplex $\bigkoszul_{\bullet,\bullet,\bullet}$. Here, $\ker(\res^H_e)$ denotes the sub-$\uR$-module of $\uR\{C_9/H\}$ generated by the kernel of restriction (in $\uR$) from $H$ to $e$, and $\uR\{C_9/H\}/\res^{C_9}_e$ is the quotient of $\uR\{C_9/H\}$ by the sub-$\uR$-module generated by restriction of the generators to the underlying level. The notation $\uR / \langle x \rangle^\bot$ denotes the quotient of $\uR$ by the Green ideal generated by $x$. 
		}\label{C9TricomplexHorizontalHomology}
	\end{figure}
	\thispagestyle{empty}
	
	\newpage
	\thispagestyle{empty}
	\begin{figure}
		\hspace*{-3cm}
		\begin{tikzcd}[ampersand replacement=\&]
			\ker(\res^{C_9}_e)/\ker(\res^{C_9}_{C_3})
				\&
			0
				\&
			0
				\&
			0
				\&
			0
				\&
			0
				\&
			0
				\&
			0
				\&
			0
				\&
			\uR\{\sfrac{C_9}{C_9}\}/\res^{C_9}_e
				\ar[ddddd, bend left=50]
				\\
			0
				\&
			0
				\&
			0
				\&
			0
				\&
			0
				\&
			0
				\&
			0
				\&
			0
				\&
			0
				\&
			0
				\\
			0
				\&
			0
				\&
			0
				\&
			0
				\&
			0
				\&
			0
				\&
			0
				\&
			0
				\&
			0
				\&
			0
				\\
			\ker(\res^{C_9}_e)
				\&
			0
				\&
			0
				\&
			0
				\&
			0
				\&
			0
				\&
			0
				\&
			0
				\&
			0
				\&
			\ker(\res^{C_9}_{C_3})/\res^{C_9}_e
				\\
				\\
			\ker(\res^{C_9}_e)/\ker(\res^{C_9}_{C_3})
				\&
			0 
				\&
			0
				\&
			0
				\&
			0
				\&
			0
				\&
			0
				\&
			0
				\&
			0
				\&
			\uR/\langle x, \nm_e^{C_3}(x) \rangle^\perp			
				\\
			0
				\&
			0
				\&
			0
				\&
			0
				\&
			0
				\&
			0
				\&
			0
				\&
			0
				\&
			0
				\&
			0
				\\
			0
				\&
			0
				\&
			0
				\&
			0
				\&
			0
				\&
			0
				\&
			0
				\&
			0
				\&
			0
				\&
			0
				\\
			\ker(\res^{C_9}_e)
				\&
			0
				\&
			0
				\&
			0
				\&
			0
				\&
			0
				\&
			0
				\&
			0
				\&
			0
				\&
			\ker(\res^{C_9}_{C_3})/\res^{C_9}_e
		\end{tikzcd}
		\caption{The vertical homology of the horizontal homology of the tricomplex $\bigkoszul_{\bullet,\bullet,\bullet}$.}\label{C9TricomplexVerticalHomology}
	\end{figure}

	\newpage
	\thispagestyle{empty}
	
	\begin{landscape}
	\begin{figure}
	\hspace*{-3cm}
	\scalebox{0.75}{
		\begin{tikzcd}[ampersand replacement=\&]
			\ker(\res^{C_9}_e)
				\&
			0
				\&
			0
				\&
			\ker(\res^{C_3}_e)
				\&
			0
				\&
			0
				\&
			\ker(\res^{C_3}_e)
				\&
			0
				\&
			0
				\&
			\uA\{\sfrac{C_9}{C_9}\}/\res^{C_9}_e
				\ar[ddddd, bend left=60]
				\\
			\ker(\res^{C_3}_e)^{\oplus 3}
			\ar[u]
				\&
			0
			\ar[u]
				\&
			0
			\ar[u]
				\&
			\ker(\res^{C_3}_e)^{\oplus 3}
			\ar[u]
				\&
			0
			\ar[u]
				\&
			0
			\ar[u]
				\&
			\ker(\res^{C_3}_e)^{\oplus 3}
			\ar[u]
				\&
			0
			\ar[u]
				\&
			0
			\ar[u]
				\&
			\uA\{\sfrac{C_9}{C_3}\}/\res^{C_3}_e
			\ar[u]
				\\
			\ker(\res^{C_3}_e)^{\oplus 3}
			\ar[u]
				\&
			0
			\ar[u]
				\&
			0
			\ar[u]
				\&
			\ker(\res^{C_3}_e)^{\oplus 3}
			\ar[u]
				\&
			0
			\ar[u]
				\&
			0
			\ar[u]
				\&
			\ker(\res^{C_3}_e)^{\oplus 3}
			\ar[u]
				\&
			0
			\ar[u]
				\&
			0
			\ar[u]
				\&
			\uA\{\sfrac{C_9}{C_3}\}/\res^{C_3}_e
			\ar[u]
				\\
			\ker(\res^{C_9}_e)
			\ar[u]
				\&
			0
			\ar[u]
				\&
			0
			\ar[u]
				\&
			\ker(\res^{C_3}_e)
			\ar[u]
				\&
			0
			\ar[u]
				\&
			0
			\ar[u]
				\&
			\ker(\res^{C_3}_e)
			\ar[u]
				\&
			0
			\ar[u]
				\&
			0
			\ar[u]
				\&
			\uA\{\sfrac{C_9}{C_9}\}/\res^{C_9}_e
			\ar[u]
				\\
				\\
			\uA\{\sfrac{C_9}{C_9}\}
				\&
			\uA\{\sfrac{C_9}{e}\}
				\&
			\uA\{4(\sfrac{C_9}{e})\}
				\&
			\uA\{9(\sfrac{C_9}{e}) + \sfrac{C_9}{C_3}\}
				\&
			\uA\{14(\sfrac{C_9}{e})\}
				\&
			\uA\{14(\sfrac{C_9}{e})\}
				\&
			\uA\{9(\sfrac{C_9}{e}) + \sfrac{C_9}{C_3}\}
				\&
			\uA\{4(\sfrac{C_9}{e})\}
				\&
			\uA\{\sfrac{C_9}{e}\}
				\&
			\uA			
				\\
			\ker(\res^{C_3}_e)
			\ar[u]
				\&
			0
			\ar[u]
				\&
			0
			\ar[u]
				\&
			\ker(\res^{C_3}_e)^{\oplus 2} \oplus (\uA\{\sfrac{C_9}{C_3}\}/\res^{C_3}_e)
			\ar[u]
				\&
			0
			\ar[u]
				\&
			0
			\ar[u]
				\&
			\ker(\res^{C_3}_e) \oplus (\uA\{\sfrac{C_9}{C_3}\} /\res^{C_3}_e)^{\oplus 2}
			\ar[u]
				\&
			0
			\ar[u]
				\&
			0
			\ar[u]
				\&
			\uA\{\sfrac{C_9}{C_3}\} / \res^{C_3}_e
			\ar[u]
				\\
			\ker(\res^{C_3}_e)
			\ar[u]
				\&
			0
			\ar[u]
				\&
			0
			\ar[u]
				\&
			\ker(\res^{C_3}_e)^{\oplus 3}
			\ar[u]
				\&
			0
			\ar[u]
				\&
			0
			\ar[u]
				\&
			\ker(\res^{C_3}_e)^{\oplus 2} \oplus \uA\{\sfrac{C_9}{C_3}\} /\res^{C_3}_e
			\ar[u]
				\&
			0
			\ar[u]
				\&
			0
			\ar[u]
				\&
			\uA\{\sfrac{C_9}{C_3}\} /\res^{C_3}_e
			\ar[u]
				\\
			\ker(\res^{C_9}_e)
			\ar[u]
				\&
			0
			\ar[u]
				\&
			0
			\ar[u]
				\&
			\ker(\res^{C_3}_e)
			\ar[u]
				\&
			0
			\ar[u]
				\&
			0
			\ar[u]
				\&
			\ker(\res^{C_3}_e)
			\ar[u]
				\&
			0
			\ar[u]
				\&
			0
			\ar[u]
				\&
			\uA\{\sfrac{C_9}{C_9}\}/ \res^{C_9}_e
			\ar[u]
		\end{tikzcd}
		}
		\caption{The horizontal homology of $\bigkoszul_{\bullet}\boxtimes_{\uR}\uA$.}\label{TorHorizontalHomology}
	\end{figure}
	
	\thispagestyle{empty}
	\end{landscape}

	\newpage
	\thispagestyle{empty}
	\begin{landscape}
	\begin{figure}
	\hspace*{-3cm}
	\scalebox{0.75}{
		\begin{tikzcd}[ampersand replacement=\&]
			\ker(\res^{C_9}_e)/\res^{C_9}_{C_3}
				\&
			0
				\&
			0
				\&
			0
				\&
			0
				\&
			0
				\&
			0
				\&
			0
				\&
			0
				\&
			\uA\{\sfrac{C_9}{C_9}\}/\res^{C_9}_{C_3}
				\ar[ddddd, bend left=60]
				\\
			0
				\&
			0
				\&
			0
				\&
			0
				\&
			0
				\&
			0
				\&
			0
				\&
			0
				\&
			0
				\&
			0
				\\
			0
				\&
			0
				\&
			0
				\&
			0
				\&
			0
				\&
			0
				\&
			0
				\&
			0
				\&
			0
				\&
			0
				\\
			\ker(\res^{C_9}_e)
				\&
			0
				\&
			0
				\&
			0
				\&
			0
				\&
			0
				\&
			0
				\&
			0
				\&
			0
				\&
			\ker(\res^{C_9}_{C_3})
				\\
				\\
			\uA\{\sfrac{C_9}{C_9}\} / \langle \res^{C_9}_e \cap \ker(\res^{C_3}_e)\rangle
				\&
			\uA\{\sfrac{C_9}{e}\}
				\&
			\uA\{4(\sfrac{C_9}{e})\}
				\&
			\uA\{9(\sfrac{C_9}{e})\} \oplus \ubbZ\{\sfrac{C_9}{C_3}\}
				\&
			\uA\{14(\sfrac{C_9}{e})\}
				\&
			\uA\{14(\sfrac{C_9}{e})\}
				\&
			\uA\{9(\sfrac{C_9}{e})\} \oplus \ubbZ\{\sfrac{C_9}{C_3}\}
				\&
			\uA\{4(\sfrac{C_9}{e})\}
				\&
			\uA\{\sfrac{C_9}{e}\}
				\&
			\uA			
				\\
			0
				\&
			0
				\&
			0
				\&
			\uL
				\&
			0
				\&
			0
				\&
			\uL^{\oplus 2}
				\&
			0
				\&
			0
				\&
			\uA\{\sfrac{C_9}{C_3}\} / \res^{C_3}_e
				\\
			0
				\&
			0
				\&
			0
				\&
			0
				\&
			0
				\&
			0
				\&
			\uL
				\&
			0
				\&
			0
				\&
			\uA\{\sfrac{C_9}{C_3}\} /\res^{C_3}_e
				\\
			\ker(\res^{C_9}_{e})
				\&
			0
				\&
			0
				\&
			0
				\&
			0
				\&
			0
				\&
			0
				\&
			0
				\&
			0
				\&
			\uA\{\sfrac{C_9}{C_9}\}/ \res^{C_9}_e
		\end{tikzcd}
		}
		\caption{The vertical homology of the horizontal homology of $\bigkoszul_{\bullet}\boxtimes_{\uR}\uA$.}\label{TorVerticalHomology}
\end{figure}

	\thispagestyle{empty}
	\end{landscape}

\clearpage
\pagestyle{plain}
\bibliographystyle{alpha}
\bibliography{references}

\end{document}

%% file: preamble.tex
\usepackage[margin=1.5in]{geometry} 
\usepackage{amsmath, amsfonts, amsthm, amstext, amssymb} 
\usepackage{hyperref} 
\usepackage{todonotes} 
\usepackage[nameinlink, capitalize]{cleveref} 
\usepackage{tikz-cd} 
\usepackage[shortlabels]{enumitem} 
\usepackage{array}
\usepackage{scalerel}
\usepackage{mathtools}
\usepackage{cancel} 
\usepackage{lscape} 

\numberwithin{equation}{section} 
\numberwithin{table}{section} 

\theoremstyle{definition}
\newtheorem{definition}[equation]{Definition}
\newtheorem{remark}[equation]{Remark}
\newtheorem{example}[equation]{Example}
\newtheorem{notation}[equation]{Notation}
\newtheorem{construction}[equation]{Construction}
\theoremstyle{plain}
\newtheorem{lemma}[equation]{Lemma}
\newtheorem{proposition}[equation]{Proposition}
\newtheorem{theorem}[equation]{Theorem}
\newtheorem{corollary}[equation]{Corollary}

\newtheorem{maintheorem}{Theorem} 
\newtheorem{maincorollary}[maintheorem]{Corollary}

\newcommand{\mhyphen}{\textnormal{-}}

\newcommand{\xto}[1]{\xrightarrow{#1}}
\newcommand{\from}{\leftarrow}
\newcommand{\xfrom}[1]{\xleftarrow{#1}}

\def\bbZ{\mathbb{Z}}
\def\z{\mathbb{Z}}

\def\cA{\mathcal{A}}

\def\cD{\mathcal{D}}

\def\cO{\mathcal{O}}
\def\cP{\mathcal{P}}

\def\uA{\underline{A}}

\def\uF{\underline{F}}

\def\uI{\underline{I}}
\def\uJ{\underline{J}}
\def\uK{\underline{K}}

\def\uL{\underline{L}}
\def\uM{\underline{M}}
\def\uN{\underline{N}}

\def\uQ{\underline{Q}}
\def\uR{\underline{R}}
\def\uS{\underline{S}}
\def\uT{\underline{T}}

\def\ubbZ{\underline{\mathbb{Z}}}
\def\uTor{\underline{\operatorname{Tor}}}
\def\uHH{\underline{\operatorname{HH}}}

\DeclareMathOperator{\Mod}{\mathcal{M}od}

\DeclareMathOperator{\Set}{\mathcal{S}et}
\DeclareMathOperator{\Mack}{\mathcal{M}ack}
\DeclareMathOperator{\Tamb}{\mathcal{T}amb}

\DeclareMathOperator{\Iso}{Iso}

\DeclareMathOperator{\Inf}{Inf}

\DeclareMathOperator{\HH}{HH}
\DeclareMathOperator{\Cone}{Cone}

\newcommand{\gen}{{\textnormal{gen}}}

\DeclareMathOperator{\id}{id} 
\DeclareMathOperator{\im}{im} 
\DeclareMathOperator{\Tor}{Tor}

\DeclareMathOperator{\Hom}{Hom}
\DeclareMathOperator{\Der}{Der}

\DeclareMathOperator{\coker}{coker} 
\DeclareMathOperator{\coeq}{coeq}
\DeclareMathOperator{\res}{res}
\DeclareMathOperator{\tr}{tr}
\DeclareMathOperator{\nm}{nm}

\DeclareMathOperator{\Inc}{Inc}
\DeclareMathOperator{\sgn}{sgn}
\DeclareMathOperator{\Tot}{Tot}

\DeclareMathOperator{\pr}{pr}
\DeclareMathOperator{\eval}{eval}
\DeclareMathOperator{\stab}{stab}

\newcommand{\bigkoszul}{\underline{\mathcal K}}
\newcommand{\ndivkoszul}{\underline{K}^\div}
\newcommand{\ndivpartial}{\partial^\div}

\makeatletter
\newcommand*{\boxwedge}{\mathbin{\mathpalette\@boxwedge{}}}
\newcommand*{\@boxwedge}[2]{%
  \sbox0{$#1\boxplus\m@th$}%
  \dimen2=.5\dimexpr\wd0-\ht0-\dp0\relax 
  \dimen@=\dimexpr\ht0+\dp0\relax
  \def\lw{.07}
  \kern\dimen2 
  \tikz[
    line width=\lw\dimen@,
    line join=round,
    x=\dimen@,
    y=\dimen@,
  ]
  \draw
    (\lw/2,0) rectangle (1-\lw,1-\lw)
    (\lw,0) -- (.5,1-\lw-\lw/2) -- (1-\lw-\lw/2 ,0)
  ;%
  \kern\dimen2 
}
\makeatother

\newcommand{\sfrac}[2]{{}^{#1}\!/_{\!#2}}
